\documentclass[12pt]{article}
\usepackage{amssymb}
\usepackage{mathrsfs}%Óиü¶àµÄÊýѧ·ûºÅ£¬±ÈÈ绨Ìå$\mathscr P$
\usepackage{pst-poly}     % From pstricks/contrib/pst-poly

\newtheorem{thm}{Theorem}[subsection]
\newtheorem{lem}[thm]{Lemma}
\newtheorem{cor}[thm]{Corollary}
\newtheorem{prop}[thm]{Proposition}
\newtheorem{con}[thm]{Conjecture}
\newtheorem{exa}[thm]{Example}

\newtheorem{sta}[thm]{Statement}
\newtheorem{obs}[thm]{Observation}
\newtheorem{prob}[thm]{Problem}
\newtheorem{que}[thm]{Question}
\newenvironment{pf}[1][Proof]{\noindent\textbf{#1.} }{\hfill\rule{1mm}{2mm}}

\makeatletter \@addtoreset{equation}{section} \makeatother

\begin{document}

\title{{\bf On Bondage Numbers of Graphs}\\
-- a survey with some comments\\
{\small(a lecture for my graduate students in seminar)}}
%\thanks {The work was supported  by NNSF of China (No.10671191).}}
\author
{Jun-Ming Xu \\ \\
{\small School of Mathematical Sciences}  \\
{\small University of Science and Technology of China}\\
{\small Hefei, Anhui, 230026, China}\\
{\small Email: xujm@ustc.edu.cn}
 }
\date{April 18, 2010}
 \maketitle

\begin{center}
\begin{minipage}{140mm}
%\vskip10pt
\begin{center} {\bf Abstract} \end{center}

The bondage number of a nonempty graph $G$ is the cardinality of a
smallest edge set whose removal from $G$ results in a graph with
domination number greater than the domination number of $G$. This
lecture gives a survey on the bondage number, including the known
results, problems and conjectures. We also summarize other types of
bondage numbers.

\end{minipage}
\end{center}

%\vskip6pt\noindent{\bf Keywords}: bondage number; domination number;
%crossing numbers; planar graphs

\vskip30pt
%\newpage

\tableofcontents
%\thispagestyle{empty}
%\markboth{Contents}{Combinatorial Network Theory}

\newpage

\section{Introduction}

For terminology and notation on graph theory not given here, the
reader is referred to~\cite{x03}. Let $G=(V,E)$ be a finite,
undirected and simple graph. We call $|V|$ and $|E|$ the order and
size of $G$, and denote them by $\upsilon=\upsilon(G)$ and
$\varepsilon=\varepsilon(G)$, respectively.

For a vertex $x$ in $G$, let $N_G(x)$ be the ({\it open}) {\it set
of neighbors} of $x$ and $N_G[x]=N[x]=N_G(x)\cup \{x\}$ be the {\it
closed set of neighbors} of $x$. For a subset $X\subset V(G)$,
$N_G(X)=(\cup_{x\in X} N_G(x))\cap \overline X$, $N_G[X]=N_G(X)\cup
X$, where $\overline X=V(G)\setminus X$. Let $E_x$ be the set of
edges incident with $x$ in $G$, that is, $E_x=\{xy\in E(G):\ y\in
N_G(x)\}$. We denote the degree of $x$ by $d_G(x)=|E_x|$. The
maximum and the minimum degree of $G$ are denoted by $\Delta(G)$ and
$\delta(G)$, respectively. The vertex of degree one is called a {\it
pendent vertex}, and the edge incident with a pendant vertex is
called a {\it pendant edge}.
%, and the distance between the vertices $x$ and $y$ by $d_G(x,y)$.

The bondage number is an important parameter of graphs which is
based upon the well-known domination number.

A subset $S\subseteq V(G)$ is called a {\it dominating set} of $G$
if $N[S]=V(G)$, i.e. every vertex $x$ in $\overline S$ has at least
one neighbor in $S$. The {\it domination number} of a graph $G$,
denoted by $\gamma(G)$, is the minimum cardinality among all
dominating sets, i.e.
 $$
 \gamma(G)=\min\{|S|:\ S\subseteq V(G), N[S]=V(G)\}.
 $$
A dominating set $S$ is called a $\gamma$-set of $G$ if
$|S|=\gamma(G)$.

The domination is so an important and classic conception that it has
become one of the most widely studied topics in graph theory, and
also is frequently studied property of networks. The domination,
with its many variations, is now well studied in graph and networks
theory. A thorough study of domination appears in \cite{hhs97a,
hhs97b}. However, the problem determining domination number for
general graphs was early proved to be NP-complete (see GT2 in
Appendix in Garey and Johnson~\cite{gj79}, 1979).
%
%In 1998, the book by Dunbar {\it et al.}\,\cite{dhtv98} renders a
%bibliography with more than 1200 papers on the domination
%number.{\bf \color{red}???}%´ËÎÄÏ×´ý²é£¡£¡.

Among various problems related with the domination number, some
focus on graph alterations and their effects on the domination
number. Here we are concerned with a particular graph alternation,
the removal of edges from a graph.

Graphs with domination numbers changed upon removal of an edge were
first investigated by Walikar and Acharya \cite{wa79} in 1979. A
graph is called {\it edge domination-critical graph} if
$\gamma(G-e)>\gamma(G)$ for every edge $e\in E(G)$. The edge
domination-critical graph was were characterized by Bauer {\it et
al.}~\cite{bhns83} in 1983, that is, a graph is edge
domination-critical if and only if it is the union of stars. The
proof is simple. The sufficiency is clear. Suppose that $S$ is a
$\gamma$-set of $G$. Then every vertex of degree at least two must
be in $S$, and no two vertices in $S$ can be adjacent. Hence $G$ is
a union of stars.

However, for lots of graphs, the domination number is out of the
range of one-edge removal. It is immediate that
$\gamma(H)\geqslant\gamma(G)$ for any spanning subgraph $H$ of $G$.
Every graph $G$ has a spanning forest $T$ with $\gamma(G)=\gamma(T)$
and so, in general, a graph will have a nonempty set of edges
$F\subseteq E(G)$ for which $\gamma(G-F)=\gamma(G)$.
%Such a set $F$ will be called an {\it inessential set} of edges in $G$.

Then it is natural for the alternation to be generalized to the
removal of several edges, which is just enough to enlarge the
domination number. That is the idea of the bondage number.

A measure of the efficiency of a domination in graphs was first
given by Bauer {\it et al.}\,\cite{bhns83} in 1983, who called this
measure as domination {\it line-stability}, defined as the minimum
number of lines (i.e. edges) which when removed from $G$ increases
$\gamma$.

In 1990, Fink {\it et al.}~\cite{fjkr90} formally introduced the
bondage number as a parameter for measuring the vulnerability of the
interconnection network under link failure. The minimum dominating
set of sites plays an important role in the network for it dominates
the whole network with the minimum cost. So we must consider whether
its function remains good under the with attack. Suppose that
someone such as a saboteur does not know which sites in the network
take part in the dominating role, but does know that the set of
these special sites corresponds to a minimum dominating set in the
related graph. Then how many links does he have to attack so that
the cost can not remains the same in order to dominate the whole
network? That minimum number of links is just the bondage number.

The {\it bondage number} $b(G)$ of a nonempty undirected graph $G$
is the minimum number of edges whose removal from $G$ results in a
graph with larger domination number. The precise definition of the
bondage number is defined as follows.
$$
b(G)=\min\{|B|:\ B\subseteq E(G),\gamma(G-B)>\gamma(G)\}.
$$
Since the domination number of every spanning subgraph of a nonempty
graph $G$ is at least as great as $\gamma(G)$, the bondage number of
a nonempty graph is well defined.

We call such an edge set $B$ that $\gamma(G-B)>\gamma(G)$ the {\it
bondage set} and the minimum one the {\it minimum bondage set}. In
fact, if $B$ is a minimum bondage set, then
$\gamma(G-B)=\gamma(G)+1$, because the removal of one single edge
can not increase the domination number by more than one. If $b(G)$
does not exist, for example empty graphs, we define $b(G)=\infty$.

It is quite difficult to compute the exact value of the bondage
number for general graphs since it strongly depends on the
domination number of the graphs. Much work focused on the bounds of
the bondage number as well as the restraints on particular classes
of graphs. The purpose of this lecture is to give a survey of
results and research methods related to these topics for graphs and
digraphs. For some results, we will give detailed proofs. For some
results and research methods, we will make some comments to develop
our study further.

The rest of the lecture is organized as follows. Section 2 gives
some preliminary results and complexity. Section 3 and Section 4
survey the study on the upper bounds and lower bounds, respectively.
The results for some special classes of graphs and planar graphs are
stated in Section 5 and Section 6, respectively. In Section 7, we
introduce some results on crossing number restraints. In Section 8
and Section 9, we are concerned about other and generalized types of
bondage numbers, respectively. In Section 10, we introduce some
results for digraphs. In the last section we introduce some results
for vertex-transitive graphs by applying efficient dominating sets.

%\vskip20pt

%\newpage

\section{Simplicity and Complexity}

As we have known from Introduction, the bondage number is an
important parameter for measuring the stability or the vulnerability
of a domination in a graph or a network. Our aim is to compute the
bondage number for any given graphs or networks. One has determined
the exact value of the bondage number for some graphs with simple
structure. For arbitrarily given graph, however, it has been proved
that determining its bondage number is NP-hard.

\subsection{Exact Values for Ordinary Graphs}

We begin our investigation of the bondage number by computing its
value for several well-known classes of graphs with simple
structure. In 1990, Fink {\it et al.}~\cite{fjkr90} proposed the
concept of the bondage number, and completely determined the exact
values of bondage numbers of some ordinary graphs, such as complete
graphs, paths, cycles and complete multipartite graphs.

By definition, to compute the exact value of bondage number for a
graph strongly depends upon its domination number. It is just that
the domination numbers for these graphs can be easily determined,
Fink {\it et al.}~\cite{fjkr90} determined the exact values of
bondage number for these graphs when they proposed the concept of
the bondage number.

%Use $K_n$, $C_n$ and $P_n$ and to denote a complete graph, cycle and
%path of order $n$, respectively. we have known that $\gamma(K_n)=1$,
%$\gamma(C_n)=\lceil\frac n3\rceil$ for $n\geqslant 3$ and
%$\gamma(P_n)=\lceil\frac n3\rceil$ for $n\geqslant 1$.

%Apart from establishing upper bounds on b(G), Fink et al. computed
%the bondage number of cycles, paths, and complete multipartite
%graphs and studied the bondage number of trees (several of these
%results can also be found in Bauer, Harary, Nieminen, and Suffel
%[l]). The purpose of this paper is to provide ties with analogous
%results for the fractional bondage number and for the discipline number.

\begin{thm}\label{thm2.1.1} \textnormal{(Fink {\it et al.}~\cite{fjkr90}, 1990)}

\noindent (a) For a complete graph $K_n$ of order $n\geqslant2$,
 $$
 \begin{array}{rl}b(K_n)=\left\lceil\frac{n}{2}\right\rceil;
 \end{array}
 $$

\noindent (b) For a path $P_n$ of order $n\geqslant2$,
   $$
   b(P_n)=\left\{\begin{array}{cl}%
   2&\ {\rm if}\ n\equiv 1({\rm mod}\ 3),\\
   1&\ {\rm otherwise};
   \end{array}\right.
   $$

\noindent (c) For a cycle $C_n$ of order $n$,
   $$
   b(C_n)=\left\{\begin{array}{cl}%
   3&\ {\rm if}\ n\equiv 1({\rm mod}\ 3),\\
   2&\ {\rm otherwise};
   \end{array}\right.
   $$

\noindent (d) For a complete $t$-partite graph
$G=K_{n_1,n_2,\ldots,n_t}$ with $n_1\leqslant
n_2\leqslant\cdots\leqslant n_t$ and $n_t>1$,
$$
b(G)=\left\{\begin{array}{cl}%
\lceil\frac{j}{2}\rceil&\ {\rm if}\ n_j=1\ {\rm and}\
n_{j+1}\geqslant2,\ {\rm for\ some}\ j,1\leqslant j< t,\\
2t-1&\ {\rm if}\ n_1=n_2=\cdots=n_t=2,\\
\sum\limits_{i=1}^{t-1}n_i&\ {\rm otherwise}.%
\end{array} \right.
$$
\end{thm}                                               % Proposition 4.1

\begin{pf}
We give the proof of the assertion (a). Clearly, $\gamma(K_n)=1$.
Let $H$ is a spanning subgraph of $K_n$ obtained by removing fewer
than $\lceil\frac n2\rceil$ edges from $K_n$. Then $H$ contains a
vertex of degree $n-1$, which can dominate all other vertices, and
hence $\gamma(H)=1$. Thus, $b(K_n)\geqslant \lceil\frac n2\rceil$.

If $n$ is even, the removal of a prefect matching from $K_n$ reduces
the degree of each vertex to $n-2$ and therefore yields a graph $H$
with $\gamma(H)=2$. If $n$ is odd, the removal of a prefect matching
from $K_n$ leaves a graph having exactly one vertex of degree $n-1$;
by removing one edge incident with this vertex, we obtain a graph
$H$ with $\gamma(H)=2$. In both cases, we can a spanning subgraph
$H$ by removal of $\lceil\frac{n}{2}\rceil$ edges from $K_n$ with
$\gamma(H)=2$. This implies
$b(K_n)\leqslant\lceil\frac{n}{2}\rceil$. Thus,
$b(K_n)=\lceil\frac{n}{2}\rceil$.

Now, we show the assertion (b). Since
$\gamma(C_n)=\gamma(P_n)=\left\lceil\frac{n}3\right\rceil$ for
$n\geqslant 3$, we see that $b(C_n)\geqslant 2$.

If $n\equiv 0,2 ({\rm mod}\,3)$, Then the graph $H$ obtained by
removing two adjacent edges from $C$, consists of an isolated vertex
and a path of order $n-1$. Thus,
 $$
 \begin{array}{rl}
 \gamma(H)=1+\gamma(P_{n-1})=1+\left\lceil\frac{n-1}3\right\rceil
 =1+\left\lceil\frac{n}3\right\rceil=1+\gamma(C_n),
 \end{array}
 $$
whence $b(C_n)\leqslant 2$, and so $b(C_n)=2$.

If $n\equiv 1 ({\rm mod}\,3)$, the removal of two edges from $C_n$
leaves a graph $H$ consisting of two paths $P_{n_1}$ and $P_{n_2}$,
where $n_1+n_2=n$. Then either $n_1\equiv n_2\equiv 2 ({\rm
mod}\,3)$, or, without loss of generality, $n_1\equiv 0 ({\rm
mod}\,3)$ and $n_2\equiv 1 ({\rm mod}\, 3)$. In the former case,
 $$
 \begin{array}{rl}
 \gamma(H) & =\gamma(P_{n_1})+\gamma(P_{n_2})=\lceil\frac{n_1}3\rceil+\lceil\frac{n_2}3\rceil\\
  &=\frac{n_1+1}3+\frac{n_2+1}3=\frac{n+2}3=\lceil\frac n3\rceil=\gamma(C_n).
  \end{array}
  $$
In the latter case,
 $$
 \begin{array}{rl}
 \gamma(H)=\frac{n_1}3 + \frac{n_2+2}3 = \frac{n+2}3=\left\lceil\frac n3\right\rceil =\gamma(C_n).
 \end{array}
 $$
In either case, when $n\equiv 1 ({\rm mod}\,3)$ we have
$b(C_n)\geqslant 3$.

Let $H$ the graph obtained from the deletion of three consecutive
edges of $C_n$. Then $H$ consists of two isolated vertices and a
path of order $n-2$. Thus,
 $$
 \begin{array}{rl}
 \gamma(H)=2+\left\lceil\frac{n-2}3\right\rceil=2+\frac{n-1}3
 =2+\left(\left\lceil\frac{n}3\right\rceil-1\right)
 =1+\gamma(C_n),
 \end{array}
 $$
so that $b(C_n)\leqslant 3$. Thus, $b(C_n)=3$.

As an immediate corollary to the assertion (b), we have the
assertion (c). The proof of the assertion (d) is left to the reader
as an exercise.
\end{pf}

\vskip6pt

Theorem~\ref{thm2.1.1} shows
$b(K_n)=\left\lceil\frac{n}{2}\right\rceil$ for an $(n-1)$-regular
graph $K_n$ of order $n\geqslant2$, $b(G)=n-1$ for an
$(n-2)$-regular graph $G$ of order $n\geqslant2$, where $G$ is a
$t$-partite graph $K_{n_1,n_2,\ldots,n_t}$ with $n_1=\cdots = n_t=2$
and $t=\frac{n}{2}$ for an even integer $n\ge 4$. For an
$(n-3)$-regular graph $G$ of order $n\ge 4$, Hu and Xu~\cite{hx11}
obtained the following result.

\begin{thm}\label{thm2.1.2}{\rm (Hu and Xu~\cite{hx11}, 2011)}\
 $b(G)=n-3$ for any {\rm ($n-3$)}-regular graph $G$ of order $n\ge 4$.
\end{thm}

The exact value of bondage number for a general graph, there is a
result as follows.

\begin{thm}\label{thm2.2} \textnormal{(Teschner \cite{t97}, 1997)}
If $G$ is a nonempty graph with a unique minimum dominating set,
then $b(G)=1$.
\end{thm}

\begin{pf}
Let $S$ be the $\gamma$-set of $G$, and let $x\notin S$.
Furthermore, let $y\in N_G(x)\cap S$. If $N_G(x)\cap S|\geqslant 2$
for each vertex $x\notin S$, then $S'=(S-\{y\})\cup\{x\}$ dominates
$G$ and $|S'|=|S|$, so that $S$ is a $\gamma$-set of $G$ as well,
which is a contradiction to the uniqueness of $S$. Thus,
$|N_G(x)\cap S|=1$ for a vertex $x\notin S$. Then
$\gamma(G-xy)>\gamma(G)$, which implies that $b(G)=1$.
\end{pf}

\vskip6pt

The following result is easy to verify.

\begin{thm}\label{thm2.2b} \textnormal{(Bauer {\it et al.}~\cite{bhns83}, 1983)}
If any vertex of a graph $G$ is adjacent with two or more pendant
vertices, then $b(G)=1$.
\end{thm}

Bauer {\it et al.}~\cite{bhns83} observed that the star is the
unique graph with the property that the bondage number is 1 and the
deletion of any edge results in the domination number increasing.
Hartnell and Rall~\cite{hr99} concluded by determining when this
very special property holds for higher bondage number. A graph is
called to be {uniformly bonded} if it has bondage number $b$ and the
deletion of any $b$ edges results in a graph with increased
domination number.

\begin{thm}\label{thm2.2c} \textnormal{(Hartnell and Rall~\cite{hr99}, 1999)}
The only uniformly bonded graphs with bondage number 2 are $C_3$ and
$P_4$. The unique graph with bondage number 3 that is uniformly
bonded is the graph $C_4$. There are no such graphs for bondage
number greater than 3.
\end{thm}

%\vskip6pt
\noindent{\bf Comments}\ As we mentioned above, to compute
the exact value of bondage number for a graph strongly depends upon
its domination number. In this sense, studying the bondage number
can greatly inspire one's research interesting to dominations.
However, determining the exact value of domination number for a
given graph is quite difficulty. In fact, even if the exact value of
the domination number for some graph is determined, it is still very
difficulty to compute the value of the bondage number for that
graph. For example, for the hypercube $Q_n$, we have
$\gamma(Q_n)=2^{n-1}$, but we have not yet determined $b(Q_n)$ for
any $n\geqslant 2$.

Perhaps Theorem~\ref{thm2.2} and Theorem~\ref{thm2.2b} provide an
approach to compute the exact value of bondage number for some
graphs by establishing some sufficient conditions for $b(G)=b$. In
fact, we will see later that Theorem~\ref{thm2.2} plays an important
role in determining the exact values of the bondage numbers for some
graphs. Thus, to study the bondage number, it is importance to
present various characterizations of graphs with a unique minimum
dominating set.

\subsection{Characterizations of Trees}
For trees, Bauer {\it et al.}~\cite{bhns83} in 1983 from the point
of view of the domination line-stability, independently, Fink {\it
et al.}~\cite{fjkr90} in 1990 from the point of view of the
domination edge-vulnerability, obtained the following result.

\begin{thm}\label{thm2.3}\ %\textnormal{\cite{ba83,Fink}}
If $T$ is a nontrivial tree, then $b(T)\leqslant2$.
\end{thm}

\begin{pf} If $\upsilon(T)=2$ then $b(T)=1$. Assume
$\upsilon(T)\geqslant3$ and let $(x_0,x_1,\ldots, x_k)$ be a longest
path in $T$. Clearly, $d_T(x_0)=d_T(x_k)=1$ and $k\geqslant 2$. If
$d_T(x_1)=2$ then $B=\{x_0x_1,x_1x_2\}$ is a bondage set, and so
$b(T)\leqslant |B|=2$. If $d_T(x_1)>2$, then $x_1$ is adjacent to
another vertex $y$ of degree one, the single edge $x_1y$ is bondage
set, and so $b(T)\leqslant 1$.
\end{pf}

\vskip6pt

It is natural to classify all trees according to their bondage
numbers. Fink {\it et al.}~\cite{fjkr90} proved that a forbidden
subgraph characterization to classify trees with different bondage
numbers is impossible, since they proved that if $F$ is a forest,
then $F$ is an induced subgraph of a tree $T$ with $b(T)=1$ and a
tree $T'$ with $b(T')=2$. However, they pointed out that the
complexity of calculating the bondage number of a tree is at most
$O(n^2)$ by methodically removing each pair of edges.

Even so, some characterizations, whether a tree has bondage number 1
or 2, have been found by several authors, see example
\cite{hr92,t97,tv00}.

First we describe the method due to Hartnell and Rall \cite{hr92},
by which all trees with bondage number $2$ can be constructed
inductively. An important tree $F_t$ in the construction is shown in
Figure~\ref{f1}. To characterize this construction, we need some
terminologies.

\vskip6pt

\begin{figure}[ht]
\begin{pspicture}(-5.2,0)(5,2)
  \cnode(3,0){3pt}{w}\rput(3.3,-0.08){$w$}
  \cnode(0,1){3pt}{u1}\cnode(1,1){3pt}{u2}\cnode(2,1){3pt}{u3}
  \cnode(3,1){3pt}{u4}\cnode(4,1){3pt}{u5}\cnode(5,1){3pt}{u6}
  \cnode(0,2){3pt}{v1}\cnode(1,2){3pt}{v2}\cnode(2,2){3pt}{v3}\cnode(3,2){3pt}{v4}
  \cnode(4,2){3pt}{v5}\cnode(5,2){3pt}{v6}
  \ncline{w}{u1}\ncline{w}{u3}\ncline{w}{u5}
  \ncline{u1}{u2}\ncline{u1}{v1}\ncline{u2}{v2}
  \ncline{u3}{u4}\ncline{u3}{v3}\ncline{u4}{v4}
  \ncline{u5}{u6}\ncline{u5}{v5}\ncline{u6}{v6}
\end{pspicture}
\caption{\label{f1}\footnotesize  Tree $F_t$}
\end{figure}
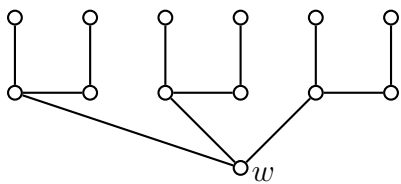

\noindent 1. {\it Attach} a path $P_n$ to a vertex $x$ of a tree
means to link $x$ and one end-vertex of the $P_n$ by an edge.

\noindent 2. {\it Attach} $F_t$ to a vertex $x$ means to link $x$
and a vertex $y$ of $F_t$ by an edge.

\vskip6pt

The following are four operations on a tree $T$:

\vskip6pt

\noindent Type 1: Attach a $P_2$ to $x\in V(T)$, where
$\gamma(T-x)=\gamma(T)$ and $x$ belongs to at least one $\gamma$-set
of $T$ (such a vertex exists, say, one end-vertex of $P_5$).

\noindent Type 2: Attach a $P_3$ to $x\in V(T)$, where
$\gamma(T-x)<\gamma(T)$.

\noindent Type 3: Attach $F_1$ to $x\in V(T)$, where $x$ belongs to
at least one $\gamma$-set of $T$.

\noindent Type 4: Attach $F_t$, $t\geqslant2$, to $x\in V(T)$, where
$x$ can be any vertex of $T$.

\vskip6pt

Let $\mathscr C=\{T:\,T$ is a tree and $T=K_1,T=P_4,T=F_t$ for some
$t\geqslant2$, or $T$ can be obtained from $P_4$ or
$F_t(t\geqslant2)$ by a finite sequence of operations of Type
$1,2,3,4\}$.

\begin{thm}\label{thm2.4} \textnormal{(Hartnell and Rall~\cite{hr92}, 1992)}
A tree has bondage number $2$ if and only if it belongs to $\mathscr
C$.
\end{thm}

Looking at different minimum dominating sets of a tree, Teschner
\cite{t97} presented a totally different characterization of the set
of trees having bondage number $1$. They defined a vertex to be {\it
universal} if it belongs to each minimum dominating set, and to be
{\it idle} if it does not belong to any minimum dominating set.

\begin{thm}\label{thm2.5} \textnormal{(Teschner~\cite{t97}, 1997)}
A tree $T$ has bondage number $1$ if and only if $T$ has a universal
vertex or an edge $xy$ satisfying

1) $x$ and $y$ are neither universal nor idle; and

2) all neighbors of $x$ and $y$\,(\,except for $x$ and $y$) are
idle.

\end{thm}

For a positive integer $k$, a subset $I\subseteq V(G)$ is called a
{\it $k$-independent set} (also called a {\it $k$-packing}) if
$d_G(x,y)>k$ for any two distinct vertices $x$ and $y$ in $I$. When
$k=1$, $1$-set is the normal independent set. The maximum
cardinality among all $k$-independent sets is called the {\it
$k$-independence number} (or {\it $k$-packing number}) of $G$,
denoted by $\alpha_k(G)$. A $k$-independent set $I$ is called an
{\it $\alpha_k$-set} if $|I|=\alpha_k(G)$. A graph $G$ is said to be
{\it $\alpha_k$-stable} if $\alpha_k(G)=\alpha_k(G-e)$ for every
edge $e$ of $G$. There are two important results on $k$-independent
sets.

\begin{prop}\label{prop2.2.4}
\textnormal{(Topp and Vestergaard~\cite{tv00}, 2000)} A tree $T$ is
$\alpha_k$-stable if and only if $T$ has a unique $\alpha_k$-set.
\end{prop}

\begin{prop}\label{prop2.2.5}
\textnormal{(Meir and Moon~\cite{mm75}, 1975)}
$\alpha_{2}(G)\leqslant\gamma(G)$ for any connected graph $G$ with
equality for any tree.
\end{prop}

Hartnell {\it et al.}~\cite{hjvw98}, independently, Topp and
Vestergaard~\cite{tv00}, also gave a constructive characterization
of trees with bondage number $2$ and, using Meir and Moon's result,
presented another characterization of those trees.

\begin{thm}\label{thm2.6} \textnormal{(Hartnell {\it et al.}~\cite{hjvw98}, 1998; Topp
and Vestergaard~\cite{tv00}, 2000)} $b(T)=2$ for a tree $T$ if and
only if $T$ has a unique $\alpha_2$-set.
\end{thm}

\begin{pf}
Let $T$ be a tree and let $e$ be any edge in $T$.

If $b(T)=2$, then $\gamma(T)=\gamma(T-e)$. By
Proposition~\ref{prop2.2.5}, we have
$\alpha_{2}(T)=\alpha_{2}(T-e)$, that is, $T$ is $\alpha_2$-stable.
It follows from Proposition~\ref{prop2.2.4} that $T$ has a unique
$\alpha_2$-set.

Conversely, if $T$ has a unique $\alpha_2$-set then, by
Proposition~\ref{prop2.2.4}, $T$ is $\alpha_2$-stable, that is,
$\alpha_{2}(T-e)=\alpha_{2}(T)$. It follows from
Proposition~\ref{prop2.2.5} that $\gamma(T-e)=\gamma(T)$, which
implies that $b(T)\geqslant 2$. By Theorem~\ref{thm2.3}, we have
$b(T)=2$.
\end{pf}

\vskip6pt

According to this characterization, Hartnell {\it et
al.}~\cite{hjvw98} presented a linear algorithm for determining the
bondage number of a tree.

\vskip6pt

\noindent{\bf Comments}\ In this subsection, we introduce three
characterizations for trees with bondage number $1$ or $2$. The
characterization in Theorem~\ref{thm2.4} is constructive,
constructing all trees with bondage number $2$, a natural and
straightforward method, by a series of graph-operations. The
characterization in Theorem~\ref{thm2.5} is a little advisable,
%ÓеãÃ÷ÖÇ
by describing the inherent property of trees with bondage number
$1$. The characterization in Theorem~\ref{thm2.6} is wonderful, by
using a strong graph-theoretic concept, $\alpha_k$-set. In fact,
this characterization is a byproduct of some results related to
$\alpha_k$-sets for trees. It is that this characterization closed
the relation between two concepts, the bondage number and the
$k$-independent set, and hence is of research hight and important
significance.

%$G$ is {\it $k$-stable} if $\gamma_k(G)=\gamma_k(G-e)$ for every
%edge $e$ of $G$.

\vskip20pt

\subsection{Complexity for General Graphs}

As mentioned above, the bondage number of a tree can be determined
within polynomial time. Indeed, in 1998, Hartnell et
at.~\cite{hjvw98} designed a linear time algorithm to compute the
bondage number of a tree. According to this algorithm, we can
determine within polynomial time the domination number of any tree
by removing each edge and verifying whether the domination number is
enlarged according to the known linear time algorithm for domination
numbers of trees.

However, it is impossible to find a polynomial time algorithm for
bondage numbers of general graphs. If such an algorithm $A$ exists,
then the domination number of any nonempty undirected graph $G$ can
be determined within polynomial time by repeatedly using $A$. Let
$G_0=G$ and $G_{i+1}=G_i-B_i$ where $B_i$ is the minimum edge set of
$G_i$ found by $A$ such that $\gamma(G_i-B_i)=\gamma(G_{i-1})+1$ for
$i=0,1,\ldots$; we can always find the minimum $B_i$ whose removal
from $G_i$ enlarges the domination number, until
$G_k=G_{k-1}-B_{k-1}$ is empty for some $k\geqslant1$, though
$B_{k-1}$ is not empty. Then
$\gamma(G)=\gamma(G_k)-k=\upsilon(G)-k$. As known to all, if $NP\ne
P$, the minimum dominating set problem is NP-complete, and so
polynomial time algorithms for the bondage number do not exist
unless $NP=P$.

In fact, Hu and Xu~\cite{hx12} have recently shown that the problem
determining the bondage number of general graphs is NP-hard.
%
%\begin{thm}
%Given a nonempty undirected graph $G$ and a positive integer
%$b\leqslant \varepsilon(G)$, determining wether or not
%$b(G)\leqslant b$ is NP-complete.
%\end{thm}
%
We first state the decision problem.

\vskip6pt\begin{prob} Consider the decision problem:

Bondage Problem

Instance: A graph $G$ and a positive integer $b$ $(\leqslant
\varepsilon(G))$.

Question: Is $b(G)\leqslant b$?
\end{prob}

Hu and Xu~\cite{hx12} showed that the bondage problem is NP-hard.
The basic way of the proof is to follow Garey and Johnson's
techniques for proving NP-hardness~\cite{gj79} by describing a
polynomial transformation from the known NP-complete problem:
$3$-satisfiability problem. To state the $3$-satisfiability problem,
we recall some terms.

Let $U$ be a set of Boolean variables. A {\it truth assignment} for
$U$ is a mapping $t: U\to\{T,F\}$. If $t(u)=T$, then $u$ is said to
be ``\,true" under $t$; if If $t(u)=F$, then $u$ is said to
be``\,false" under $t$. If $u$ is a variable in $U$, then $u$ and
$\bar{u}$ are {\it literals} over $U$. The literal $u$ is true under
$t$ if and only if the variable $u$ is true under $t$; the literal
$\bar{u}$ is true if and only if the variable $u$ is false.

A {\it clause} over $U$ is a set of literals over $U$. It represents
the disjunction of these literals and is {\it satisfied} by a truth
assignment if and only if at least one of its members is true under
that assignment. A collection $\mathscr C$ of clauses over $U$ is
{\it satisfiable} if and only if there exists some truth assignment
for $U$ that simultaneously satisfies all the clauses in $\mathscr
C$. Such a truth assignment is called a {\it satisfying truth
assignment} for $\mathscr C$. The $3$-satisfiability problem is
specified as follows.

\begin{center}
\begin{minipage}{130mm}
\setlength{\baselineskip}{24pt}

\vskip6pt\noindent {\bf $3$-satisfiability problem}:%(3SAT)}

\noindent {\bf Instance:}\ {\it A
%set $U=\{u_1,u_2,\cdots,u_n\}$ of finite variables, and a
collection $\mathscr{C}=\{C_1,C_2,\ldots,C_m\}$ of clauses over a
finite set $U$ of variables such that $|C_j| =3$ for $j=1,
2,\ldots,m$.}

\noindent {\bf Question:}\ {\it Is there a truth assignment for $U$
that satisfies all the clauses in $\mathscr{C}$?}

\end{minipage}
\end{center}

\begin{lem} \textnormal{(Theorem 3.1 in~\cite{gj79})}
The $3$-satisfiability problem is NP-complete.
\end{lem}

\begin{thm}\label{thm2.8}
\textnormal{(Hu and Xu~\cite{hx12}, 2012)} The bondage problem is
NP-hard.
\end{thm}

\begin{pf}
We show the NP-hardness of the bondage problem by transforming the
$3$-satisfiability problem to it in polynomial time.

Let $U=\{u_1,u_2,\ldots,u_n\}$ and $\mathscr{C}=\{C_1,C_2,
\ldots,C_m\}$ be an arbitrary instance of the $3$-satisfiability
problem. We will construct a graph $G$ and take a positive integer
$k$ such that $\mathscr{C}$ is satisfiable if and only if $b(G)\leq
k$. Such a graph $G$ can be constructed as follows.

For each $i=1,2,\ldots,n$, corresponding to the variable $u_i\in U$,
associate a triangle $T_i$ with vertex-set $\{u_i,\bar{u_i},v_i\}$.
For each $j=1,2,\ldots,m$, corresponding to the clause
$C_j=\{x_j,y_j,z_j\}\in \mathscr{C}$, associate a single vertex
$c_j$ and add an edge-set $E_j=\{c_jx_j, c_jy_j, c_jz_j\}$. Finally,
add a path $P=s_1s_2s_3$, join $s_1$ and $s_3$ to each vertex $c_j$
with $1\leqslant j\leqslant m$ and set $k=1$.
%Let $D$ be a $\gamma$-set in $G$.

Figure~\ref{f1} shows an example of the graph obtained when
$U=\{u_1,u_2,u_3,u_4\}$ and $\mathscr{C}=\{C_1,C_2,C_3\}$, where
$C_1=\{u_1,u_2,\bar{u_3}\}, C_2=\{\bar{u_1},u_2,u_4\},
C_3=\{\bar{u_2},u_3,u_4\}$.

\begin{figure}[ht]
\begin{center}
\begin{pspicture}(-5,-1.1)(5,6.7)

\cnode*(0,-.6){3pt}{s2}\rput(0,-1){$s_2$}
\cnode(-1,0){3pt}{s1}\rput(-1,-.4){$s_1$}
\cnode(1,0){3pt}{s3}\rput(1.1,-.4){$s_3$} \ncline{s2}{s1}
\ncline[linewidth=1.5pt]{s2}{s3}

\cnode(0,1.9){3pt}{c2}\rput(0,2.3){$c_2$}
\ncline[linewidth=1.5pt]{c2}{s1} \ncline{c2}{s3}
\cnode(-2.5,2){3pt}{c1}\rput(-2.8,1.7){$c_1$}
\ncline[linewidth=1.5pt]{c1}{s1} \ncline{c1}{s3}
\cnode(2.5,2){3pt}{c3}\rput(2.8,1.7){$c_3$}
\ncline[linewidth=1.5pt]{c3}{s1} \ncline{c3}{s3}

\cnode(-4.5,4){3pt}{u1}\rput(-4.7,4.3){$u_1$}
\cnode(-3,4){3pt}{u1'}\rput(-2.8,4.3){$\bar{u}_1$}
\ncline[linewidth=1.5pt]{u1}{u1'}
\cnode*(-2,4){3pt}{u2}\rput(-2.2,4.3){$u_2$}
\cnode(-0.5,4){3pt}{u2'}\rput(-0.3,4.3){$\bar{u}_2$}
\ncline[linewidth=1.5pt]{u2}{u2'}
\cnode(0.5,4){3pt}{u3}\rput(0.3,4.3){$u_3$}
\cnode*(2,4){3pt}{u3'}\rput(2.2,4.3){$\bar{u}_3$}
\ncline[linewidth=1.5pt]{u3}{u3'}
\cnode*(3,4){3pt}{u4}\rput(2.8,4.3){$u_4$}
\cnode(4.5,4){3pt}{u4'}\rput(4.7,4.3){$\bar{u}_4$}
\ncline[linewidth=1.5pt]{u4}{u4'}
\cnode*(-3.75,5.3){3pt}{v1}\rput(-3.75,5.6){$v_1$}
\ncline[linewidth=1.5pt]{v1}{u1} \ncline{v1}{u1'}
\cnode(-1.25,5.3){3pt}{v2}\rput(-1.25,5.6){$v_2$}
\ncline[linewidth=1.5pt]{v2}{u2} \ncline{v2}{u2'}
\cnode(1.25,5.3){3pt}{v3}\rput(1.25,5.6){$v_3$}
\ncline[linewidth=1.5pt]{v3}{u3} \ncline{v3}{u3'}
\cnode(3.75,5.3){3pt}{v4}\rput(3.75,5.6){$v_4$}
\ncline[linewidth=1.5pt]{v4}{u4} \ncline{v4}{u4'}

\ncline{c1}{u1} \ncline{c1}{u2} \ncline{c1}{u3'} \ncline{c2}{u1'}
\ncline{c2}{u2} \ncline{c2}{u4} \ncline{c3}{u2'} \ncline{c3}{u3}
\ncline{c3}{u4}
\end{pspicture}
\caption{\label{f1}\footnotesize An instance of the bondage problem
resulting from an instance of the $3$-satisfiability problem, in
which $U=\{u_1,u_2,u_3,u_4\}$ and
$\mathscr{C}=\{\{u_1,u_2,\bar{u_3}\},\{\bar{u_1},u_2,u_4\},
\{\bar{u_2},u_3,u_4\}\}$. Here $k=1$ and $\gamma=5$, where the set
of bold points is a $\gamma$-set.}
\end{center}
\end{figure}
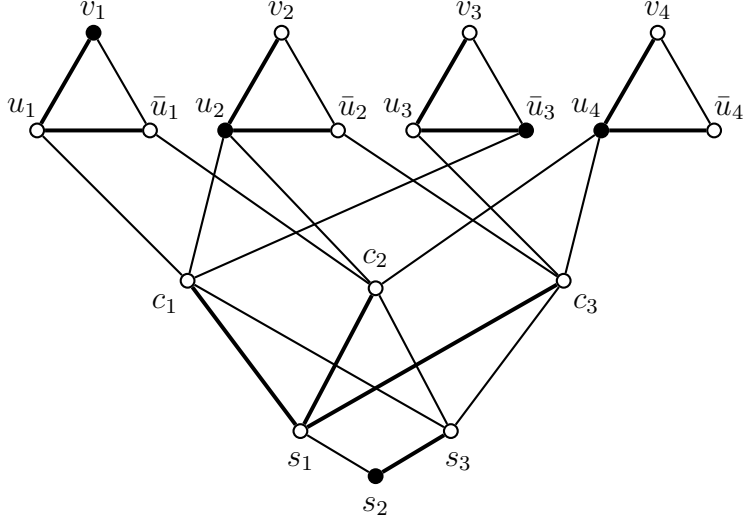

To prove that this is indeed a transformation, we must show that
$b(G)=1$ if and only if there is a truth assignment for $U$ that
satisfies all the clauses in $\mathscr{C}$.
%$\mathscr{C}$ is satisfiable.
This aim can be obtained by proving the following four claims.

\begin{description}

\item [Claim 3.1]
{\it $\gamma(G)\geq n+1$. Moreover, if $\gamma(G)=n+1$, then for any
$\gamma$-set $D$ in $G$, $D\cap V(P)=\{s_2\}$ and $|D\cap V(T_i)|=1$
for each $i=1,2,\ldots,n$, while $c_j\notin D$ for each
$j=1,2,\ldots,m$.}

\begin{pf}
Let $D$ be a $\gamma$-set of $G$. By the construction of $G$, the
vertex $s_2$ can be dominated only by vertices in $P$, which implies
$|D\cap V(P)|\geq 1$; for each $i=1,2,\ldots,n$, the vertex $v_i$
can be dominated only by vertices in $T_i$, which implies $|D\cap
V(T_i)|\geq 1$. It follows that $\gamma(G)=|D|\geq n+1$.

Suppose that $\gamma(G)=n+1$. Then $|D\cap V(P)|=1$ and $|D\cap
V(T_i)|=1$ for each $i=1,2,\ldots,n$. Consequently, $c_j\notin D$
for each $j=1,2,\ldots,m$. If $s_1\in D$, then $|D\cap V(P)|=1$
implies that $D\cap V(P)=\{s_1\}$, and so $s_3$ could not be
dominated by $D$, a contradiction. Hence $s_1\notin D$. Similarly
$s_3\notin D$ and, thus, $D\cap V(P)=\{s_2\}$ since $|D\cap
V(P)|=1$.
\end{pf}

\item [Claim 3.2]
{\it $\gamma(G)= n+1$ if and only if $\mathscr{C}$ is satisfiable.}

\begin{pf}
Suppose that $\gamma(G)=n+1$ and let $D$ be a $\gamma$-set of $G$.
By Claim 3.1, for each $i=1,2,\ldots,n$, $|D\cap V(T_i)|=1$, it
follows that $D\cap V(T_i)=\{u_i\}$ or $D\cap V(T_i)=\{\bar{u_i}\}$
or $D\cap V(T_i )=\{v_i\}$. Define a mapping $t: U\to \{T,F\}$ by
 \begin{equation}\label{e3.1}
 t(u_i)=\left\{
 \begin{array}{l}
 T \ \ {\rm if}\ u_i\in D \ {\rm or} \ v_i\in D, \\
 F \ \ {\rm if}\ \bar {u_i}\in D,
\end{array}
 \right.
 \ i=1,2,\ldots,n.
 \end{equation}

We will show that $t$ is a satisfying truth assignment for
$\mathscr{C}$. It is sufficient to show that every clause in
$\mathscr{C}$ is satisfied by $t$. To this end, we arbitrarily
choose a clause $C_j\in\mathscr{C}$ with $1\leqslant j\leqslant m$.
%Without loss of generality,
Since the corresponding vertex $c_j$ in $G$ is adjacent to neither
$s_2$ nor $v_i$ for any $i$ with $1\leqslant i\leqslant n$, there
exists some $i$ with $1\leqslant i\leqslant n$ such that $c_j$ is
dominated by $u_i\in D$ or $\bar{u}_i\in D$. Suppose that $c_j$ is
dominated by $u_i\in D$. Since $u_i$ is adjacent to $c_j$ in $G$,
the literal $u_i$ is in the clause $C_j$ by the construction of $G$.
Since $u_i\in D$, it follows that $t(u_i)=T$ by (\ref{e3.1}), which
implies that the clause $C_j$ is satisfied by $t$. Suppose that
$c_j$ is dominated by $\bar{u}_i\in D$. Since $\bar{u}_i$ is
adjacent to $c_j$ in $G$, the literal $\bar{u}_i$ is in the clause
$C_j$. Since $\bar{u}_i\in D$, it follows that $t(u_i)=F$ by
(\ref{e3.1}). Thus, $t$ assigns $\bar{u}_i$ the truth value $T$,
that is, $t$ satisfies the clause $C_j$. By the arbitrariness of $j$
with $1\leqslant j\leqslant m$, we show that $t$ satisfies all the
clauses in $\mathscr{C}$, that is, $\mathscr{C}$ is satisfiable.

Conversely, suppose that $\mathscr{C}$ is satisfiable, and let $t:
U\to \{T,F\}$ be a satisfying truth assignment for $\mathscr{C}$.
Construct a subset $D'\subseteq V(G)$ as follows. If $t(u_i)=T$,
then put the vertex $u_i$ in $D'$; if $t(u_i)=F$, then put the
vertex $\bar{u_i}$ in $D'$. Clearly, $|D'|=n$. Since $t$ is a
satisfying truth assignment for $\mathscr{C}$, for each
$j=1,2,\ldots,m$, at least one of literals in $C_j$ is true under
the assignment $t$. It follows that the corresponding vertex $c_j$
in $G$ is adjacent to at least one vertex in $D'$ since $c_j$ is
adjacent to each literal in $C_j$ by the construction of $G$. Thus
$D'\cup \{s_2\}$ is a dominating set of $G$, and so $\gamma(G)\leq
|D'\cup \{s_2\}|=n+1$. By Claim 3.1, $\gamma(G)\geq n+1$, and so
$\gamma(G)=n+1$.
\end{pf}

\item [Claim 3.3]
{\it $\gamma(G-e)\leqslant n+2$ for any $e\in E(G)$.}

\begin{pf}
Let $E_1=\{s_2s_3,s_1c_j,u_i\bar{u_i},u_iv_i,: i=1,2,\ldots,n;
j=1,2,\ldots,m\}$ (induced by heavy edges in Figure~\ref{f1}) and
let $E_2=E(G)\setminus E_1$. Assume $e\in E_2$. Let
$D'=\{u_1,u_2,\ldots,u_n,s_1, s_2\}$. Clearly, $D'$ is a dominating
set of $G-e$ since every vertex not in $D'$ is incident with some
vertex in $D'$ via an edge in $E_1$. Hence, $\gamma(G-e)\leq
|D'|=n+2$. Now assume $e\in E_1$. Let $D''=\{u_1,u_2,\ldots,u_n,s_2,
s_3\}$. If $e$ is either $s_2s_3$ or incident with the vertex $s_1$,
then $D''$ is a dominating set of $G-e$, clearly. If $e$ is either
$u_i\bar{u_i}$ or $u_iv_i$ for some $i$ ($1\leqslant i\leqslant n$),
then we use the vertex either $v_i$ or $\bar{u_i}$ instead of $u_i$
in $D''$ to obtain $D'''$; and hence $D'''$ is a dominating set of
$G-e$. These facts imply that $\gamma(G-e)\leqslant n+2$.
\end{pf}

\item [Claim 3.4]
{\it $\gamma(G)=n+1$ if and only if $b(G)=1$.}

\begin{pf}
Assume $\gamma(G)=n+1$ and consider the edge $e=s_1s_2$. Suppose
$\gamma(G)=\gamma(G-e)$. Let $D'$ be a $\gamma$-set in $G-e$. It is
clear that $D'$ is also a $\gamma$-set of $G$. By Claim 3.1 we have
$c_j\notin D'$ for each $j=1,2,\ldots,m$ and $D'\cap V(P)=\{s_2\}$.
But then $s_1$ is not dominated by $D'$, a contradiction. Hence,
$\gamma(G)<\gamma(G-e)$, and so $b(G)=1$.

Now, assume $b(G)=1$. By Claim 3.1, we have that $\gamma(G)\geq
n+1$. Let $e'$ be an edge such that $\gamma(G)<\gamma(G-e')$. By
Claim 3.3, we have that $\gamma(G-e')\leqslant n+2$. Thus, $n +
1\leqslant \gamma(G)< \gamma(G-e')\leqslant n+2$, which yields
$\gamma(G)=n+1$.
\end{pf}

\end{description}

By Claim 3.2 and Claim 3.4, we prove that $b(G)=1$ if and only if
there is a truth assignment for $U$ that satisfies all the clauses
in $\mathscr{C}$. Since the construction of the bondage instance is
straightforward from a $3$-satisfiability instance, the size of the
bondage instance is bounded above by a polynomial function of the
size of $3$-satisfiability instance. It follows that this is a
polynomial transformation.

The theorem follows.
\end{pf}

\vskip6pt

\noindent{\bf Comments}\ Theorem~\ref{thm2.8} shows that we are
unable to find a polynomial time algorithm to determine bondage
numbers of general graphs unless $NP=P$. At the same time, this
result also shows that the following study is of important
significance.

\begin{itemize}

 \item
 Find approximation polynomial algorithms with performance ratio as
small as possible.

 \item
Find the lower and upper bounds with difference as small as
possible.

 \item
Determine exact values for some graphs, specially
 well-known networks.

 \end{itemize}

Unfortunately, we can not proved whether or not determining the
bondage is NP-problem, since for any subset $B\subset E(G)$, it is
not clear that there is a polynomial algorithm to verify
$\gamma(G-B)>\gamma(G)$. Since the problem of determining the
domination number is NP-complete, we conjecture that it is not in
$NP$. This is a worthwhile task to study further.

However, Hartnell et at.~\cite{hjvw98} designed a linear time
algorithm to compute the bondage number of a tree.  Motivated by
this fact, we can made an attempt to consider whether there is a
polynomial time algorithm to compute the bondage number for some
special classes of graphs such as planar graphs, Cayley graphs, or
graphs with some restrictions of graph-theoretical parameters such
as degree, diameter, connectivity, domination number and so on.

\vskip20pt
%\newpage
\section{Upper Bounds}

By Theorem~\ref{thm2.8}, since we can not find a polynomial time
algorithm for determining the exact values of bondage numbers of
general graphs, it is weightily significative to establish some
sharp bounds of the bondage number of a graph. In this section, we
survey several known upper bounds of the bondage number in terms of
some other graph-theoretical parameters.

\subsection{Most Basic Upper Bounds}

Along with the exact values of bondage numbers for some ordinary
graphs computed, several general upper bounds were also derived. In
this subsection, we will survey some simple and important upper
bounds in terms of the sum of degrees of two vertices with distance
1 or 2. To show the simpleness of these upper bounds, we give their
proofs. We start this subsection with an easy observation.

\begin{lem}\label{lem3.1.1}
\textnormal{(Teschner~\cite{t97}, 1997)} Let $H$ be a spanning
subgraph obtained by removing $k$ edges from a graph $G$. Then
$b(G)\leqslant b(H)+k$.
\end{lem}                                                     % Lemma 3.1

\begin{pf}
Let $E'=E(G)\setminus E(H)$ and $B$ be a bondage set of $H$. Then
$\gamma(G-E'-B)=\gamma(H-B)>\gamma(H)\geqslant\gamma(G)$ and so
$b(G)\leqslant|B|+|E'|=b(H)+k$.
\end{pf}

\vskip6pt If we select a spanning subgraph $H$ such that $b(H)=1$,
then Lemma~\ref{lem3.1.1} yields some upper bounds on the bondage
number of a graph.

\begin{thm}\label{thm3.1.2}
\textnormal{(Bauer {\it et al.}~\cite{bhns83}, 1983)} If there
exists at least one vertex $x$ in a graph $G$ such that
$\gamma(G-x)\geqslant\gamma(G)$, then $b(G)\leqslant
d_G(x)\leqslant\Delta(G)$.
\end{thm}                                                     % thm3.1.2

\begin{pf}
Let $x$ be a vertex in $G$ such that
$\gamma(G-x)\geqslant\gamma(G)$. Then $x$ is not in any $\gamma$-set
of $G$, and so it is dominated by some $y\in N_G(x)$. Let
$H=G-E_x+xy$. Then $H$ is a spanning subgraph obtained by removing
$k$ edges from $G$, where $k=d(x)-1$. It is clear that the removal
of the edge $xy$ from $H$ results in increase of the domination
number of $H$, and so $b(H)=1$. The result follows from Lemma
\ref{lem3.1.1} immediately.
\end{pf}

\vskip6pt

The following early result obtained by Bauer {\it et
al.}~\cite{bhns83} and Fink {\it et al.}~\cite{fjkr90},
respectively, can be derived from Lemma~\ref{lem3.1.1}.

\begin{thm}\label{thm3.1.3}
$b(G)\leqslant d(x)+d(y)-1$ for any two adjacent vertices $x$ and
$y$ in a graph $G$, that is,
 $$
 b(G)\leqslant\min\limits_{xy\in E(G)}\{d_G(x)+d_G(y)-1\}.
 $$
\end{thm}                                                     % thm3.1.3

\begin{pf}
Let $E_{xy}$ denote the set of edges that are incident with at least
one of $x$ and $y$, but not both, and let $H=G-E_{xy}$. Then $H$ is
a spanning subgraph obtained by removing $k$ edges from $G$, where
$k=d(x)+d(y)-2$. It is clear that the removal of the edge $xy$ from
$H$ results in increase of the domination number of $H$, and so
$b(H)=1$. The result follows from Lemma \ref{lem3.1.1} immediately.
\end{pf}

\vskip6pt

\vskip6pt This theorem gives a natural corollary obtained by several
authors.

\begin{cor}\label{con3.1.4}
\textnormal{(Bauer {\it et al.}~\cite{bhns83}, 1983; Fink {\it et
al.}~\cite{fjkr90}, 1990)} If $G$ is a graph without isolated
vertices, then $b(G)\leqslant \Delta(G)+\delta(G)-1$.
\end{cor}                                                 % Cor3.4

\vskip6pt

In 1999, Hartnell and Rall~\cite{hr99} extended
Theorem~\ref{thm3.1.3} to the following more general case, which can
be also derived from Lemma \ref{lem3.1.1} by adding
``$H=G-E_x-E_y+(x,z,y)$ if $d_G(x,y)=2$, where $(x,z,y)$ is a path
of length $2$ in $G$\," in the above proof of
Theorem~\ref{thm3.1.3}.

\begin{thm}\label{thm3.1.5}\textnormal{(Hartnell and Rall~\cite{hr99}, 1999)}
\ $b(G)\leqslant d(x)+d(y)-1$ for any distinct two vertices $x$ and
$y$ in a graph $G$ with $d_G(x,y)\leqslant2$, that is,
 $$
 b(G)\leqslant\min\limits_{d_G(x,y)\leqslant2}\{d_G(x)+d_G(y)-1\}.
 $$
\end{thm}                                                     % thm3.1.5

%\begin{pf}
%For every $x\in V(G)$, let $E_x$ be the set of edges incident with
%$x$. Let $H=G-E_x\cup E_y+P$, where $P$ a path in $G$ and $P=(x,y)$
%if $d_G(x,y)=1$ or $P=(x,z,y)$ if $d_G(x,y)=2$. It is clear that $P$
%is a connected component of $H$. The removal of an edge of $P$
%results in increase of the domination number of $H$, and so
%$b(H)=1$. The result follows form Lemma \ref{lem3.1.1}, in which
%$k=d(x)+d(y)-2$.
%\end{pf}
%
%
%\vskip6pt

\begin{cor}\label{con3.1.6} \textnormal{(Fink {\it et al.}~\cite{fjkr90}, 1990)}
If a vertex of a graph $G$ is adjacent with two or more vertices of
degree one, then $b(G)= 1$.
\end{cor}

We remark that the bounds stated in Corollary~\ref{con3.1.4} and
Theorem~\ref{thm3.1.5} are sharp. As indicated by
Theorem~\ref{thm2.1.1}, one class of graphs in which the bondage
number achieves these bounds is the class of cycles whose orders are
congruent to 1 modulo 3.

On the other hand, Hartnell and Rall~\cite{hr94} sharpened the upper
bound in Theorem~\ref{thm3.1.3} as follows, which can be also
derived from Lemma~\ref{lem3.1.1}.
%The following is another fundamental result, which is an improvement of .

\begin{thm}\label{thm3.1.7} \textnormal{(Hartnell and Rall~\cite{hr94}, 1994)}
$b(G)\leqslant d_G(x)+d_G(y)-1-|N_G(x)\cap N_G(y)|$ for any two
adjacent vertices $x$ and $y$ in a graph $G$, that is,
 $$
 b(G)\leqslant\min\limits_{xy\in E(G)}\{d_G(x)+d_G(y)-1-|N_G(x)\cap N_G(y)\}.
 $$
\end{thm}                                                     % thm3.1.7

%\begin{pf} Ô­ÎÄÖеÄÖ¤Ã÷
%Let $u$ be a vertex in $G$ and $x\in N_G(u)$ achieve the minimum in
%the statement of the theorem. If $F_u=\{uy:\ y\in N_G(u)\}$, then in
%the graph $H=G-F_u-E_G(\{x\},V-N[u])$, $u$ is an isolated vertex and
%all of neighbors of $x$ in $H$, if any, lie in $N_G(u)$.
%
%Let $A$ be a $\gamma$-set of $H$ then $u\in A$ and $A\cap N(u)\ne
%\emptyset$ since, in particular, $A$ must dominates $x$ in $H$.
%However, $A-\{u\}$ is a dominating set for $G$, and so
%$\gamma(G)\leqslant |A|-1$. Hence $b(G)\leqslant |F_u\cup
%E_G(\{x\},V-N[u])|$ as claimed.
%\end{pf}

\begin{pf}
Let $F_y=\{yu:\ u\in N_G(y)\setminus N_G(x)\}$ and $H=G-E_x\cup
F_y+xy$. Then $H$ is a spanning subgraph obtained by removing $k$
edges from $G$, where $k=|E_x\cup F_y|+1=d(x)+d(y)-|N_G(x)\cap
N_G(y)|-2$. It is clear that at least one of $x$ and $y$ is in any
$\gamma$-set of $H$, and so the removal of $e$ results in increase
of the domination number of $H$, and so $b(H)=1$. The result follows
from Lemma \ref{lem3.1.1} immediately.
\end{pf}

\vskip6pt These results give simple but important upper bounds on
the bondage number of a graph, and is also the foundation of almost
all results on bondage numbers upper bounds obtained till now.

By careful consideration of the nature of the edges from the
neighbors of $x$ and $y$, Wang~\cite{wang96} further refined the
bound in Theorem~\ref{thm3.1.7}. For any edge $xy\in E(G)$, $N_G(y)$
contains the following four subsets.

1) $T_1(x,y)=N_G[x]\cap N_G(y)$;

2) $T_2(x,y)=\{w\in N_G(y):N_G(w)\subseteq N_G(y)-x\}$;

3) $T_3(x,y)=\{w\in N_G(y):N_G(w)\subseteq N_G(z)-x$ for some $z\in
N_G(x)\cap N_G(y)\}$; and

4) $T_4(x,y)=N_G(y)\setminus(T_1(x,y)\cup T_2(x,y)\cup T_3(x,y))$.

\vskip20pt
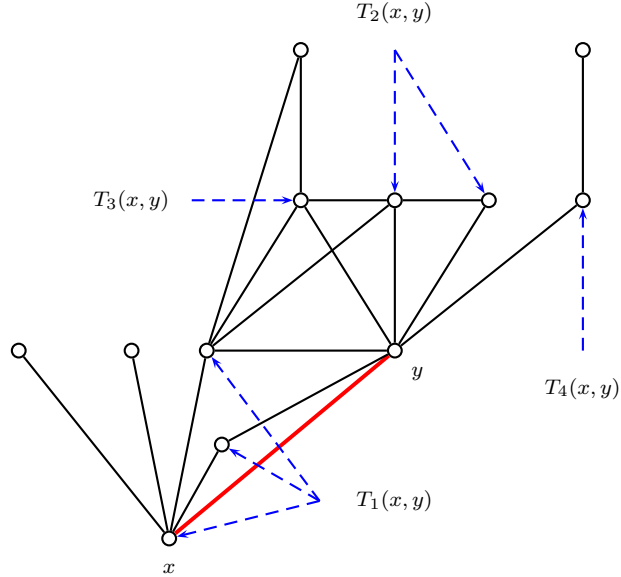
\begin{figure}[ht]
\begin{pspicture}(-.8,2)(12,9)
  \cnode(5,2.5){3pt}{x}\rput(5,2.1){\scriptsize $x$}
  \cnode(8,5){3pt}{y}\rput(8.3,4.7){\scriptsize $y$}
  \cnode(3,5){3pt}{x1}\cnode(4.5,5){3pt}{x2}\cnode(5.5,5){3pt}{x3}\cnode(5.7,3.75){3pt}{x4}
  \cnode(6.75,7){3pt}{y1}\cnode(8,7){3pt}{y2}\cnode(9.25,7){3pt}{y3}\cnode(10.5,7){3pt}{y4}
  \cnode(6.75,9){3pt}{z1}\cnode(10.5,9){3pt}{z2}
  \ncline[linecolor=red,linewidth=1.5pt]{x}{y}\ncline{x}{x1}\ncline{x}{x2}\ncline{x}{x3}\ncline{x}{x4}
  \ncline{y}{x3}\ncline{y}{x4}\ncline{y}{y1}\ncline{y}{y2}\ncline{y}{y3}\ncline{y}{y4}
  \ncline{x3}{y1}\ncline{x3}{y2}\ncline{x3}{z1}
  \ncline{y2}{y1}\ncline{y2}{y3}
  \ncline{y1}{z1}\ncline{y4}{z2}
  \psset{linecolor=blue,arrows=->,linestyle=dashed}
  \pnode(7,3){T1}\rput(8,3){\scriptsize $T_1(x,y)$}
  \ncline{T1}{x}\ncline{T1}{x3}\ncline{T1}{x4}
  \pnode(8,9){T2}\rput(8,9.5){\scriptsize $T_2(x,y)$}
  \ncline{T2}{y2}\ncline{T2}{y3}
  \pnode(5.3,7){T3}\rput(4.5,7){\scriptsize $T_3(x,y)$}
  \ncline{T3}{y1}
  \pnode(10.5,5){T4}\rput(10.5,4.5){\scriptsize $T_4(x,y)$}
  \ncline{T4}{y4}
  \end{pspicture}
\caption{\label{f2}\footnotesize  Illustration of $T_i(x,y)$ for
$i=1,2,3,4$}
\end{figure}

%\vskip6pt

The illustrations of $T_1(x,y)$, $T_2(x,y)$, $T_3(x,y)$ and
$T_4(x,y)$ are shown in Figure \ref{f2} (corresponding vertices
pointed by dashed arrows).

\begin{thm}\label{thm3.1.8} \textnormal{(Wang\,\cite{wang96},
1996)} Let $G$ be a nonempty graph. Then
 $$
 b(G)\leqslant\min\limits_{xy\in E(G)}\{d_G(x)+|T_4(x,y)|\}.
 $$
\end{thm}                                                    %thm3.1.8

\begin{pf}
Let $H=G-E_x\cup T_4(x,y)+xy$. Then $H$ is a spanning subgraph
obtained by removing $k$ edges from $G$, where $k=|E_x\cup
F_y|+1=d(x)+|T_4(x,y)|-1$. Without loss of generality, assume
$d_y(G)\leqslant d_G(x)$. If $y$ is not in any $\gamma$-set of $H$,
then $x$ must be in every $\gamma$-set of $G$, and so $y$ is not in
any $\gamma$-set of $G$. By Theorem~\ref{thm3.1.2}, $b(G)\leqslant
d_G(y)\leqslant d_G(x)$. Assume that $y$ is in some $\gamma$-set of
$H$ below. Then $x$ can be dominated by $y$ in $H$. Thus, the edge
$\{xy\}$ is a bondage set of $H$, and so $b(H)=1$. By
Lemma~\ref{lem3.1.1}, $b(G)\leqslant 1+k\leqslant
d_G(x)+|T_4(x,y)|$.
\end{pf}

\vskip6pt

The graph shown in Figure~\ref{f2} shows that the upper bound given
in Theorem~\ref{thm3.1.8} is better than that in
Theorem~\ref{thm3.1.5} and Theorem~\ref{thm3.1.7}, for the upper
bounds obtained from these two theorems are $d_G(x)+d_G(y)-1=11$ and
$d_G(x)+d_G(y)-|N_G(x)\cap N_G(y)|=9$, respectively, while the upper
bound given by Theorem \ref{thm3.1.8} is $d_G(x)+|T_4(x,y)|=6$.

\vskip6pt

The following result is also an improvement of
Theorem~\ref{thm3.1.3}, in which $t=2$.

\begin{thm}\label{thm3.1.9} \textnormal{(Teschner \cite{t97},
1997)} If $G$ contains a complete subgraph $K_t$ with $t\geqslant
2$, then $b(G)\leqslant \min\limits_{xy\in
E(K_t)}\{d_G(x)+d_G(y)-t+1\}$.
\end{thm}

\begin{pf}
Let $xy\in E(K_t)$, $F_y$ be the set of edges incident with $y$ in
$G$ but in $K_t$, and let $H=G-F_y+xy$. Any dominating set $S$ in
$H-x$ includes a vertex of $N_H[y]$ because $y$ has to be dominated.
But then $S$ is also a dominating set of $H$ and therefore
$\gamma(H)\leqslant\gamma(H-x)$. By Theorem~\ref{thm3.1.2},
$b(H)\leqslant d_H(x)$. Thus, by Lemma~\ref{lem3.1.1},
$b(G)\leqslant b(H)+|F_y|+1\leqslant d_G(x)+d_G(y)-(t-1)$.
\end{pf}

\vskip6pt Following Fricke {\it et al.}~\cite{Fhhhl02}, a vertex $x$
of a graph $G$ is {\it $\gamma$-good} if $x$ belongs to some
$\gamma$-set of $G$ and {\it $\gamma$-bad} if $x$ belongs to no
$\gamma$-set of $G$. Let $A(G)$ be the set of $\gamma$-good
vertices,
%with degree at least one,
and let $B(G)$ be the set of $\gamma$-bad vertices in $G$. Clearly,
$\{A(G),B(G)\}$ is a partition of $V(G)$. Note there exists $x\in A$
such that $\gamma(G-x)=\gamma(G)$, say, one end-vertex of $P_5$.
Samodivkin~\cite{s08} presented some sharp upper bounds for $b(G)$
in terms of $\gamma$-good and $\gamma$-bad vertices of $G$.

\begin{thm}\label{thm3.1.10}
\textnormal{(Samodivkin~\cite{s08}, 2008)} Let $G$ be a graph.

(i)\ Let $C(G)=\{x\in V(G): \gamma(G-x)\geqslant\gamma(G)\}$. If
$C(G)\ne\emptyset$, then
 $$
 b(G)\leqslant\min\{d_G(x)+\gamma(G)-\gamma(G-x): x\in C(G)\}.
 $$

(ii)\ If $B\ne\emptyset$, then
 $$
 b(G)\leqslant\min\{|N_G(x)\cap A|: x\in B(G)\}.
 $$
\end{thm}

\begin{pf}
Notice that if $x$ is an isolated vertex in $G$ then $x\in A(G)$,
and if $\gamma(G-x)>\gamma(G)$ then $x$ is not an isolated vertex
and is in every $\gamma$-set of $G$.

(i)\ Let $x\in C(G)$ and let $\gamma(G-x)=\gamma(G)+p$. Then
$p\geqslant 0$. If $p=0$, then $b(G)\leqslant d_G(x)$ by
Theorem~\ref{thm3.1.2}.

Now assume $p\geqslant 1$. Then $\gamma(G-x)>\gamma(G)$. By the
above explanation, it follows that $x$ is in every $\gamma$-set of
$G$. Let $S$ be a $\gamma$-set of $G$. Then $S'=(S-\{x\})\cup
N_G(x)$ is a dominating set of $G-x$ which implies
$\gamma(G)+p=\gamma(G-x)\leqslant |S'|=\gamma(G)-1+d_G(x)$. Hence
$1\leqslant p\leqslant d_G(x)-1$. Let $F\subseteq E_G(x)=E_x$ with
$|F|=d_G(x)-p$. Then
 $$
\gamma(G-F)\geqslant\gamma(G-E_x)-p=\gamma(G-x)+1-p=\gamma(G)+1,
$$
which implies $b(G)\leqslant |F|=d_G(x)+\gamma(G)-\gamma(G-x)$.

(ii)\ Let $x$ be any vertex in $B(G)$. Then $N_G(x)\cap
A\ne\emptyset$. Let $y\in N_G(x)\cap A$ such that
$\gamma(G-xy)=\gamma(G)$. Such an edge does exist since $x$ is
$\gamma$-bad. Notice that every $\gamma$-set of $G-xy$ is a
$\gamma$-set of $G$. Thus, $A(G-xy)\subseteq A(G)$ and
$B(G-xy)\supseteq B(G)$.

To prove (ii), we only need to prove that
$\gamma(G-E_G(x,A))>\gamma(G)$. Assume to the contrary that there is
some $x\in B(G)$ such that $\gamma(G-E_G(x,A))=\gamma(G)$. Let
$H=G-E_G(x,A)$. By the above discussion, $B(H)\supseteq B(G)$ which
implies $N_H[x]\subseteq B(H)$. But this is clearly impossible.
\end{pf}

\begin{thm}\label{thm3.1.11}
\textnormal{(Samodivkin~\cite{s08}, 2008)} Let $G$ be a graph. If
$A^+(G)=\{x\in A(G): d_G(x)\geqslant 1\ {\rm and}\
\gamma(G-x)<\gamma(G)\}\ne\emptyset$, then
 $$
 b(G)\leqslant\min\limits_{x\in A^+(G),y\in B(G-x)}\{d_G(x)+
|N_G(y)\cap A(G-x)|\}.
$$
\end{thm}

\begin{pf}
Let $x\in A^+(G)$ and $S$ be a $\gamma$-set of $G-x$. Then clearly
no neighbor of $x$ is in $S$, which implies $\emptyset\ne
N_G(x)\subseteq B(G-x)$. Since $\gamma(G-E_x)=\gamma(G)$ it follows
that $b(G)\leqslant d_G(x)+b(G-x)$.

By the assertion (ii) in Theorem~\ref{thm3.1.10}, $b(G-x)\leqslant
|N_G(y)\cap A(G-x)|$ for any $y\in B(G-x)$. Hence $b(G)\leqslant
d_G(x)+|N_G(y)\cap A(G-x)$.
\end{pf}

\begin{prop}\label{prop3.1.12}\textnormal{(Samodivkin~\cite{s08}, 2008)}
Under the notation of Theorem~\ref{thm3.1.8}, if $x\in A^+(G)$, then
$(T_1(x,y)-\{x\})\cup T_2(x,y)\cup T_3(x,y)\subseteq N_G(y)\setminus
B(G-x)$.
\end{prop}

By Proposition~\ref{prop3.1.12}, if $x\in A^+(G)$, then
 $$
 \begin{array}{rl}
 d_G(x)+\min\limits_{y\in N_G(x)}\{|T_4(x,y)|\}
 &\geqslant d_G(x)+\min\limits_{y\in N_G(x)}\{|N_G(y)\cap A(G-x)|\}\\
 & \geqslant d_G(x)+\min\limits_{y\in B(G-x)}\{|N_G(y)\cap A(G-x)|\}.
 \end{array}
 $$
Hence Theorem~\ref{thm3.1.8} can be seen to follow from
Theorem~\ref{thm3.1.11}. Any graph $G$ with $b(G)$ achieving the
upper bound of some of Theorem~\ref{thm3.1.8} can be used to show
that the bound of Theorem~\ref{thm3.1.11} is sharp.

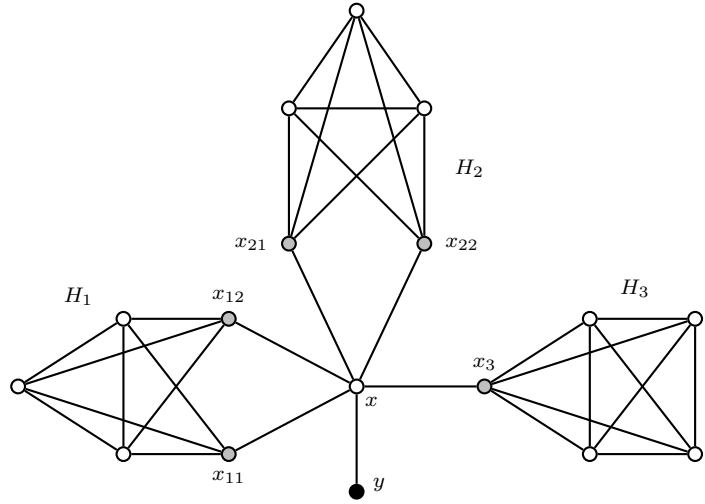
\begin{figure}[h]  % h: here; t: top; b: bottom
%\psset{unit=.9}     % scale
\psset{arrowsize=.14}

\begin{center}
\begin{pspicture}(-4.6,0)(4.6,7)
\cnode(-4.5,1.4){3pt}{1} \cnode(-3.1,0.5){3pt}{2}
\cnode(-3.1,2.3){3pt}{3}
\cnode[fillstyle=solid,fillcolor=lightgray](-1.7,0.5){3pt}{4}
\rput(-1.7,0.2){\scriptsize$x_{11}$}
\cnode[fillstyle=solid,fillcolor=lightgray](-1.7,2.3){3pt}{5}\rput(-1.7,2.6){\scriptsize$x_{12}$}
\cnode*(0,0){3pt}{6}\rput(.3,.1){\scriptsize$y$}
\cnode(0,1.4){3pt}{7}\rput(.2,1.2){\scriptsize$x$}
\cnode[fillstyle=solid,fillcolor=lightgray](1.7,1.4){3pt}{8}\rput(1.7,1.7){\scriptsize$x_3$}
\cnode(3.1,0.5){3pt}{9} \cnode(3.1,2.3){3pt}{10}
\cnode(4.5,0.5){3pt}{11} \cnode(4.5,2.3){3pt}{12}
\cnode[fillstyle=solid,fillcolor=lightgray](-.9,3.3){3pt}{13}\rput(-1.4,3.3){\scriptsize$x_{21}$}
\cnode[fillstyle=solid,fillcolor=lightgray](.9,3.3){3pt}{14}\rput(1.4,3.3){\scriptsize$x_{22}$}
\cnode(-.9,5.1){3pt}{15} \cnode(.9,5.1){3pt}{16}
\cnode(0,6.4){3pt}{17}

\rput(1.5,4.3){\scriptsize$H_2$}\rput(-3.7,2.6){\scriptsize$H_1$}\rput(3.7,2.7){\scriptsize$H_3$}

\ncline{1}{2} \ncline{1}{3}  \ncline{1}{4} \ncline{1}{5}
\ncline{2}{3} \ncline{2}{4} \ncline{2}{5} \ncline{3}{4}
\ncline{3}{5} \ncline{7}{4} \ncline{7}{5} \ncline{7}{6}
\ncline{7}{8} \ncline{7}{13} \ncline{7}{14} \ncline{8}{9}
\ncline{8}{10}  \ncline{8}{11} \ncline{8}{12} \ncline{9}{10}
\ncline{9}{11} \ncline{9}{12} \ncline{10}{11} \ncline{10}{12}
\ncline{11}{12} \ncline{17}{16} \ncline{17}{15}  \ncline{17}{14}
\ncline{17}{13} \ncline{16}{15} \ncline{16}{14} \ncline{16}{13}
\ncline{15}{14} \ncline{15}{13}

\end{pspicture}

\caption{\label{f3}\footnotesize The graph $G_2$}
\end{center}

\end{figure}

%\vskip6pt
%Example 2.3.
Let $t\geqslant 2$ be an integer. Samodivkin~\cite{s08} constructed
a very interesting graph $G_t$ to show that the upper bound in
Theorem~\ref{thm3.1.11} is better than the known bounds. Let
$H_0,H_1,H_2,\ldots,H_{t+1}$ be mutually vertex-disjoint graphs such
that $H_0\cong K_2$, $H_{t+1}\cong K_{t+3}$ and $H_i\cong K_{t+3}-e$
for each $i=1, 2,\ldots,t$. Let $V(H_0)=\{x,y\}$, $x_{t+1}\in
V(H_{t+1})$ and $x_{i1}, x_{i2}\in V(H_i), x_{i1}x_{i2}\notin
E(H_i)$ for each $i=1,2,\dots,t$. The graph $G_t$ is defined as
follows.
 $$
 \begin{array}{rl}
 & V(G_t)=\cup_{k=0}^{t+1}V(H_k)\ \ {\rm and}\\
 & E(G_t)=(\cup^{t+1}_{k=0}E(H_k))\cup(\cup^t_{i=1}\{xx_{i1},
 xx_{i2}\})\cup\{xx_{t+1}\}.
 \end{array}
 $$
Such a constructed graph $G_t$ is shown in Figure~\ref{f3} when
$t=2$.

Observe that $\gamma(G_t)=t+2$, $A(G_t)=V(G_t)$, $d_{G_t}(x)=2t+2$,
$d_{G_t}(x_{t+1})=t+3$, $d_{G_t}(y)=1$
%, $\lambda(G_t)=1$
and $d_{G_t}(z)=t+2$ for each $z\in V(G_t-\{x,y,x_{t+1}\})$.
Moreover, $\gamma(G-y)<\gamma(G)$ and $\gamma(G_t-z) =\gamma(G_t)$
for any $z\in V(G_t)-\{y\}$. Hence each of the bounds stated in
theorems \ref{thm3.1.2} - \ref{thm3.1.9} is greater than or equals
$t+2$.

Consider the graph $G_t-xy$. Clearly $\gamma(G_t-xy)=\gamma(G_t)$
and
 $$
 B(G_t-xy)=B(G_t-y)=\{x\}\cup V(H_{t+1}-x_{t+1})\cup(\cup^t_
{k=1}\{x_{k_1}, x_{k_2}\}).
 $$
Therefore, $N_{G_t}(x)\cap G(G_t-y)= \{x_{t+1}\}$ which implies that
the upper bound stated in Theorem~\ref{thm3.1.11}
%(iii)
is equals to
$d_{G_t}(y)+|\{x_{t+1}\}|=2$. Clearly $b(G_t)=2$ and hence this
bound is sharp for $G_t$.

\vskip6pt

From the graph $G_t$, we obtain the following statement immediately.

\begin{prop}\label{prop3.1.13}
For every integer $t\geqslant 2$, there is a graph $G$ such that the
difference between any upper bound stated in theorems \ref{thm3.1.2}
- \ref{thm3.1.9} and the upper bound of Theorem~\ref{thm3.1.11}
%, provided $G=G_t$,
is equal to $t$.
\end{prop}

\noindent{\bf Comments}\ Although Theorem~\ref{thm3.1.11} supplies
us with the upper bound that is closer to $b(G)$ for some graph $G$
than what any one of theorems \ref{thm3.1.2} - \ref{thm3.1.9}
provides, it is not easy to determine the sets $A^+(G)$ and $B(G)$
mentioned in Theorem~\ref{thm3.1.11} for an arbitrary graph $G$.
Thus the upper bound given in Theorem~\ref{thm3.1.11} is of
theoretical importance, but not applied since, until now, we have
not found a new class of graphs whose bondage numbers are determined
by Theorem~\ref{thm3.1.11}.

\vskip6pt

The above-mentioned upper bounds on the bondage number are involved
in only degrees of two vertices. Hartnell and Rall~\cite{hr99}
established an upper bound of $b(G)$ in terms of the numbers of
vertices and edges of $G$. For any connected graph $G$, let $\bar
\delta(G)$ represent the average degree of vertices in $G$. Hartnell
and Rall first discovered the following proposition.

\begin{prop}\label{prop3.1.14}
For any connected graph $G$, there exist two vertices $x$ and $y$
with distance at most two and, with the property that
$d_G(x)+d_G(y)\leqslant 2 \bar \delta(G)$.
\end{prop}

Using Proposition~\ref{prop3.1.14} and Theorem~\ref{thm3.1.5},
Hartnell and Rall gave the following bound.

\begin{thm}\label{thm3.1.15}\textnormal{(Hartnell and Rall~\cite{hr99}, 1999)}
For any connected graph $G$, with $n$ vertices and $m$ edges,
$b(G)\leqslant \frac{4m}{n}-1$.
\end{thm}

\begin{pf}
Let $G$ be a graph satisfying the hypothesis. By
Proposition~\ref{prop3.1.14}, there are two vertices $x$ and $y$
with distance at most two and, with the property that
$d_G(x)+d_G(y)\leqslant 2 \bar \delta(G)$. By
Theorem~\ref{thm3.1.5}, we have that
 $$
 b(G)+1 \leqslant d_G(x)+d_G(y)\leqslant 2\bar\delta(G),
 $$
from which, we have that
 $$
 4m(G)=2\, n\,\bar\delta(G)\geqslant n\, (b(G)+1).
 $$
That is, $b(G)\leqslant \frac{4m}{n}-1$.
\end{pf}

\begin{cor}
$b(G)\leqslant 2\bar\delta(G)-1$ for any connected graph $G$.
\end{cor}

\begin{cor}\label{cor3.1.17}
$m(G)\geqslant \frac{n}{4}(b(G)+1)$ for any connected graph $G$ of
order $n$.
\end{cor}

Hartnell and Rall~\cite{hr99} gave examples to show that for each
value of $b(G)$, the lower bound given in the
Corollary~\ref{cor3.1.17} is sharp for some values of $n$.

If $b(G)=1$, simply take $n=2$ (necessary for $G$ to be connected)
and $G$ isomorphic to $K_2$. If $b(G)=2$, consider $n=4$ and $G$
isomorphic to $P_4$.

For $b(G)=k$ with $2<k<\frac n2$, let $G$ be the graph on $n=4m$
vertices constructed as follows. Start with a $k$-graph $H$ with
order $2m$. In fact, if $k$ is even, then let $H$ be the circulant
graph $G(2m; \pm \{1,\ldots,\lfloor\frac k2\rfloor\})$; if $k$ is
odd, then let $H$ be $G(2m; \pm \{1,\ldots,\lfloor\frac
k2\rfloor\},\frac n2)$. Observe that each vertex is of degree $k$.
Now attach a leaf to each of the $2m$ vertices of $H$ to form $G$,
that is, $G=H\circ K_1$ with $b(G)=k+1$ (see Theorem~\ref{thm5.1}).

\vskip20pt
%\newpage

\subsection{Bounds Implied by Connectivity}

Use $\kappa(G)$ and $\lambda(G)$ to denote the vertex-connectivity
and the edge-connectivity of a connected graph $G$, respectively,
which are the minimum numbers of vertices and edges whose removal
result in $G$ disconnected. The famous Whitney's inequality can be
stated as $\kappa(G)\leqslant\lambda(G)\leqslant\delta(G)$ for any
graph or digraph $G$. Corollary~\ref{con3.1.4} was improved by
several authors as follows.

\begin{thm}\label{thm3.2.1} \textnormal{(Hartnell and
Rall~\cite{hr94}, 1994 and Teschner~\cite{t97}, 1997)} If $G$ is a
connected graph, then $b(G)\leqslant\Delta(G)+\lambda(G)-1$.
\end{thm}                                                   % Theorem 2.4

\begin{pf}
Let $G$ be a connected graph with edge-connectivity $\lambda(G)$ and
$F$ be $\lambda$-cut of $G$. Then $H=G-F$ is a spanning subgraph of
$G$.

If there is a vertex $x$ incident with some edge in $F$ such that
$\gamma(H-x)\geqslant\gamma(H)$, then $b(H)\leqslant d_H(x)$ by
Theorem~\ref{thm3.1.2}. By Lemma~\ref{lem3.1.1}, we have
$b(G)\leqslant d_H(x)+|F|=\Delta(G)+\lambda(G)-1$.

Assume now that any $\gamma$-set of $H$ contains all vertices
incident with edges in $F$. Arbitrarily choose an edge $xy\in F$.
Then there exists a vertex $z\in N_H(x)\setminus\{y\}$ that is
dominated only by $x$. Thus $\gamma(H-E'_x)>\gamma(H)$, where $E'_x$
is the set of edges incident with $x$ in $H$, and so $b(H)\leqslant
|E'_x|$. By Lemma~\ref{lem3.1.1}, we have $b(G)\leqslant
b(H)+|F|\leqslant d_H(x)+|F|= \Delta(G)+\lambda(G)-1$.
\end{pf}

\vskip6pt The upper bound given in Theorem~\ref{thm3.2.1} can be
attained. For example, a cycle $C_{3k+1}$ of order $3k+1$ with
$k\geqslant 1$, $b(C_{3k+1})=3$ by Theorem~\ref{thm2.1.1}. Since
$\kappa(C_{3k+1})=\lambda(C_{3k+1})=2$, we have
$\kappa(C_{3k+1})+\lambda(C_{3k+1})-1=2+2-1=3= b(C_{3k+1})$.

Motivated by Corollary~\ref{con3.1.4}, Theorems~\ref{thm3.2.1} and
the Whitney's inequality, Dunbar {\it et al.}~\cite{dhtv98}
naturally proposed the following conjecture.
%If $G$ is a connected nonempty graph, then
%$b(G)\leqslant\Delta(G)+\kappa(G)-1$.

\begin{con}\label{con3.2.2}
%\textnormal{(Dunbar, Haynes, Teschner and Volkmann~\cite{dhtv98}, 1998)}
If $G$ is a connected graph, then
$b(G)\leqslant\Delta(G)+\kappa(G)-1$.
\end{con}                                                % Conjecture 2.5

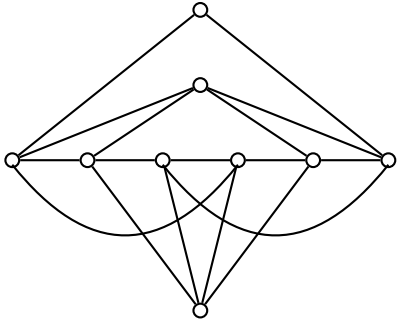
\begin{figure}[h]
\begin{pspicture}(-.3,1.7)(15,6.)
  \cnode(5,4){3pt}{v1}\cnode(6,4){3pt}{v2}\cnode(7,4){3pt}{v3}
  \cnode(8,4){3pt}{w3}\cnode(9,4){3pt}{w2}\cnode(10,4){3pt}{w1}
  \cnode(7.5,6){3pt}{u1}\cnode(7.5,5){3pt}{u2}\cnode(7.5,2){3pt}{u3}
  %\rput(7.5,6.5){\footnotesize the vertex of degree two}
  \ncline{v1}{v2}\ncline{w2}{w1}\ncline{w2}{w3}\ncline{w3}{v3}\ncline{v2}{v3}
  \ncline{u1}{v1}\ncline{u1}{w1}
  \ncline{u2}{v1}\ncline{u2}{w1}\ncline{u2}{v2}\ncline{u2}{w2}
  \ncline{u3}{v2}\ncline{u3}{v3}\ncline{u3}{w2}\ncline{u3}{w3}
  \parabola{-}(5,3.94)(6.5,3)\parabola{-}(10,3.94)(8.5,3)
\end{pspicture}
\caption{\label{f4}\footnotesize A graph $H$ with $\gamma(H)=3$ and
$b(H)=5$}
\end{figure}

However, Liu and Sun\,\cite{liu03} presented a counterexample to
this conjecture. They first constructed a graph $H$ showed in
Figure~\ref{f4} with $\gamma(H)=3$ and $b(H)=5$. Then let $G$ be the
disjoint union of two copies of $H$ by identifying two vertices of
degree two. They proved $b(G)\geqslant5$. Clearly, $G$ is a
$4$-regular graph with $\kappa(G)=1$ and $\lambda(G)=2$, and so
$b(G)\leqslant 5$ by Theorem~\ref{thm3.2.1}. Thus,
$b(G)=5>4=\Delta(G)+\kappa(G)-1$.

With a suspicion of the relationship between the bondage number and
the vertex-connectivity of a graph, the following conjecture is
proposed.

\begin{con}\label{con3.2.3}
\textnormal{(Liu and Sun\,\cite{liu03}, 2003)} For any positive
integer $r$, there exists a connected graph $G$ such that
$b(G)\geqslant\Delta(G)+\kappa(G)+r$.
\end{con}

%{\bf(Question: What graphs make the gap between $\kappa(G)$ and
%$\lambda(G)$ arbitrarily large?)}

To the knowledge of the author, until now no results have been known
about this conjecture.

\vskip6pt\noindent{\bf Comments}\ We conclude this subsection with
following comments.

From Theorem~\ref{thm3.2.1}, if Conjecture~\ref{con3.2.3} holds for
some connected graph $G$, then $\lambda(G)>\kappa(G)+r$, which
implies that $G$ is of large edge-connectivity and small
vertex-connectivity.

Use $\xi(G)$ to denote the minimum edge-degree of $G$, that is,
 $$
 \xi(G)=\min\limits_{xy\in E(G)}\{d_G(x)+d_G(y)-2\}.
 $$
Theorem~\ref{thm3.1.3} implies the following result.

\begin{prop}\label{prop3.2.4} \
$b(G)\leqslant \xi(G)+1$ for any graph $G$.
\end{prop}

Use $\lambda'(G)$ to denote the restricted edge-connectivity of a
connected graph $G$, which is the minimum number of edges whose
removal result in $G$ disconnected and no isolated vertices.
Esfahanian and Hakimi [A. H. Esfahanian and S. L. Hakimi, On
computing a conditional edge-connectivity of a graph. Information
Processing Letters, 27 (1988), 195-199] showed the following
result.

\begin{prop}\label{prop3.2.5}
If $G$ is neither $K_{1,n}$ nor $K_3$, then $\lambda'(G)\leqslant
\xi(G)$.
\end{prop}

Combining Proposition~\ref{prop3.2.4} and
Proposition~\ref{prop3.2.5}, we propose a conjecture as follows.

\begin{con}\label{con3.2.6}\ If a connected $G$ is neither $K_{1,n}$ nor $K_3$, then
$b(G)\leqslant \delta(G)+\lambda'(G)-1$.
\end{con}

A cycle $C_{3k+1}$ also satisfies
$b(C_{3k+1})=3=\delta(C_{3k+1})+\lambda'(C_{3k+1})-1$ since
$\lambda'(C_{3k+1})=\delta(C_{3k+1})=2$ for any integer $k\geqslant
1$.

For the graph $H$ shown in Figure~\ref{f4}, $\lambda'(H)=4$ and
$\delta(H)=2$, and so $b(H)=5=\delta(H)+\lambda'(H)-1$.

For the $4$-regulae graph $G$ constructed by Liu and
Sun\,\cite{liu03} obtained from the disjoint union of two copies of
the graph $H$ showed in Figure~\ref{f4} by identifying two vertices
of degree two, we have $b(G)\geqslant5$. Clearly, $\lambda'(G)=2$.
Thus, $b(G)=5=\delta(G)+\lambda'(G)-1$.

For the $4$-regulae graph $G_t$ constructed by
Samodivkin~\cite{s08}, see Figure~\ref{f3} for $t=2$, we have
$b(G_t)=2$. Clearly, $\delta(G_t)=1$ and $\lambda'(G)=2$. Thus,
$b(G)=2=\delta(G)+\lambda'(G)-1$.

These examples show that if Conjecture~\ref{con3.2.6} is true then
the given upper bound is tight.

\vskip20pt
%\newpage

\subsection{Bounds Implied by Degree Sequence}

Now let us return to Theorem~\ref{thm3.1.5}, from which
Teschner\cite{t97} obtained some other bounds in terms of the degree
sequences. The {\it degree sequence} $\pi(G)$ of a graph $G$ with
vertex-set $V=\{x_1,x_2,\ldots,x_n\}$ is the sequence
$\pi=(d_1,d_2,\ldots,d_n)$ with $d_1\leqslant
d_2\leqslant\cdots\leqslant d_n$, where $d_i=d_G(x_i)$ for each
$i=1,2,\ldots,n$. The following result is essentially a corollary of
Theorem~\ref{thm3.1.5}.

%For a positive integer $k$, use $\alpha_k(G)$ to
%denote the maximum number of vertices any distinct vertices of which
%has distance larger than $k$ in $G$.

%For a positive integer $k$, a set $I\subseteq (G)$ is called a {\it
%k-independent set} if $d_G(x,y)>k$ for any two distinct vertices $x$
%and $y$ in $I$. The maximum cardinality of all $k$-independent sets
%is called the {\it k-independence number} of $G$, denoted by
%$\alpha_k(G)$ (see, for example, Topp and Vestergaard~\cite{tv00}).

\begin{thm}\label{thm3.3.1} \textnormal{(Teschner~\cite{t97}, 1997)}
Let $G$ be a nonempty graph with degree sequence $\pi(G)$. If
$\alpha_2(G)=t$, then $b(G)\leqslant d_t+d_{t+1}-1$.
\end{thm}                                                     % Lemma 2.7

\begin{pf}
Let $I=\{x_1,x_2,\ldots, x_t,x_{t+1}\}$, the set of vertices
corresponding the first $t+1$ elements in the degree sequence
$\pi(G)$. If there are two vertices $x$ and $y$ in $I$ with
$d_G(x,y)\leqslant 2$, then the lemma follows by
Theorem~\ref{thm3.1.5} immediately. Otherwise, $I$ is a
$2$-independent set, and so $\alpha_2(G)\geqslant |I|=t+1$, a
contradiction.
\end{pf}

\vskip6pt

Combining Theorem~\ref{thm3.3.1} with Proposition~\ref{prop2.2.5}
(that is, $\alpha_2(G)=\gamma(G)$), we have the following corollary.

\begin{cor}\label{cor3.3.2} \textnormal{(Teschner~\cite{t97}, 1997)}
Let $G$ be a nonempty graph with the degree sequence $\pi(G)$. If
$\gamma(G)=\gamma$, then $b(G)\leqslant d_{\gamma}+d_{\gamma+1}-1$.
\end{cor}                                                 % Corollary 2.8

In \cite{t97}, Teschner showed that these two bounds are sharp for
arbitrarily many graphs. Let $H=C_{3k+1}+\{x_1x_4,x_1x_{3k-1}\}$,
where $C_{3k+1}$ is a cycle $(x_1,x_2,\ldots,x_{3k+1},x_1)$ for any
integer $k\geqslant 2$. Then $\gamma(H)=k+1$ and so $b(H)\leqslant
2+2-1=3$ by Corollary~\ref{cor3.3.2}. Since $C_{3k+1}$ is a spanning
subgraph of $H$ and $\gamma(G)=\gamma(H)$, Lemma \ref{lem3.1.1}
yields that $b(H)\geqslant b(C_{3k+1})=3$. Then $b(H)=3$.
%$b(H)\leqslant b(G)$ ¼ûLemma \ref{lem3.1.1}ϵĺì×Ö½áÂÛ.

\vskip6pt

Hartnell and Rall~\cite{hr99} established an upper bound of the
bondage number $b(G)$ in terms of order and the sum of all degrees.
For a connected graph $G$ with order $\upsilon$, let $\mu(G)=\frac
1\upsilon\sum\limits_{x\in V}d_G(x)$, called the average degree of
$G$.

\vskip6pt\begin{prop}\label{prop3.3.3} \textnormal{(Hartnell and
Rall~\cite{hr99}, 1999)}\ For a connected graph $G$, there exists a
pair of vertices $x$ and $y$ such that $d_G(x,y)\leqslant 2$ and
$d_G(x)+d_G(y)\leqslant 2\,\mu(G)$.
\end{prop}

\begin{pf}\ Assume that there is a connected graph $G$ such that
the proposition is false. Let $X=\{x\in V(G):\
d_G(x)\leqslant\mu(G)\}$ and $Y=\{x\in V(G):\ d_G(y)>\mu(G)\}$.

By our assumption, $X$ is an independent set in $G$. Hence, each
$x\in X$ has only vertices in $Y$ as its neighbors. Also each $y\in
Y$ has at most one vertex of $X$ as its neighbor otherwise, if there
were two, they would contradict our assumption. These facts imply
that $G$ has a matching $M$ that saturate every vertex in $X$ and
$|X|=|Y|$.

By our assumption, $d_G(x)+d_G(y)>2\, \mu(G)$ for every $xy\in M$,
where $x\in X$ and $y\in Y$. Thus,
 $$
 \sum\limits_{v\in V(G)}d_G(u)=\sum\limits_{xy\in
 M}(d_G(x)+d_G(y))>2\,|X|\mu(G)=\upsilon(G)\mu(G)=\sum\limits_{v\in V(G)}d_G(u),
 $$
a contradiction. The proposition follows.
\end{pf}

Combining Proposition~\ref{prop3.3.3} with Theorem~\ref{thm3.1.5},
they obtained the following result.

\begin{thm}\label{thm3.3.4} \textnormal{(Hartnell and Rall~\cite{hr99},
1999)} Let $G$ be a connected graph. Then the bondage number
 $
 b(G)\leqslant 2\,\mu(G)-1.
 $
\end{thm}

\begin{pf}
Let $G$ be a connected graph with order $\upsilon$ and average
degree $\mu$. By Proposition~\ref{prop3.3.3}, there is such two
vertices, say $x$ and $y$, that $d_G(x,y)\leqslant 2$ and
$d_G(x)+d_G(y)\leqslant 2\,\mu(G)$. By Theorem~\ref{thm3.1.5}, we
immediately have
 $$
 b(G)\leqslant d_G(x)+d_G(y)-1\leqslant 2\,\mu(G)-1.
 $$
The theorem follows.
\end{pf}

Note that the number of edges $2\,\varepsilon(G)=\upsilon\,\mu(G)$.
Theorem~\ref{thm3.3.4} implies the following bound in terms of
vertex-number $\upsilon(G)$ and edge-number $\varepsilon(G)$.
 \begin{equation}\label{e2.1}
 \begin{array}{rl}
 b(G)\leqslant \frac{4\,\varepsilon(G)}{\upsilon(G)}-1.
 \end{array}
 \end{equation}

Hartnell and Rall~\cite{hr99} also observed that for each value of
$b(G)$, the upper bound given in Eq.~(\ref{e2.1}) is sharp for some
values of $\upsilon$.

If $b(G)=1$ or $2$, simply take $G=K_2$ or $P_4$, respectively. For
$b(G)=k>2$, let $H$ be a circulant undirected graph with order $2m$
and degree $k-1$, and let $G$ be a graph obtained from $H$ by
attaching a leaf to each of the $2m$ vertices. It is easy to see
that $\varepsilon(G)=m(k+1)$ and $b(G)=k$.

\vskip6pt

\noindent{\bf Comments}\ Although various of upper bounds have been
establish as the above, we find that the appearance of these bounds
is essentially based upon the local structures of a graph, precisely
speaking, the structures of the neighborhoods of two vertices within
distance $2$. Even if these bounds can be achieved by some special
graphs, it is more often not the case. The reason lies essentially
in the definition of the bondage number, which is the minimum value
among all bondage sets, an integral property of a graph. While it
easy to find upper bounds just by choosing some bondage set, the gap
between the exact value of the bondage number and such a bound
obtained only from local structures of a graph is often large. For
example, a star $K_{1,\Delta}$, however large $\Delta$ is,
$b(K_{1,\Delta})=1$. Therefore one has been longing for better
bounds upon some integral parameters. However, as what we will see
below, it is difficult to establish such upper bounds.

\subsection{Bounds in $\gamma$-critical Graphs}

A graph $G$ is called a {\it vertex domination-critical
graph}\,(\,{\it vc-graph} or {\it $\gamma$-critical} for short) if
$\gamma(G-x)<\gamma(G)$ for any vertex $x\in V(G)$, proposed by
Brigham, Chinn and Dutton~\cite{bcd88} in 1988.

Several families of graphs are known to be $\gamma$-critical. From
definition, it is clear that if $G$ is a $\gamma$-critical graph,
then $\gamma(G)\geqslant 2$. The class of $\gamma$-critical graphs
with $\gamma=2$ is characterized as follows.

\begin{prop}\label{prop3.4.1}
\textnormal{(Brigham, Chinn and Dutton~\cite{bcd88}, 1988)} A graph
$G$ with $\gamma(G)=2$ is a $\gamma$-critical graph if and only if
$G$ is a complete graph $K_{2t}$ ($t\geqslant2$) with a perfect
matching removed.
\end{prop}

A more interesting family is composed of the $n$-critical graphs
$G_{m,n}$ defined for $m, n\geqslant 2$ by the circulant undirected
graph $G(N,\pm S)$, where $N=(m+1)(n-1)+1$ and
$S=\{1,2,\ldots,\lfloor m/2\rfloor\}$.

The reason why the $\gamma$-critical graphs are of special interest
in this context is easy to see that they play an important role in
the study of the bondage number. For instance, it immediately
follows from Theorem~\ref{thm3.1.2} that if $b(G)>\Delta(G)$ then
$G$ is a $\gamma$-critical graph. The $\gamma$-critical graphs are
defined exactly in this way. In order to find graphs $G$ with a high
bondage number (i.e. higher than $\Delta(G))$ and beyond its general
upper bounds for the bondage number we therefore have to look at
$\gamma$-critical graphs.

The bondage numbers of some $\gamma$-critical graphs have been
examined by several authors, see for example~\cite{s08, s09, ts95,
t97}. From Theorem~\ref{thm3.1.2} we know that the bondage number of
a graph $G$ is bounded from above by $\Delta(G)$ if $G$ is not a
$\gamma$-critical graph. For $\gamma$-critical graphs it is more
difficult to find an upper bound. We will see that the bondage
numbers of $\gamma$-critical graphs in general are not even bounded
from above by $\Delta+c$ for any fixed natural number $c$.

%There are some basic properties for $\gamma$-critical graphs.

%If a graph $G$ is a counterexample to
%Conjecture~\ref{con3.15}, it must be a $\gamma$-critical graph by
%Theorem~\ref{thm3.1.2}. It is why the vertex domination-critical
%graphs are of special interest in the literature.
%It is that this fact motivates one's special interest in $\gamma$-critical graphs.

In this subsection we introduce some upper bounds for the bondage
number of a $\gamma$-critical graph. By Proposition~\ref{prop3.4.1},
we can easily prove the following result.

\begin{thm}\label{thm3.4.2}
If $G$ is a $\gamma$-critical graph with $\gamma(G)=2$, then
$b(G)\leqslant \Delta+1$.
\end{thm}

\begin{pf}
By Proposition~\ref{prop3.4.1}, $G$ is a complete $K_{2t}$
($t\geqslant2$) with a perfect matching removed. Since
$\gamma(G)=2$, let $\{x,y\}$ be a $\gamma$-set of $G$. Furthermore,
Let $H=G-E_x$. Since $G$ is a $\gamma$-critical graph, $\{x,y\}$ is
the unique $\gamma$-set of $H$, and so $b(H)=1$ by
Theorem~\ref{thm2.2}. By Lemma~\ref{lem3.1.1}, we have
$b(G)\leqslant b(H)+|E_x|=1+\Delta(G)$.
\end{pf}

\vskip6pt In Section 4, by Theorem~\ref{thm4.6}, we will see the
equality in Theorem~\ref{thm3.4.2} holds, that is, $b(G)=\Delta+1$
if $G$ is a $\gamma$-critical graph with $\gamma(G)=2$.

%In 1997, Teschner \cite{t97} obtained an upper bound of the bondage
%number for a $\gamma$-critical graph.

\vskip6pt \begin{thm}\label{thm3.4.3}
\textnormal{(Teschner~\cite{t97}, 1997)} Let $G$ be a
$\gamma$-critical graph with degree sequence $\pi(G)$. Then
$b(G)\leqslant\max\{\Delta(G)+1,d_1+d_2+\cdots+d_{\gamma-1}\}$,
where $\gamma=\gamma(G)$.
\end{thm}                                                   % Theorem 3.19

\begin{pf}
Since $G$ is a $\gamma$-critical graph, we have $\gamma(G)\geqslant
2$. Let $t=d_1+d_2+\cdots+d_{\gamma-1}$. To prove the theorem, it is
sufficient to show that if $b(G)>t(G)$ then $b(G)\leqslant
\Delta(G)+1$.

Let $E_i=E_{x_i}(G)$, where $d_G(x_i)=d_i$ for each
$i=1,2,\ldots,\gamma-1$. Let $F=E_1\cup E_2\cup\ldots\cup
E_{\gamma-1}$ and let $U=\{x_1,x_2,\ldots,x_{\gamma-1}\}$. Since
$|F|\leqslant t$, we see by the hypothesis that
$\gamma(G-F)=\gamma(G)$, that is, $U$ is contained in any
$\gamma$-set of $G$, and so $\gamma(G-U)=1$ since $|U|=\gamma-1$.
Hence, we deduce $\Delta(G)\geqslant\upsilon(G)-\gamma(G)$. From the
obvious inequality $\Delta(G)\leqslant\upsilon(G)-\gamma(G)$ we
conclude $\Delta(G)=\upsilon(G)-\gamma(G)$.

Let $x$ be a vertex of maximum degree $\Delta(G)$.
%Then we know by Brigham, Chinn and Dutton~\cite{bcd88} (Theorem 7) that
%$G$ has certain properties, among others the property, that
It is easy to see that $G-x$ has a unique minimum dominating set.
Let $H=G-E_x$. Then $H$ must also have a unique minimum dominating
set. Then by Theorem~\ref{thm2.2} we have $b(H)=1$, and so
$b(G)\leqslant\Delta(G)+1$.
\end{pf}

\vskip6pt

As we mention above, if $G$ is a $\gamma$-critical graph with
$\gamma(G)=2$ then $b(G)=\Delta+1$, which shows that the bound given
in Theorem~\ref{thm3.4.3} can be attained for $\gamma=2$. However,
we have not known whether this bound is tight for general
$\gamma\geqslant 3$. Theorem~\ref{thm3.4.3}
%generalizes the second conclusion in Theorem \ref{thm3.3.1} (b) for a nonempty graph, and
gives the following
corollary immediately.

\begin{cor}\label{cor3.4.4}\textnormal{(Teschner~\cite{t97}, 1997)}
Let $G$ be a $\gamma$-critical graph with degree sequence $\pi(G)$.
If $\gamma(G)=3$, then
$b(G)\leqslant\max\{\Delta(G)+1,\delta(G)+d_2\}$.
\end{cor}

From Theorem~\ref{thm3.4.2} and Corollary~\ref{cor3.4.4}, we have
$b(G)\leqslant 2\,\Delta(G)$ if $G$ is a $\gamma$-critical graph
graph with $\gamma(G)\leqslant 3$. The following result shows that
this bound is not tight.

%The following result shows that Conjecture~\ref{con3.19} is valid
%for some $\gamma$-critical graphs.

\begin{thm}\label{thm3.4.5} \textnormal{(Teschner~\cite{ts95}, 1995)}
Let $G$ be a $\gamma$-critical graph graph with $\gamma(G)\leqslant
3$. Then $b(G)\leqslant\frac{3}{2}\,\Delta(G)$.
\end{thm}

Until now, we have not known whether the bound given in
Theorem~\ref{thm3.4.5} is tight or not. We state two conjectures on
$\gamma$-critical graphs proposed by Samodivkin~\cite{s08}. The
first of them was motivated by Theorem~\ref{thm3.1.11} and
Theorem~\ref{thm3.1.8}.

\begin{con}\label{con3.4.6}
\textnormal{(Samodivkin~\cite{s08}, 2008)} For every connected
nontrivial $\gamma$-critical graph $G$,
$$
\begin{array}{rl}
 \min\limits_{x\in A^+(G),\,y\in B(G-x)}\{d_G(x)+
|N_G(y)\cap A(G-x)|\}\leqslant\frac 32\Delta(G).
\end{array}
$$
\end{con}

To state the second conjecture, we need the following result on
$\gamma$-critical graphs.

\begin{prop}\label{prop3.4.7}
%\textnormal{(Brigham, Chinn and Dutton~\cite{bcd88}, 1988; Fulman,
%Hanson and Mac- Gillivray~\cite{fhm95}, 1995)}
If $G$ is a $\gamma$-critical graph then $\upsilon(G)\leqslant
(\Delta(G)+1)(\gamma(G)-1)+1$, moreover, if the equality holds, then
$G$ is regular.
\end{prop}

The bound of Proposition~\ref{prop3.4.7} is the best possible in the
sense that equality holds for the infinite class of
$\gamma$-critical graphs $G_{m,n}$ defined in the beginning of this
subsection. In Proposition~\ref{prop3.4.7}, the first result is due
to Brigham, Chinn and Dutton~\cite{bcd88} in 1988; the second is due
to Fulman, Hanson and MacGillivray~\cite{fhm95} in 1995.

For a $\gamma$-critical graph $G$ with $\gamma=3$, by
Proposition~\ref{prop3.4.7}, $\upsilon(G)\leqslant 2\Delta(G)+3$.

\begin{thm}\label{thm3.4.8}
\textnormal{(Teschner~\cite{ts95}, 1995)} If $G$ is a
$\gamma$-critical graph with $\gamma=3$ and $\upsilon(G)=
2\Delta(G)+3$, then $b(G)\leqslant\Delta +1$.
\end{thm}

We have not yet known whether the equality in Theorem~\ref{thm3.4.8}
holds or not. However, Samodivkin proposed the following conjecture.

%Conjecture 4.4.
\begin{con}\label{con3.4.9}
\textnormal{(Samodivkin~\cite{s08}, 2008)} If $G$ is a
$\gamma$-critical graph with $(\Delta(G)+1)(\gamma-1)+1$ vertices
then $b(G)=\Delta(G)+1$.
\end{con}

%V. Samodivkin~\cite{s08, s09}

In general, based on Theorem~\ref{thm3.4.5}, Teschner proposed the
following conjecture.

\begin{con}\label{con3.4.10}
\textnormal{(Teschner~\cite{ts95}, 1995)} If $G$ is a
$\gamma$-critical graph then $b(G)\leqslant\frac 32\,\Delta(G)$ for
$\gamma\geqslant 4$.
\end{con}

\noindent{\bf Comments}\ We conclude this subsection with some
comments.

Graphs which are minimal or critical with respect to a given
property or parameter frequently play an important role in the
investigation of that property or parameter. Not only are such
graphs of considerable interest in their own right, but also a
knowledge of their structure often aids in the development of the
general theory. In particular, when investigating any finite
structure, a great number of results are proven by induction.
Consequently it is desirable to learn as much as possible about
those graphs that are critical with respect to a given property or
parameter so as to aid and abet such investigations.
%(abet vt.½ÌËô, É¿¶¯, °ïÖú, Ö§³Ö)

In this subsection we survey some results on the bondage number for
$\gamma$-critical graphs. Although these results are not very
perfect, it provides a feasible method to approach the bondage
number from different angles. In particular, the methods given in
Teschner~\cite{ts95} worthily further explore and develop.
%{\color{red}?????????}

The following proposition is maybe useful for us to further
investigate the bondage number of a $\gamma$-critical graph.

\begin{prop}\label{prop3.4.11}
\textnormal{(Brigham, Chinn and Dutton~\cite{bcd88}, 1988)} If $G$
has a non-isolated vertex $x$ such that the induced subgraph
$G[N_G(x)]$ is complete, then $G$ is not $\gamma$-critical.
\end{prop}

\begin{pf}  Let $y\in N_G(x)$. Any minimum dominating set of $G-y$
includes a vertex of $N_G[v]$ and hence must also be a minimum
dominating set of $G$.
\end{pf}

\vskip6pt

This simple fact shows that any $\gamma$-critical graph contains no
vertices of degree one.

\vskip20pt

\subsection{Bounds Implied by Domination}

In preceding subsection, we introduce some upper bounds of the
bondage numbers for $\gamma$-critical graphs by consideration of
dominations. In this subsection, we introduce related results for
general graphs.

\begin{thm}\label{thm3.5.1}
\textnormal{(Fink  {\it et al.}~\cite{fjkr90}, 1990)} Let $G$ be a
connected graph of order $\upsilon\,(\geqslant2)$. Then

1) $b(G)\leqslant \upsilon-1$;

2) $b(G)\leqslant(\gamma(G)-1)\Delta(G)+1$ if $\gamma(G)\geqslant2$;

3) $b(G)\leqslant \upsilon-\gamma(G)+1$.
\end{thm}                                                   % Theorem 2.11

\begin{pf}
Let $x$ and $y$ be adjacent vertices with $d_G(x)\leqslant d_G(y)$.
If $b(G)\leqslant d_G(x)$, then $b(G)\leqslant\upsilon-1$, so we
suppose that $b(G)>d_G(x)$. Then $\gamma(G-E_x)=\gamma(G)$ and
$\gamma(G-x)=\gamma(G)-1$. Also, if $D$ denotes the union of all
$\gamma$-sets in $G-x$, then $x$ is in $G$ not adjacent any vertex
in $D$. Hence, $|E_x|\leqslant \upsilon-1-|D|$ and $y\notin D$. Now,
if $F_y$ denotes the set of edges from $y$ to a vertex in $D$, then
since $y\notin D$ we must have
$\gamma(G-x-F_y)>\gamma(G-x)=\gamma(G)-1$. Thus $\gamma(G-(E_x\cup
F_y))>\gamma(G)$ and we see that $b(G)\leqslant|E_x\cup
F_y|=|E_x|+|F_y|\leqslant(\upsilon-1-|D|)+|D|=\upsilon-1$. This
proves the first assertion. The proofs of other two assertions are
omitted here, and left the reader as an exercise.
\end{pf}

\vskip6pt

While the bound $b(G)\leqslant\upsilon -1$ is not particularly good
for many classes of graphs (e.g. trees and most cycles), it is an
attainable bound. For example, if $G$ is a complete $t$-partite
graph $G=K_{2,2,\ldots,2}$, then the three bounds on $b(G)$ in
Theorem~\ref{thm3.5.1} are sharp.

Teschner~\cite{ts95} consider $\gamma(G)=1$ and $\gamma(G)=2$. The
next result is almost trivial but useful, the proof is similar to
the proof of Theorem~\ref{thm2.1.1} (a).

\begin{lem}\label{lem3.5.2}
Let $G$ be a graph with order $n$ and $\gamma(G)=1$, $t$ be the
number of vertices of degree $n-1$. Then $b(G)=\lceil t/2\rceil$.
\end{lem}

Since $t\leqslant \Delta(G)$ clearly, Lemma~\ref{lem3.5.2} yields
the following result immediately.

\begin{thm}\label{thm3.5.3}
\textnormal{(Teschner~\cite{ts95}, 1995)} $b(G)\leqslant\frac
12\Delta +1$ for any graph $G$ with $\gamma(G)=1$.
\end{thm}

\begin{thm}\label{thm3.5.4}
\textnormal{(Teschner~\cite{ts95}, 1995)} $b(G)\leqslant\Delta+1$
for any graph $G$ with $\gamma(G)=2$.
\end{thm}

\begin{pf}
If there is a vertex $x$ in $G$ such that $\gamma(G-x)\geqslant
\gamma(G)$ then, by Theorem~\ref{thm3.1.2}, $b(G)\leqslant
\Delta(G)$. Suppose now that $G$ is a $\gamma$-critical graph. By
Theorem~\ref{thm3.4.2}, $b(G)\leqslant \Delta+1$.
\end{pf}

\vskip6pt

\noindent{\bf Remarks}\ For a complete graph $K_{2k+1}$,
$b(K_{2k+1})=k+1=\frac 12\Delta +1$, which shows that the upper
bound given in Theorem~\ref{thm3.5.3} can be attained. For a graph
$G$ with $\gamma(G)=2$, by Theorem~\ref{thm4.6} later, the upper
bound given in Theorem~\ref{thm3.5.4} can be also attained by a
$2$-critical graph (see Proposition~\ref{prop3.4.1}).

\vskip10pt

\subsection{Two Conjectures}

In 1990, when Fink {\it et al.}~\cite{fjkr90} introduced the concept
of the bondage number, they proposed the following conjecture.
%which has not been solved yet.

\begin{con}\label{con3.6.1}
%\textnormal{(Fink, Jacobson, Kinch and Roberts~\cite{fjkr90}, 1990)}
If $G$ is a nonempty graph, then $b(G)\leqslant\Delta(G)+1$.
\end{con}                                                % Conjecture 2.12

\begin{figure}[ht]
\begin{pspicture}(-1.65,1)(10,5)
  \cnode(3,1){3pt}{a}\cnode(9,1){3pt}{b}\cnode(9,5){3pt}{c}\cnode(3,5){3pt}{d}
  \cnode(6,2){3pt}{e}\cnode(7,3){3pt}{f}\cnode(6,4){3pt}{g}\cnode(5,3){3pt}{h}
  \cnode(6,3){3pt}{i}
  \ncline{a}{b}\ncline{b}{c}\ncline{c}{d}\ncline{d}{a}
  \ncline{i}{e}\ncline{i}{f}\ncline{i}{g}\ncline{i}{h}
  \ncline{e}{a}\ncline{e}{b}
  \ncline{f}{b}\ncline{f}{c}
  \ncline{g}{c}\ncline{g}{d}
  \ncline{h}{d}\ncline{h}{a}
  \ncline{a}{b}
  \psbezier{-}(5.95,2.05)(5.5,2.5)(5.5,3.5)(5.95,3.95)
  \psbezier{-}(7,3.06)(6.5,3.5)(5.5,3.5)(5,3.06)
  \end{pspicture}
\caption{\label{f5}\footnotesize A graph with $\Delta=4$ and $b=6$}
\end{figure}
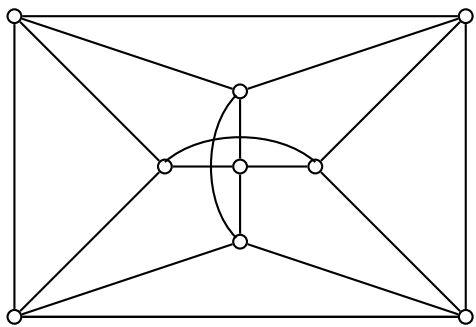

Although these results partially support Conjecture~\ref{con3.6.1},
Teschner \cite{ts93} in 1993 found a counterexample to this
conjecture, the cartesian product $K_3\times K_3$, as shown in
Figure~\ref{f5}, which shows $b(G)=\Delta(G)+2$.

\vskip6pt

If a graph $G$ is a counterexample to Conjecture~\ref{con3.6.1}, it
must be a $\gamma$-critical graph by Theorem~\ref{thm3.1.2}. It is
why the vertex domination-critical graphs are of special interest in
the literature.

\vskip6pt

Now we return to Conjecture~\ref{con3.6.1}. Hartnell and Rall
\cite{hr94} and Teschner \cite{ts96}, independently, proved that
$b(G)$ can be much greater than $\Delta(G)$ by showing the following
result.

\begin{thm}\label{thm3.6.2} \textnormal{(Hartnell and Rall
\cite{hr94}, 1994; Teschner \cite{ts96}, 1996)} For an integer
$n\geqslant 3$, let $G_n$ be the cartesian product $K_n\times K_n$.
Then $b(G_n)=3(n-1)=\frac 32 \Delta(G_n)$.
\end{thm}

This theorem shows that there exist no upper bounds of the form
$b(G)\leqslant\Delta(G)+c$ for any integer $c$. Teschner \cite{ts95}
proved that $b(G)\leqslant\frac{3}{2}\,\Delta(G)$ for any graph $G$
with $\gamma(G)\leqslant 2$ (see Theorem~\ref{thm3.5.3} and
Theorem~\ref{thm3.5.4}) and for some class of graphs with
$\gamma(G)=3$, and proposed the following conjecture.

\begin{con}\label{con3.6.3} \textnormal{(Teschner \cite{ts95}, 1995)}
$b(G)\leqslant\frac{3}{2}\,\Delta(G)$ for any graph $G$.
\end{con}                                                % Conjecture 2.15

We believe that this conjecture is valid, but so far there are no
much work about it.

\section{Lower Bounds}

Since the bondage number is defined as the smallest number of edges
whose removal results in increase of domination number, each
constructive method that creates a concrete bondage set leads to an
upper bound on the bondage number. For that reason it is hard to
find lower bounds. Nevertheless, there are still a few lower bounds
obtained by Teschner \cite{t97}, the first one of them can be got in
terms of its spanning subgraph.
%(In the next section we will
%mention the exact values of bondage numbers of some graphs.)

\begin{thm}\label{thm4.1} \textnormal{(Teschner~\cite{t97}, 1997)}
Let $H$ be a spanning subgraph of a nonempty graph $G$. If
$\gamma(H)=\gamma(G)$ then $b(H)\leqslant b(G)$.
\end{thm}

\begin{pf}
Let $B\subseteq E(G)$ be a bondage set of $G$ with $|B|=b(G)$. We
necessarily have $\gamma(H-B)>\gamma(G-K)>\gamma(G)=\gamma(H)$ and,
therefore, $b(H) \leqslant |B|=b(G)$.
\end{pf}

\vskip6pt

By Theorem~\ref{thm2.1.1}, $b(C_n)=3$ if $n\equiv 1\,({\rm mod}\,3)$
and $b(C_n)=2$ otherwise, $b(P_n)=2$ if $n\equiv 1\,({\rm mod}\,3)$
and $b(P_n)=1$ otherwise. From these results and
Theorem~\ref{thm4.1}, we get the following two corollaries.

\begin{cor}
If $G$ is hamiltonian with order $n\geqslant 3$ and
$\gamma(G)=\lceil n/3\rceil$, then $b(G)\geqslant 2$ and in addition
$b(G)\geqslant 3$ if $n\equiv 1\,({\rm mod}\,3)$.
\end{cor}

\begin{cor}
If $G$ of order $n\geqslant 2$ has a hamiltonian path and
$\gamma(G)=\lceil n/3\rceil$, then $b(G)\geqslant 2$ if $n\equiv
1\,({\rm mod}\,3)$.
\end{cor}

The {\it vertex covering number} $\beta(G)$ of $G$ is the minimum
number of vertices that are incident with all edges in $G$.
%For any graph $G$, a vertex set $V$ is called a {\it vertex
%covering} if every edge of $G$ is incident with some vertex in $V$.
%The minimum cardinality of all vertex coverings is the {\it vertex
%covering number}, denoted by $\beta(G)$.
If $G$ has no isolated vertices, then $\gamma(G)\leqslant\beta(G)$
clearly. In \cite{v94}, Volkmann gave a lot of graphs with
$\beta=\gamma$.
% since any vertex covering is a dominating set.

\begin{thm}\label{thm4.2} \textnormal{(Teschner~\cite{t97}, 1997)}
Let $G$ be a graph. If $\gamma(G)=\beta(G)$, then

1) $b(G)\geqslant\delta(G)$;

2) $b(G)\geqslant\delta(G)+1$ if $G$ is a $\gamma$-critical graph.
\end{thm}                                                   % Theorem 3.2

\begin{pf} Let $G$ be a graph with
$\beta(G)=\gamma(G)$.

1) Assume $\delta(G)\geqslant2$, without loss of generality. Let
$B\subseteq E(G)$ with $|B|\leqslant\delta(G)-1$. Then
$\delta(G-B)\geqslant1$ and so
$\gamma(G-B)\leqslant\beta(G-B)\leqslant\beta(G)=\gamma(G)$. Thus,
$B$ is not a bondage set of $G$, and so $b(G)\geqslant\delta(G)$.

2) It is clear from the proof of 1) that for any bondage set $B$,
$E_x\subseteq B$ for some vertex $x$, and so $x$ is an isolated
vertex in $G-B$. On the other hand, if $G$ is a $\gamma$-critical
graph, then $\gamma(G-x)<\gamma(G)$, which implies that
$\gamma(G-E_x)\leqslant\gamma(G)$. Thus $b(G)\geqslant\delta(G)+1$.
\end{pf}

\vskip6pt

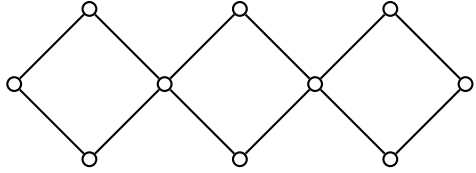
\begin{figure}[ht]
\begin{pspicture}(-2.5,1)(5,3.5)
  \cnode(3,1){3pt}{u1}\cnode(5,1){3pt}{u2}\cnode(7,1){3pt}{u3}
  \cnode(2,2){3pt}{v1}\cnode(4,2){3pt}{v2}\cnode(6,2){3pt}{v3}\cnode(8,2){3pt}{v4}
  \cnode(3,3){3pt}{w1}\cnode(5,3){3pt}{w2}\cnode(7,3){3pt}{w3}
  \ncline{v1}{u1}\ncline{v1}{w1}
  \ncline{v2}{u1}\ncline{v2}{w1}\ncline{v2}{u2}\ncline{v2}{w2}
  \ncline{v3}{u2}\ncline{v3}{w2}\ncline{v3}{u3}\ncline{v3}{w3}
  \ncline{v4}{u3}\ncline{v4}{w3}
  \end{pspicture}
\caption{\label{f6}\footnotesize A $\gamma$-critical graph with
$\beta=\gamma=4,\delta=2,b=3$}
\end{figure}

The graph shown in Figure~\ref{f6} shows that the bound given in
Theorem \ref{thm4.2} is sharp. For the graph $G$, it is a
$\gamma$-critical graph and $\gamma(G)=\beta(G)=4$. By
Theorem~\ref{thm4.2}, we have $b(G)\geqslant 3$. On the other hand,
$b(G)\leqslant 3$ by Theorem~\ref{thm3.1.5}. Thus, $b(G)=3$.

%\cnode(-3.5,2){3pt}{1}

\begin{prop}\label{prop4.3} \textnormal{(Sanchis~\cite{san91}, 1991)}
Let $G$ be a graph $G$ of order $\upsilon\,(\geqslant 6)$. If $G$
has no isolated vertices and $3\leqslant\gamma(G)\leqslant
\upsilon/2$, then
$\varepsilon(G)\leqslant{{\upsilon-\gamma(G)+1}\choose{2}}$.
\end{prop}                                                   % Lemma 3.3

Using the idea in the proof of Theorem~\ref{thm4.2}, every upper
bound for $\gamma(G)$ can lead to a lower bound for $b(G)$. In this
way Teschner \cite{t97} obtained another lower bound from
Proposition~\ref{prop4.3}.

\begin{thm}\label{thm4.6} \textnormal{(Teschner~\cite{t97}, 1997)}
Let $G$ be a graph $G$ of order $\upsilon\,(\geqslant 6)$, and
$2\leqslant\gamma(G)\leqslant \upsilon/2-1$. Then

1)
$b(G)\geqslant\min\{\delta(G),\varepsilon-{{\upsilon-\gamma(G)}\choose{2}}
\}$;

2)
$b(G)\geqslant\min\{\delta(G)+1,\varepsilon-{{\upsilon-\gamma(G)}\choose{2}}
\}$ if $G$ is a $\gamma$-critical graph.
\end{thm}                                                   % Theorem 3.4

\begin{pf} Let $G$ be a graph with
$2\leqslant\gamma(G)\leqslant \upsilon/2-1$ and $B$ be a minimum
bondage set. Then $\gamma(G-B)=\gamma(G)+1$ and so
$3\leqslant\gamma(G-B)\leqslant \upsilon/2$.

If $G-B$ has isolated vertices, then $b(G)\geqslant\delta(G)$ or
$b(G)\geqslant\delta(G)+1$ when $G$ is a $\gamma$-critical graph.

If $G-B$ has no isolated vertices, then Proposition~\ref{prop4.3}
yields that
 $$
 \begin{array}{rl}
\varepsilon(G-B)\leqslant {{\upsilon(G-B)-\gamma(G-B)+1}\choose{2}}=
{{\upsilon(G)-\gamma(G)}\choose{2}}.
\end{array}
 $$
Thus $b(G)=|B|=\varepsilon(G)-\varepsilon(G-B)\geqslant \varepsilon-
{{\upsilon-\gamma(G)}\choose{2}}$.
\end{pf}

\vskip6pt

The lower bound in Theorem~\ref{thm4.6} is sharp for a class of
$\gamma$-critical graphs with domination number 2. By
Proposition~\ref{prop3.4.1}, $G$ is a complete graph $K_{2t}$
($t\geqslant2$) with a perfect matching removed. Then
$b(G)\geqslant\delta(G)+1=\Delta(G)+1$ by Theorem~\ref{thm4.6} and
$b(G)\leqslant\Delta(G)+1$ by Theorem~\ref{thm3.4.2}, and so
$b(G)=\Delta(G)+1=2t-1$.

\vskip6pt

So far as we know, there are no more lower bounds. In view of
applications of the bondage number, a network is vulnerable if its
bondage number is small while it is stable if its bondage number is
large. Therefore better lower bounds let us learn better stability
of networks from this point of view. In our opinion, it is of great
significance to seek more lower bounds for various classes of
graphs.

\section{Results on Graphs-operations}

Generally speaking, it is quite difficult to determine the exact
value of the bondage number for a given graph since it strongly
depends on the dominating number of the graph. Thus, determining
bondage numbers for some special graphs is interesting if the
dominating numbers of those graphs are known or can be easily
determined. In this section, we will introduce results on bondage
numbers for some special classes of graphs.

\vskip20pt

\subsection{Cartesian Product Graphs}

Let $G_1=(V_1,E_1)$ and $G_2=(V_2,E_2)$ be two graphs. The {\it
cartesian product} of $G_1$ and $G_2$ is an undirected graph,
denoted by $G_1\times  G_2$, where $V(G_1\times G_2)=V_1\times V_2$,
two distinct vertices $x_1x_2$ and $y_1y_2$, where $x_1,y_1\in
V(G_1)$ and $x_2, y_2\in V(G_2)$, are linked by an edge in
$G_1\times G_2$ if and only if either $x_1=y_1$ and $x_2y_2\in
E(G_2)$, or $x_2=y_2$ and $x_1y_1\in E(G_1)$. The cartesian product
is a very effective method for constructing a larger graph from
several specified small graphs.

\begin{thm} \label{thm5.7}\textnormal{(Dunbar {\it et al.}~\cite{dhtv98}, 1998)}
Let $G=C_n\times K_2$ with $n\geqslant3$. Then
$$
b(G)=\left\{\begin{array}{cl}%
2&\ {\rm if}\ n\equiv0\ {\rm or}\ 1\,({\rm mod}\ 4),\\
3&\ {\rm if}\ n\equiv3\,({\rm mod}\ 4),\\
4&\ {\rm if} \ n\equiv2\,({\rm mod}\ 4).
\end{array}\right.
$$
\end{thm}                                               % Proposition 4.2

For the Cartesian product $C_n\times C_m$ of two cycles $C_n$ and
$C_m$, where $n\geqslant 4$ and $m\geqslant 3$, Klavzar and
Seifter~\cite{ks95} determined $\gamma(C_n\times C_m)$ for some $n$'
and $m$'s. For example, $\gamma(C_n\times C_4)=n$ for $n\geqslant 3$
and $\gamma(C_n\times C_3)=n-\lfloor\frac n4\rfloor$ for $n\geqslant
4$. Kim~\cite{k03} determined $b(C_3\times C_3)=6$ and $b(C_4\times
C_4)=4$. For a general $n\geqslant 4$, the exact values of the
bondage numbers of $C_n\times C_3$ and $C_n\times C_4$ were
determined as follows.

\begin{thm}\label{thm5.8a} \textnormal{(Sohn, Yuan and Jeong~\cite{syj07}, 2007)}
\ For $n\geqslant 4$,
 $$
 b(C_n\times C_3)=\left\{
 \begin{array}{ll}
 2\ & {\rm if}\ n\equiv 0\,({\rm mod}\ 4),\\
 4\ & {\rm if}\ n\equiv 1\ {\rm or}\ 2\,({\rm mod}\ 4),\\
 5\ & {\rm if}\ n\equiv 3 \,({\rm mod}\ 4).
 \end{array}\right.
 $$
\end{thm}

\begin{thm}\label{thm5.8} \textnormal{(Kang, Sohn and Kim~\cite{ksk05}, 2005)}
\ $b(C_n\times C_4)=4$ for $n\geqslant 4$.
\end{thm}

For larger $m$ and $n$, Huang and Xu~\cite{hx08} obtained the
following result, see Theorem~\ref{thm10.8} for more details.

\begin{thm}\label{thm5.9} \textnormal{(Huang and
Xu~\cite{hx08}, 2008)}\ $b(C_{5i}\times C_{5j})=3$ for any positive
integers $i$ and $j$.
\end{thm}

Cao, Yuan and Moo~\cite{cym10} determined that for $n\geqslant 5$,
 $$
 b(C_n\times C_5)\ \left\{
 \begin{array}{ll}
 =3\ & {\rm if}\ n\equiv 0\ {\rm or}\ 1\,({\rm mod}\ 5),\\
 =4\ & {\rm if}\ n\equiv 2\ {\rm or}\ 4\,({\rm mod}\ 5),\\
 \leqslant 7\ & {\rm if}\  n\equiv 3\,({\rm mod}\ 5).
 \end{array}\right.
 $$

For the Cartesian product $P_n\times P_m$ of two paths $P_n$ and
$P_m$,
%\begin{prop}\label{prop5.5} {\rm
Jacobson and Kinch~\cite{jk83} determined $\gamma(P_n \times
P_2)=\lceil \frac{n+1}{2}\rceil$, $\gamma(P_n \times
P_3)=n-\lfloor\frac{n-1}{4}\rfloor$ and $\gamma(P_n \times P_4)=
n+1$ if $n=1,2,3,5,6, 9$, and $n$ otherwise. The bondage number of
$P_n\times P_m$ for $n\geqslant 4$ and $2\leqslant m\leqslant 4$ was
determined as follows (Li~\cite{l94} also determined $b(P_n \times
P_2)$).

\begin{thm}\label{thm5.10} {\rm (Hu, Cao and Xu~\cite{hcx09}, 2009)} For
$n\geq 4$,
 $$
 \begin{array}{rl}
 & b(P_n \times P_2)=\left\{ \begin{array}{l}
 1,\ {\rm if}\ n\equiv 1\,({\rm mod}\,2)  \\
 2, \ {\rm if}\ n\equiv 0\,({\rm mod}\,2);
 \end{array}
 \right.\\
 &b(P_n\times P_3)=\left\{ \begin{array}{l}
 1,\ {\rm if}\ n\equiv 1\ {\rm or}\ 2 ({\rm mod}\,4)  \\
2, \ {\rm if}\ n\equiv 0\ {\rm or}\ 3({\rm mod}\,4);\ {\rm and}
\end{array}
 \right.\\
 & b(P_n\times P_4)=1\ {\rm for}\ n\geq 14.
 \end{array}
$$
\end{thm}

From the proof of Theorem~\ref{thm5.10}, we find that if
$P_n=\{0,1,\ldots,n-1\}$ and $P_m=\{0,1,\ldots,m-1\}$, then the
removal of the vertex $(0,0)$ in $P_n\times P_m$ does not change the
domination number. If $m$ increase, the effect of $(0,0)$ for the
domination number will be smaller and smaller from the probability.
Therefore we expect it is possible that $\gamma (P_n\times
P_m-(0,0))=\gamma(P_n\times P_m)$ for $m\geq 5$ and give the
following conjecture.

\begin{con}\label{con5.11}\
$b(P_n\times P_m)\leqslant 2$ for $m\geq 5$.
\end{con}

%\vskip6pt
\subsection{Block Graphs and Cactus Graphs}

In this subsection, we introduce some results for corona graphs,
block graphs and cactus graphs.

The {\it corona} $G_1\circ G_2$, proposed by Frucht and Harary [R.
Frucht, F. Harary, On the corona of two graphs, Aequationes Math. 4
(1970), 322-324], is the graph formed from a copy of $G_1$ and
$\upsilon(G_1)$ copies of $G_2$ by joining the $i$-th vertex of
$G_1$ to the $i$-th copy of $G_2$. In particular, we are concerned
with the corona $H\circ K_1$, the graph formed by adding a new
vertex $v_i$ and a new edge $u_iv_i$ for every vertex $u_i$ in $H$.
Carlson and Develin~\cite{cd06} determined the bondage number of
$H\circ K_1$.

\begin{thm}\label{thm5.1} \textnormal{(Carlson and Develin~\cite{cd06}, 2006)}
\ $b(H\circ K_1)=\delta(H)+1$.
\end{thm}

\begin{pf}\
Let $G=H\circ K_1$, $\{u_1,u_2,\ldots, u_n\}$ be the vertex-set of
$H$ and $\{v_1,v_2,\ldots, v_n\}$ be the corresponding vertices
added in the construction of the corona. That is, for each
$i=1,2,\ldots,n$, the vertex $u_i$ is adjacent to $v_i$ via the edge
$e_i$. Then it is clear that $\gamma(G)=n$. In particular, any
minimal dominating set of $G$ contains exactly one of the vertices
$u_i$ and $v_i$ for each $i=1,2,\ldots,n$.

To show that $b(G)\geqslant\delta(H)+1$, we must show that if we
remove any $\delta(H)$ edges, the domination number of $G$ is
unchanged. Suppose, without loss of generality, that we remove $k$
edges $e_1, \ldots, e_k$ and $\delta(H)-k$ edges from $H$. Consider
the set $S=\{v_1, \ldots, v_k, u_{k+1}, \ldots, u_n\}$. Then
$|S|=n=\gamma(G)$. We claim $S$ is a dominating set of $G-S$. It is
clear that all $v_i$'s are in $N_G[S]$. Similarly, $u_j\in S$ for
$j>k$, so the only thing we need to check is that $u_j$ is adjacent
to an element of $S$ for $j\leqslant k$. Since
$d_H(u_j)\geqslant\delta(H)\geqslant k$, $u_j$ is adjacent to at
least one element of $\{u_{k+1}, \ldots, u_n\}$, completing the
proof of this direction.

To show $b(G)\leqslant\delta(H)+1$, take any vertex of minimum
degree in $H$, and delete its pendant edge and the pendant edges
incident on its $\delta(H)$ neighbors. It is easy to see that
deletion of these $\delta(H)+1$ pendant edges of $G$ increases the
domination number of the graph.
\end{pf}

\vskip6pt

A {\it block graph} is a graph whose blocks are complete graphs.
Each block in a {\it cactus graph} is either a cycle or a $K_2$. If
each block of a graph is either a complete graph or a cycle, then we
call this graph a {\it block-cactus graph}. Teschner \cite{t97b}
first studied the bondage numbers for these graphs.

\begin{thm}\label{thm5.2}  \textnormal{(Teschner~\cite{t97b}, 1997)}
$b(G)\leq\Delta(G)$ for any block graph $G$
\end{thm}                                                   % Theorem 4.70

Teschner~\cite{t97b} characterized all block graphs with
$b(G)=\Delta(G)$. At the same paper, Teschner found that
$\gamma$-critical graphs were instrumental
%(instrumental-ÓаïÖúµÄ)
in determining bounds for the bondage number of cactus and block
graphs and obtained the following result.

\begin{thm}\label{thm5.3}  \textnormal{(Teschner~\cite{t97b}, 1997)}
$b(G)\leqslant3$ for any nontrivial cactus graph $G$.
\end{thm}                                                   % Theorem 4.71

This bound can be achieved by $H\circ K_1$ where $H$ is a nontrivial
cactus graph without vertices of degree one by Theorem~\ref{thm5.1}.
In 1998, Dunbar {\it et al.}~\cite{dhtv98} proposed the following
problem.

\begin{prob}\label{prob5.4}
Characterize all cactus graphs with bondage number 3.
\end{prob}                                                   % Problem 4.6

Some upper bounds for block-cactus graphs were also obtained.

\begin{thm}\label{thm5.5}  \textnormal{(Teschner~\cite{t97b}, 1997)}
%; Teschner and Volkmann~\cite{tsvk})}
Let $G$ be a connected block-cactus graph with at least two blocks.
Then $b(G)\leqslant\Delta(G)$.
\end{thm}                                                   % Theorem 4.7

\begin{pf} Let $B$ be an end-block of $G$ with a cut vertex $x$ and
$y\in V(B)-\{x\}$. If $B$ is complete then $N_G[y]\subseteq N_G[x]$
and so $x$ is not critical. By Theorem~\ref{thm3.1.2},
$b(G)\leqslant\Delta(G)$. If $B$ is a cycle then there exist two
vertices $x$ and $y$ of degree $2$ with $d_G(x,y)\leqslant 2$, and
so $b(G)\leqslant3\leqslant\Delta(G)$ by Theorem~\ref{thm3.1.5}.
\end{pf}

\begin{thm}\label{thm5.6} \textnormal{(Dunbar {\it et al.}~\cite{dhtv98}, 1998)}
Let $G$ be a connected block-cactus graph which is neither a cactus
graph nor a block graph. Then $b(G)\leqslant\Delta(G)-1$.
\end{thm}                                                   % Theorem 4.8

%\subsection{For Regular and Vertex-transitive Graphs}

\section{Results on Planar Graphs}

From Section 2, we have seen that the bondage number for a tree has
been completely solved. Moreover, a linear time algorithm to compute
the bondage number of a tree was designed by Hartnell et
at.~\cite{hjvw98}. It is quite to consider the bondage number for a
planar graph. In this section, we will state some results and
problems on the bondage number for a planar graph.

Recall some results on planar graphs used in this section. For any
planar graph $G$,
 $$
\delta(G)\leqslant 5\ \ {\rm and}\ \ \varepsilon(G)\leqslant
3\upsilon(G)-6,
 $$
where $\upsilon(G)$ is the number of vertices and $\varepsilon(G)$
is the number of edges in $G$.

The well-known Euler's formula is stated as follows. For a connected
planar graph $G$,
 $$
\upsilon(G)-\varepsilon(G)+\phi(G)=2,
 $$
where $\phi(G)$ is the number of faces in any embedding of $G$ in
the plane or the sphere.

\subsection{Conjecture on the Bondage Number}

As mentioned in Section 3, the bondage number can be much larger
than the maximum degree. But for a planar graph $G$, the bondage
number $b(G)$ can not exceed $\Delta(G)$ too much. It is clear that
$b(G)\leqslant\Delta(G)+4$ by Corollary~\ref{con3.1.4} since
$\delta(G)\leqslant5$ for any planar graph $G$. In 1998, Dunbar {\it
et al.}~\cite{dhtv98} posed the following conjecture.

\begin{con}\label{con6.1}
%\textnormal{(Haynes, Teschner and Volkmann~\cite{dhtv98}, 1998)}
If $G$ is a planar graph, then $b(G)\leqslant\Delta(G)+1$.
\end{con}                                                % Conjecture 5.1

Because of its attraction, it immediately became the focus of
attention as soon as this conjecture is proposed. In fact, the main
aim concerning the research of the bondage number for a planar graph
is focused on this conjecture.

It has been mentioned in Theorem~\ref{thm2.1.1} and
Theorem~\ref{thm5.7} that $b(C_{3k+1})=3=\Delta+1$, and
$b(C_{4k+2}\times K_2)=4=\Delta+1$. It is easy to see that
$b(K_6-M)=5=\Delta+1$ where $M$ is a perfect matching of the
complete graph $K_6$. These examples show that if
Conjecture~\ref{con6.1} is true then the upper bound is best
possible for $2\leqslant\Delta\leqslant 4$.

Here, we note that it is sufficient to prove this conjecture for
connected planar graphs, since the bondage number of a disconnected
graph is simply the minimum of the bondage numbers of its
components.

\vskip6pt

\subsection{Bounds Implied by Maximum Degree}

The first paper attacking this conjecture is due to Kang and Yuan
\cite{ky00}, which confirmed the conjecture for every connected
planar graph $G$ with $\Delta\geqslant 7$. The proofs mainly base on
Theorem~\ref{thm3.1.5}, Theorem~\ref{thm3.1.7} and the following
lemma, which is a simple observation.

\begin{lem}\label{lem6.2} \textnormal{(Kang and Yuan
\cite{ky00}, 2000)} Let $G$ be a planar graph and $x\in V(G)$ with
$d_G(x)\geqslant 2$. Then there is a subset $F\subseteq \{yz:\ y,z
\in N_G(x)$ and $yz \notin E(G)\}$ such that $H=G+F$ is still a
planar graph and $H[N_G(x)]$ is connected when $d_G(x)=2$ and
$2$-connected when $d_G(x)\geqslant 3$.
\end{lem}                                                     % Lemma 5.2

Given a planar graph $H_0$ and its independent set
$X=\{x_1,\ldots,x_k\}$, by Lemma~\ref{lem6.2}, for each
$i=1,\ldots,k$, there exists $F_i\subseteq \{yz|\ y,z\in
N_{H_0}(x_i),y\neq z,yz\notin E(H_{i-1})\}$ such that
$H_i=H_{i-1}+F_i$ is planar and $H_i[N_{H_0}(x_i)]$ is connected
when $d_{H_0}(x_i)=2$ and $2$-connected when
$d_{H_0}(x_i)\geqslant3$. From these constructions Kang and Yuan
proved that if $b(G)\geqslant9$ for a planar graph $G$, then a
contradiction can be led to the well-known inequality
$\varepsilon\leqslant 3\upsilon-6$ for a planar graph. Based on
Theorem~\ref{thm3.1.5}, Theorem~\ref{thm3.1.7} and Euler's formula
on a connected planar graph, Kang and Yuan proved that
$b(G)\leqslant\Delta(G)+2$.

\begin{thm}\label{thm6.3} \textnormal{(Kang and Yuan
\cite{ky00}, 2000)}\ If $G$ is a connected planar graph then
$b(G)\leqslant\min\{8, \Delta(G)+2\}$.
\end{thm}                                                     % Theorem 5.3

%In 2006, Carlson and Develin~\cite{cd06} also gave a proof of
%$b(G)\leqslant\Delta(G)+2$, it is more intuitive than Kang and
%Yuan¡¯s proof, which simply proceeds by exhaustion.

Obviously, in view of Theorem~\ref{thm6.3}, Conjecture~\ref{con6.1}
is true for any connected planar graph with $\Delta\geqslant 7$.

%ÒÔÏÂÐðÊö²»¶Ô£¡
%{Thus, to Conjecture~\ref{con6.1}, we only need to
%consider the remaining cases, $5\leqslant \Delta\leqslant 6$, that
%is, Conjecture~\ref{con6.1} can be stated as follows.
%
%\begin{con}\label{con6.1a}
%If $G$ is a planar graph with $\Delta=5$ or $6$, then
%$b(G)\leqslant\Delta(G)+1$.
%\end{con}
%}

%Our proof is very topological in nature, and should extend to graphs
%embedded in other manifolds; it is more intuitive than Kang and
%Yuan¡¯s proof, which simply proceeds by exhaustion.

%\subsection{Bounds Implied by Forbidden Graphs}

%A graph $H$ is called a minor of a graph $G$, if $H$ can be obtained
%from $G$ by deleting a vertex, deleting an edge, or contracting an edge.

%\begin{lem}\label{lem6.3c} \textnormal{(Dirac~\cite{d60}, 1960)}\
%If the minimum degree of a graph $G$ is at least 3; then $G$ has a
%minor insomorphic to $K_4$.\end{lem}

%\bibitem{d60}
%G. A. Dirac, In abstraken Graphen vorhandene vollstandige 4-Graphen
%und ihre Unterteilungen, Math. Nath. 22 (1960), 61-85.

%\begin{thm}\label{thm6.3b} \textnormal{(Kang and Yuan
%\cite{ky00}, 2000)}\ If a connected graph $G$ has no minor
%isomorphic to $K_4$, then $b(G)\leqslant3$.
%\end{thm}                                                     % Theorem 5.3

%The fact that $b(C_4)=3$ shows that the result of this theorem is
%best possible.

\subsection{Bounds Implied by Degree-conditions}

As we have seen from Theorem~\ref{thm6.3}, to attack
Conjecture~\ref{con6.1}, we only need to consider connected planar
graphs with maximum $\Delta\leqslant 6$. Thus, studying the bondage
number of planar graphs by degree-conditions is of significance. The
first result on bounds implied by degree-conditions was obtained by
Kang and Yuan~\cite{ky00}.

\begin{thm} \label{thm6.4}\textnormal{(Kang and Yuan
\cite{ky00}, 2000)}\ If $G$ is a connected planar graph without
vertices of degree $5$, then $b(G)\leqslant7$.
%
%1)\ $b(G)\leqslant7$ if $G$ has no vertices of degree five;
%
%2)\ $b(G)\leqslant3$ if $G$ does not contain $K_4$.
\end{thm}                                                     % Theorem 5.4

%\begin{thm} \textnormal{\cite{ky00}}
%If $G$ is a connected planar graph then
%$b(G)\leqslant\Delta(G)+2$.
%\end{thm}                                                     % Theorem 5.5

%We say a graph H is a minor of a graph G, if H can be obtained from
%G by the following three operations: delete a vertex, delete an
%edge, or contract an edge.

%\begin{thm} \textnormal{\cite{ky00}}
%If $G$ is a connected planar graph with no minor isomorphic to
%$K_4$, then $b(G)\leqslant3$.
%\end{thm}                                                     % Theorem 5.6

%Since an outer-planar graph has no minor isomorphic to $K_4$, then
%\begin{cor} \textnormal{\cite{ky00}}
%If $G$ is an outer-planar graph then $b(G)\leqslant3$.
%\end{cor}                                                     % Corollary 5.7

%The second bound of Theorem~\ref{thm6.4} is best possible in view of $b(C_4)=3$.
%Whether the bound of Theorem 5.5 is sharp remains open.

As further applications of Lemma~\ref{lem6.2}, Fischermann,
Rautenbach and Volkmann~\cite{frv03} generalized
Theorem~\ref{thm6.4} as follows.

\begin{thm}\label{thm6.5} \textnormal{(Fischermann et al~\cite{frv03},
2003)}\ Let $G$ is a connected planar graph and $X$ be the set of
vertices of degree $5$ which have distance at least $3$ to vertices
of degree $1,2$ and $3$. If all vertices in $X$ not adjacent with
vertices of degree $4$ are independent and not adjacent to vertices
of degree $6$, then $b(G)\leqslant7$.
\end{thm}                                                       % Theorem 5.11

Clearly, if $G$ has no vertices of degree $5$, then $X=\emptyset$.
Thus, Theorem~\ref{thm6.5} yields Theorem~\ref{thm6.4}.

\vskip6pt

Use $\upsilon_i=\upsilon_i(G)$ to denote the number of vertices of
degree $i$ in $G$ for each $i=1,2,$ $\ldots,\Delta(G)$. Using
Theorem~\ref{thm3.1.7}, Fischermann {\it et al.} obtained the
following two theorems.

\begin{thm}\label{thm6.6}
\textnormal{(Fischermann {\it et al.}~\cite{frv03}, 2003)}\ Let $G$
be a connected planar graph. If
%there are two vertices $x$ and $y$ such that $d_G(x,y)\leqslant2$ and
%$d_G(x)+d_G(y)\leqslant8$, or if
$\upsilon_5<2\upsilon_2+3\upsilon_3+2\upsilon_4+12$, then
$b(G)\leqslant7$, and if $G$ contains no vertices of degree $4$ and
$5$, then $b(G)\leqslant6$.
\end{thm}                                                       % Theorem 5.12

\begin{thm}\label{thm6.7}
\textnormal{(Fischermann {\it et al.}~\cite{frv03}, 2003)}\ Let $G$
be a connected planar graph. Then $b(G)\leqslant\Delta(G)+1$ if

\noindent1) $\Delta(G)=6$ and every edge $e=xy$ with $d_G(x)=5$ and
$d_G(y)=6$ is contained in at most one triangle; or

\noindent2) $\Delta(G)=5$ and no triangle contains an edge $e=xy$
with $d_G(x)=5$ and $4\leqslant d_G(y)\leqslant5$.
\end{thm}                                                   % Theorem 5.14

\subsection{Bounds Implied by Girth-conditions}

The girth $g(G)$ of a graph $G$ is the length of the shortest cycle
in $G$. If $G$ has no cycles we define $g(G)=\infty$.

Combining Theorem~\ref{thm6.3} with Theorem~\ref{thm6.7}, we find
that if a planar graph contains no triangles and has maximum degree
$\Delta\geqslant 5$ then Conjecture~\ref{con6.1} holds. This fact
motivated Fischermann {\it et al.}~\cite{frv03} attempting to attack
Conjecture~\ref{con6.1} by grith restraints. They showed that the
conjecture is valid for all connected planar graphs of girth
$g(G)\geqslant 4$ and maximum degree $\Delta(G)\geqslant 5$ as well
as for all not $3$-regular graphs of girth $g(G)\geqslant 5$. In
this subsection, we will introduce their results and some opened
problems.

First, they improved the inequality $\varepsilon\leqslant
3\upsilon-6$ for a planar graph as follows.

\begin{lem}\label{lem6.8}
If $G$ is a planar graph with finite girth $g(G)\geqslant 3$, then
 $$
  \begin{array}{rl}
 \varepsilon(G)\leqslant\frac{g(G)(\upsilon(G)-2)-c(G)}{g(G)-2},
 \end{array}
 $$
where $c(G)$ is the number of cut-edges in $G$.
\end{lem}                                                         % Lemma 5.8

\begin{pf}\ Since every cut-edge is on the boundary of exactly
one face and every noncut-edge is on the boundary of two faces, we
reduce that $g(G)\phi (G)\leqslant 2\,\varepsilon(G)-c(G)$. Applying
Euler's Formula, the result follows.
\end{pf}

\vskip6pt

Then, by Lemma~\ref{lem6.8} and considering the relations among
$\upsilon_i$'s, the following result for planar graphs with girth
restraints is obtained.

\begin{thm}\label{thm6.9} \textnormal{(Fischermann {\it et al.}~\cite{frv03},
2003)}\ If $G$ is a connected planar graph then
 \begin{equation}\label{e6.1}
 b(G)\leqslant \left\{\begin{array}{ll}
 6, \ & {\rm if}\  g(G)\geqslant4;\\
 5, \ & {\rm if}\  g(G)\geqslant5;\\
 4, \ & {\rm if}\  g(G)\geqslant6;\\
 3, \ & {\rm if}\  g(G)\geqslant8.
 \end{array}\right.
 \end{equation}
\end{thm}                                                       % Theorem 5.9

\begin{pf}
Using Corollary~\ref{cor3.1.17} and Lemma~\ref{lem6.8}, we can give
a simple proof of the first three conclusions. Since the function
$f(x)=\frac{x}{x-2}$ is monotonically decreasing when $x\geqslant
3$, we have that
 \begin{equation}\label{e6.4.2}
 \begin{array}{rl}
\frac{g(G)}{g(G)-2}\leqslant \frac{k}{k-2}\ \ {\rm if}\
g(G)\geqslant k\geqslant 3.
 \end{array}
 \end{equation}
It follows from Corollary~\ref{cor3.1.17}, Lemma~\ref{lem6.8} and
(\ref{e6.4.2}) that
 $$%\begin{equation}\label{e6.4.3}
 \begin{array}{rl}
 \frac{\upsilon(G)}{4}(b(G)+1)\leqslant \varepsilon(G)\leqslant
 \frac{k}{k-2}\,(\upsilon(G)-2),
 \end{array}
 $$%\end{equation}
that is,
 \begin{equation}\label{e6.4.3}
 \begin{array}{rl}
 b(G)<\frac{3k+2}{k-2}.
 \end{array}
 \end{equation}
Substituting $k\geqslant 4,5,6$ into (\ref{e6.4.3}), respectively,
yields the desired conclusions.
\end{pf}

The first result in Theorem~\ref{thm6.9} shows that
Conjecture~\ref{con6.1} is valid for all connected planar graphs
with $g(G)\geqslant 4$ and $\Delta(G)\geqslant 5$. It is easy to
verify that the second result in Theorem~\ref{thm6.9} implies that
Conjecture~\ref{con6.1} is valid for all not $3$-regular graphs of
girth $g(G)\geqslant 5$, which is stated the following corollary.

\begin{cor}\label{con6.10}
%\textnormal{(Fischermann {\it et al.}~\cite{frv03}, 2003)}\
Let $G$ be a connected planar graph with $g(G)\geqslant 5$.
If $G$ is not $3$-regular, then $b(G)\leqslant \Delta(G)+1$.
\end{cor}

Since $b(C_{3k+1})=3$ for $k\geqslant 3$ (see
Theorem~\ref{thm2.1.1}), the last bound in Theorem~\ref{thm6.9} is
tight. Whether other bounds in Theorem~\ref{thm6.9} are tight
remains open. In 2003, Fischermann {\it et al.}~\cite{frv03} the
following conjecture.

\begin{con}\label{con6.11}
%\textnormal{(Fischermann {\it et al.}~\cite{frv03}, 2003)}\
If $G$ is a connected planar graph, then $b(G)\leqslant7$.
Furthermore,
$$
b(G)\leqslant \left\{\begin{array}{ll}
 5, \ & {\rm if}\  g(G)\geqslant4;\\
 4, \ & {\rm if}\  g(G)\geqslant5.
 \end{array}\right.
$$
\end{con}                                                    % Conjecture 6.3.4

We conclude this subsection with a question on bondage numbers of
planar graphs.

\begin{que}\label{que6.12}
\textnormal{(Fischermann {\it et al.}~\cite{frv03}, 2003)}\ Is there
a planar graph $G$ with $6\leqslant b(G)\leqslant8$?
\end{que}

In 2006, Carlson and Develin~\cite{cd06} showed that the corona
$G=H\circ K_1$ for a planar graph $H$ with $\delta(H)=5$ has the
bondage number $b(G)=\delta(H)+1=6$ (see Theorem~\ref{thm5.1}).
Since the minimum degree of planar graphs is at most $5$, then
$b(G)$ can attach $6$. If we take $H$ as the graph of the
icosahedron (see Figure~\ref{f7a} (a)), then $G=H\circ K_1$ is such
an example.

The question for the existence of planar graphs with bondage number
$7$ or $8$ remains open.

%We note that there exist planar graphs H with $\delta(H) = 5. Taking the
%corona G = H . K1 gives a planar graph with b(G) = 6.

\subsection{Comments on the Conjectures}

Conjecture~\ref{con6.1} is true for all connected planar graphs with
minimum degree $\delta\leqslant 2$ by Theorem~\ref{thm3.1.5}, or
maximum degree $\Delta\geqslant 7$ by Theorem~\ref{thm6.3}, or not
$\gamma$-critical planar graphs by Theorem~\ref{thm3.1.2}. Thus, to
attack Conjecture~\ref{con6.1}, we only need to consider connected
critical planar graphs with degree-restriction
$3\leqslant\delta\leqslant \Delta\leqslant 6$.

%holds with Theorem~\ref{thm6.7}, we find that if a planar graph
%contains no triangles and maximum degree $\Delta\geqslant 5$
% and has maximum degree $5$ or $6$ then Conjecture~\ref{con6.1} is holds.

Recalling and anatomizing %(×ÐϸÑо¿)
the proofs of all results mentioned in the preceding subsections on
the bondage number for connected planar graphs, we find that they
strongly depend upon Theorem~\ref{thm3.1.5} or
Theorem~\ref{thm3.1.7}. In other words, a basic way used in the
proofs is to find two vertices $x$ and $y$ with distance at most two
in a considered planar graph $G$ such that
 $$
 d_G(x)+d_G(y)\ {\rm or}\ d_G(x)+d_G(y)-|N_G(x)\cap N_G(y)|,
 $$
which bounds $b(G)$, is as small as possible. Let
 \begin{equation}\label{e6.2a}
 B(G)=\min\limits_{x,y\in V(G)}\left\{
 \begin{array}{l}
 \{d_G(x)+d_G(y)-1:1\leqslant d_G(x,y)\leqslant2\}\cup\\
 \{d_G(x)+d_G(y)-|N_G(x)\cap N_G(y)|-1:d(x,y)=1\}
 \end{array}\right\}.
 \end{equation}
Then, by Theorem~\ref{thm3.1.5} and Theorem~\ref{thm3.1.7}, we have
 \begin{equation}\label{e6.2}
 b(G)\leqslant B(G).
 \end{equation}
The proofs given in Theorem~\ref{thm6.3} and Theorem~\ref{thm6.9}
indeed imply the following stronger results.

\begin{thm}\label{thm6.5.1}
If $G$ is a connected planar graph then
$B(G)\leqslant\min\{8,\Delta(G)+2\}$.
\end{thm}                                                      % Theorem  6.3'

\begin{thm}\label{thm6.5.2}
If $G$ is a connected planar graph then
$$
B(G)\leqslant \left\{\begin{array}{ll}
 6, \ & {\rm if}\  g(G)\geqslant4;\\
 5, \ & {\rm if}\  g(G)\geqslant5;\\
 4, \ & {\rm if}\  g(G)\geqslant6;\\
 3, \ & {\rm if}\  g(G)\geqslant8.
 \end{array}\right.
$$
\end{thm}                                                       % Theorem 6.9'

Thus, using Theorem~\ref{thm3.1.5} or Theorem~\ref{thm3.1.7}, if we
can prove Conjecture~\ref{con6.1} and Conjecture~\ref{con6.11}, then
we can prove the following statement.

\begin{sta}\label{sta6.5.3}\ If $G$ is a connected planar
graph, then
$$
B(G)\leqslant \left\{\begin{array}{ll}
 \Delta(G)+1;\\
 7;\\
 5, \ & {\rm if}\  g(G)\geqslant4;\\
 4, \ & {\rm if}\  g(G)\geqslant5.
 \end{array}\right.
$$
\end{sta}

It follows from (\ref{e6.2}) that Statement~\ref{sta6.5.3} implies
Conjecture~\ref{con6.1} and Conjecture~\ref{con6.11}.
%This fact seems to means that the conjectures can be proved by
%Theorem~\ref{thm3.1.5} and Theorem~\ref{thm3.1.7}.
However, none of conclusions in Statement~\ref{sta6.5.3} is true by
the following examples.

To construct these counterexamples, we recall the operation of {\it
subdividing an edge}, i.e, replacing the edge $xy$ by a $2$-path
$xvy$ through a new vertex $v$, called the {\it subdividing vertex}
or {\it s-vertex\,} for short. We say {\it subdividing} the edge
$xy$ {\it twice} if $xy$ is replaced by a $3$-path $xv_1v_2y$.

% Figure 3.13
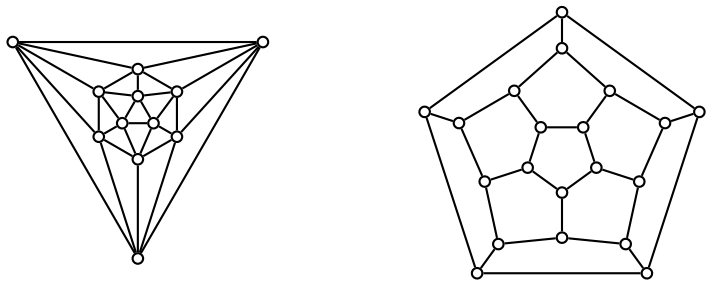
\begin{figure}[h]  % h: here; t: top; b: bottom
\psset{unit=1.2}%2,arrowsize=.13}     % scale

\hskip3cm

\begin{pspicture}(-8.5,-.6)(3.6,1.6)
\rput{54}{%
\SpecialCoor\degrees[5]
\multido{\i=0+1}{5}{\rput(.4;\i){\cnode{.07}{1\i}}}
\multido{\i=0+1}{5}{\rput(.9;\i){\cnode{.07}{2\i}}} }
\rput{18}{%
\SpecialCoor\degrees[5]
\multido{\i=0+1}{5}{\rput(1.2;\i){\cnode{.07}{3\i}}}
\multido{\i=0+1}{5}{\rput(1.6;\i){\cnode{.07}{4\i}}}
}%
\multido{\i=0+1}{5}{\ncline{1\i}{2\i}}
\multido{\i=0+1}{5}{\ncline{3\i}{4\i}}
\multido{\i=0+1}{5}{\ncline{2\i}{3\i}}
\ncline{20}{31}\ncline{21}{32}\ncline{22}{33}\ncline{23}{34}\ncline{24}{30}
\multido{\i=1+3}{2}{\ncline{\i0}{\i1}\ncline{\i1}{\i2}
\ncline{\i2}{\i3}\ncline{\i3}{\i4}\ncline{\i4}{\i0}}
\rput(0.,-1.8){\scriptsize (b)\ dodecahedron}
\end{pspicture}
%\newline
\begin{pspicture}(-3.8,-1.1)(-1.4,-2.)
\SpecialCoor%
\cnode(.2;-30){.07}{11}\cnode(.2;90){.07}{12}\cnode(.2;210){.07}{13}
\cnode(.5;30){.07}{21}\cnode(.5;90){.07}{22}\cnode(.5;150){.07}{23}
\cnode(.5;210){.07}{24}\cnode(.5;270){.07}{25}\cnode(.5;330){.07}{26}
\cnode(1.6;30){.07}{31}\cnode(1.6;150){.07}{32}\cnode(1.6;270){.07}{33}
\ncline{11}{12}\ncline{12}{13}\ncline{13}{11}
\ncline{21}{22}\ncline{22}{23}\ncline{23}{24}
\ncline{24}{25}\ncline{25}{26}\ncline{26}{21}
\ncline{31}{32}\ncline{32}{33}\ncline{33}{31}
\ncline{11}{21}\ncline{11}{25}\ncline{11}{26}
\ncline{12}{21}\ncline{12}{22}\ncline{12}{23}
\ncline{13}{23}\ncline{13}{24}\ncline{13}{25}
\ncline{31}{21}\ncline{31}{22}\ncline{31}{26}
\ncline{32}{22}\ncline{32}{23}\ncline{32}{24}
\ncline{33}{24}\ncline{33}{25}\ncline{33}{26}
\rput(0.,-2.3){\scriptsize (a)\ icosahedron}
\end{pspicture}

\caption{\label{f7a} \footnotesize  The two regular planar graphs}
\end{figure}

\begin{exa}\label{exa6.4.12}
For any $5$-regular planar graph $H$, subdivide each edge of $H$.
Then for each vertex $v\in V(H)$ there are five subdividing vertices
adjacent to $v$. Link these five vertices to form a cycle and keep
the planarity. The resulting graph $G$ has $\Delta(G)=6$ and
$B(G)=8$, which disproves the first two conclusions in
Statement~\ref{sta6.5.3}.
\end{exa}

%\begin{comment}
\begin{figure}[h]
\psset{unit=0.8}
\begin{pspicture}(-5,-3)(4,2)
\SpecialCoor \degrees[5]%
\dotnode(0;0){u}\rput(.3;.5){\scriptsize$u$}%
\dotnode(2;0){u1}\rput(2.1;-.1){\scriptsize$u_1$}%
\dotnode(2;1){u2}\rput(2.1;.9){\scriptsize$u_2$}%
\dotnode(2;2){u3}\rput(2.1;1.9.){\scriptsize$u_3$}%
\dotnode(2;3){u4}\rput(2.1;2.9){\scriptsize$u_4$}%
\dotnode(2;4){u5}\rput(2.1;3.9){\scriptsize$u_5$}%
\ncline{u}{u1}\ncline{u}{u2}\ncline{u}{u3}\ncline{u}{u4}\ncline{u}{u5}
\psline[doubleline=true,doublesep=.1,arrows=->](3;0)(4;0)
\rput(2.8;3.77){\scriptsize$u$ and $N_G(u)$ in $H$}%
\end{pspicture}
\begin{pspicture}(-3.5,-3)(3,2)
\SpecialCoor \degrees[5]%
\dotnode(0;0){u}\rput(.3;.5){\scriptsize$u$}%
\cnode(1;0){.08}{v1}\rput(1.2;-.18){\scriptsize$v_1$}%
\cnode(1;1){.08}{v2}\rput(1.2;.82){\scriptsize$v_2$}%
\cnode(1;2){.08}{v3}\rput(1.2;1.82){\scriptsize$v_3$}%
\cnode(1;3){.08}{v4}\rput(1.2;2.82){\scriptsize$v_4$}%
\cnode(1;4){.08}{v5}\rput(1.2;3.82){\scriptsize$v_5$}%
\dotnode(2;0){u1}\rput(2.1;-.1){\scriptsize$u_1$}%
\dotnode(2;1){u2}\rput(2.1;.9){\scriptsize$u_2$}%
\dotnode(2;2){u3}\rput(2.1;1.9.){\scriptsize$u_3$}%
\dotnode(2;3){u4}\rput(2.1;2.9){\scriptsize$u_4$}%
\dotnode(2;4){u5}\rput(2.1;3.9){\scriptsize$u_5$}%
\ncline{u}{v1}\ncline{u}{v2}\ncline{u}{v3}\ncline{u}{v4}\ncline{u}{v5}
\ncline{u1}{v1}\ncline{u2}{v2}\ncline{u3}{v3}\ncline{u4}{v4}\ncline{u5}{v5}
\ncline{v1}{v2}\ncline{v2}{v3}\ncline{v3}{v4}\ncline{v4}{v5}\ncline{v5}{v1}
\rput(2.8;3.77){\scriptsize$u,N_G(u)$ and $s$-vertices in $G$}%
\end{pspicture}
\caption{\label{f7}\footnotesize Example~\ref{exa6.4.12}}
\end{figure}
%\end{comment}

\begin{pf}
The existence of $H$ is guaranteed by the icosahedron (see
Figure~\ref{f7a} (a)). The construction is showed as
Figure~\ref{f7}. Assume $V(H)=\{u_1,\ldots,u_n\}$. Let
$S=\{v_1,\ldots,v_m\}$ be the set of subdividing vertices. Then
$V(G)=V(H)\cup S$. It is easy to observe that, $d_G(u_i)=5$ for
$i=1,\ldots,n$ for each $i=1,\ldots,n$, and $d_G(v_j)=6$ for each
$j=1,\ldots,m$. Furthermore, $d_G(u_i,u_j)=2$ ($i\ne j$) and
$|N_G(x)\cap N_G(y)|\leqslant2$ for every edge $xy$ of $G$. Thus
$B(G)=5+6-1-2=8=\Delta+2$.
\end{pf}

\begin{exa}\label{exa6.4.13}
For any $3$-regular planar graph $H$ with $g(H)\geqslant4$, we can
subdivide each edge twice and then link the subdividing vertices
properly such that the resulting graph $G$ is a planar graph with
$\Delta(G)=4$, $g(G)\geqslant4$ and $B(G)=6$, which disproves the
first and the third conclusions in Statement~\ref{sta6.5.3}.
\end{exa}

%\begin{comment}
\begin{figure}[h]
\begin{pspicture}(-2.2,-1)(3,4)
\dotnode(0,0){u1}%\rput(-.2,-.2){\scriptsize$u_1$}%
\dotnode(3,0){u2}%\rput(3.2,-.2){\scriptsize$u_2$}%
\dotnode(3,3){u3}%\rput(3.3,3){\scriptsize$u_3$}%
\dotnode(0,3){u4}%\rput(-.3,3){\scriptsize$u_4$}%
\dotnode(1,1){u5}%\rput(1.15,.8){\scriptsize$u_5$}%
\dotnode(2,1){u6}%\rput(1.9,.8){\scriptsize$u_6$}%
\dotnode(2,2){u7}%\rput(1.9,2.2){\scriptsize$u_7$}%
\dotnode(1,2){u8}%\rput(1.15,2.2){\scriptsize$u_8$}%
\ncline{u1}{u2}\ncline{u2}{u3}\ncline{u3}{u4}\ncline{u4}{u1}
\ncline{u5}{u6}\ncline{u6}{u7}\ncline{u7}{u8}\ncline{u8}{u5}
\ncline{u1}{u5}\ncline{u2}{u6}\ncline{u3}{u7}\ncline{u4}{u8}
\psline[doubleline=true,doublesep=.1,arrows=->](4,1.5)(5.5,1.5)
\end{pspicture}
\begin{pspicture}(-4,-1)(3,3)
\dotnode(0,0){u1}%\rput(-.2,-.2){\scriptsize$u_1$}%
\dotnode(3,0){u2}%\rput(3.2,-.2){\scriptsize$u_2$}%
\dotnode(3,3){u3}%\rput(3.3,3){\scriptsize$u_3$}%
\dotnode(0,3){u4}%\rput(-.3,3){\scriptsize$u_4$}%
\dotnode(1,1){u5}%\rput(1.15,.8){\scriptsize$u_5$}%
\dotnode(2,1){u6}%\rput(1.9,.8){\scriptsize$u_6$}%
\dotnode(2,2){u7}%\rput(1.9,2.2){\scriptsize$u_7$}%
\dotnode(1,2){u8}%\rput(1.15,2.2){\scriptsize$u_8$}%
\ncline{u1}{u2}\ncline{u2}{u3}\ncline{u3}{u4}\ncline{u4}{u1}
\ncline{u5}{u6}\ncline{u6}{u7}\ncline{u7}{u8}\ncline{u8}{u5}
\ncline{u1}{u5}\ncline{u2}{u6}\ncline{u3}{u7}\ncline{u4}{u8}
\pscurve(2,0)(-.5,-.5)(0,1)\pscurve(0,1)(-.5,3.5)(1,3)
\pscurve(1,3)(3.5,3.5)(3,2)\pscurve(3,2)(3.5,-.5)(2,0)
\cnode[fillstyle=solid](1,0){.07}{12}\cnode[fillstyle=solid](2,0){.07}{21}
\cnode[fillstyle=solid](3,1){.07}{23}\cnode[fillstyle=solid](3,2){.07}{32}
\cnode[fillstyle=solid](2,3){.07}{34}\cnode[fillstyle=solid](1,3){.07}{43}
\cnode[fillstyle=solid](0,2){.07}{41}\cnode[fillstyle=solid](0,1){.07}{14}
\cnode[fillstyle=solid](1.33,1){.07}{56}\cnode[fillstyle=solid](1.67,1){.07}{65}
\cnode[fillstyle=solid](2,1.33){.07}{67}\cnode[fillstyle=solid](2,1.67){.07}{76}
\cnode[fillstyle=solid](1.67,2){.07}{78}\cnode[fillstyle=solid](1.33,2){.07}{87}
\cnode[fillstyle=solid](1,1.67){.07}{85}\cnode[fillstyle=solid](1,1.33){.07}{58}
\cnode[fillstyle=solid](.33,.33){.07}{15}\cnode[fillstyle=solid](.67,.67){.07}{51}
\cnode[fillstyle=solid](2.67,.33){.07}{26}\cnode[fillstyle=solid](2.33,.67){.07}{62}
\cnode[fillstyle=solid](2.67,2.670){.07}{37}\cnode[fillstyle=solid](2.33,2.33){.07}{73}
\cnode[fillstyle=solid](.33,2.67){.07}{48}\cnode[fillstyle=solid](.67,2.33){.07}{84}
\ncline{12}{51}\ncline{51}{65}\ncline{65}{26}\ncline{26}{12}
\ncline{62}{76}\ncline{76}{37}\ncline{37}{23}\ncline{23}{62}
\ncline{73}{34}\ncline{34}{48}\ncline{48}{87}\ncline{87}{73}
\ncline{84}{58}\ncline{58}{15}\ncline{15}{41}\ncline{41}{84}
\ncline{56}{67}\ncline{67}{78}\ncline{78}{85}\ncline{85}{56}
\end{pspicture}
\caption{\label{f8}\footnotesize Example~\ref{exa6.4.13}}
\end{figure}
%\end{comment}

\begin{pf}
The construction is showed as Figure~\ref{f8} where $H$ is the cube.
The resulting graph $G$ is a planar graph with $g(G)\geqslant4$.
Note that all subdividing vertices have degree $4$ and
$d_G(u,v)\geqslant3$ for any two vertices $u,v\in V(H)$. Thus
$B(G)=4+3-1=6$. It is easy to verify that the result holds for any
$3$-regular planar graph.
\end{pf}

\begin{exa}\label{exa6.4.14}
The dodecahedron (see Figure~\ref{f7a} (b)) $G$ is a $3$-regular
planar graph with $g(G)=5$ and $B(G)=5$, which disproves the fist
and the late conclusions in Statement~\ref{sta6.5.3}.
\end{exa}

Examples \ref{exa6.4.12}$\,\sim\,$\ref{exa6.4.14} disprove
Statement~\ref{sta6.5.3}. As a result, we can state the following
conclusion.

\vskip6pt\begin{thm}\label{thm6.5.10} It is unable to prove
Conjecture \ref{con6.1}, and Conjecture \ref{con6.11}, if they are
right, using Theorem~\ref{thm3.1.5} and Theorem~\ref{thm3.1.7}.
\end{thm}

{\color{blue}\bf Therefore, a new method is need to prove these
conjectures, if they are right.}

\vskip6pt

%\noindent{\bf Remark} A careful argument for the dominations in
%Example \ref{exa6.4.12}$\,\sim\,$\ref{exa6.4.14} shows that their bondage
%numbers are not large enough to disprove Conjecture
%\ref{con1.1}$\,\sim\,$\ref{con1.4}. However, all the known upper
%bounds (cf. \cite{db98,wang96}) can not prove that these examples
%satisfy Conjecture \ref{con1.1}$\,\sim\,$\ref{con1.4}. Hence these
%conjectures can not be proved only by the known upper bounds.

\subsection{Minimum Counterexamples to the Conjectures}

As mentioned in the preceding subsection, now that we can not prove
these conjectures, then we may consider to disprove them. If one of
these conjectures is invalid, then there exists a {\it minimum
counterexample} $G$ with respect to $\upsilon(G)+\varepsilon(G)$.
Huang and Xu~\cite{hx12} investigated the property of the minimum
counterexample. Let $G_1$ be the possible existing minimum
counterexamples to Conjecture~\ref{con6.1}, $G_2,G_3$ and $G_4$ be
the possible existing minimum counterexamples to three conjectures
in Conjecture \ref{con6.11}, respectively. From the above
discussions, $G_i$ is a connected planar graph with $3\leqslant
\Delta\leqslant 6$ for each $i=1,2,3,4$. In particular,
 \begin{itemize}
 \item $G_1$ is $\gamma$-critical and $b(G_1)=\Delta(G_1)+2$ (by Theorem~\ref{thm3.1.2}).

 \item $b(G_2)=8$, $g(G_3)=3$ and $G_2$ contains vertices of degree $5$ (by Theorem~\ref{thm6.9} and Theorem \ref{thm6.4}).

 \item $b(G_3)=6$ and $g(G_3)=4$ (by Theorem \ref{thm6.9}).

 \item $b(G_4)=5$ and $g(G_4)=5$ (by Theorem \ref{thm6.9}).
\end{itemize}

In order to obtain further properties of these minimum
counterexamples, we consider how the bondage number changes under
some operation of a graph $G$ which decreases
$\upsilon(G)+\varepsilon(G)$ and preserves planarity. A simplest
operation satisfying this requirement is the edge deletion.

%\begin{lemma}\textnormal{\cite{ts97}}\label{lem3.1.1}
%Let $H$ be a spanning subgraph obtained by removing $k$ edges from
%$G$. Then $b(G)\leqslant b(H)+k$.
%\end{lemma}
%\begin{comment}
\begin{lem}\label{lem6.6.1}
Let $G$ be any graph and $e\in E(G)$. Then $b(G-e)\geqslant b(G)-1$.
Moreover, $b(G-e)\leqslant b(G)$ if $\gamma(G-e)=\gamma(G)$.
\end{lem}

\begin{pf}
Let $E'\subseteq E(G-e)$ with $|E'|=b(G-e)$. Then
$\gamma(G-e-E')>\gamma(G-e)$, and so $b(G)\leqslant|E'\cup
\{e\}|=b(G-e)+1$.

Assume $\gamma(G-e)=\gamma(G)$ and $E''\subseteq E(G)$ with
$|E''|=b(G)$. Then $\gamma(G-E'')>\gamma(G)$.

If $e\notin E''$ then
$\gamma(G-e-E'')\geqslant\gamma(G-E'')>\gamma(G)$.

If $e\in E''$ then
$\gamma(G-e-E'\setminus\{e\})=\gamma(G-E'')>\gamma(G)$.

Thus $b(G-e)\leqslant|E''|=b(G)$.
\end{pf}
%\end{comment}

\begin{lem}\label{lem6.6.2}
For any edge $e$ in $G_i$, $b(G_i-e)=b(G_i)-1$, $i=1,2,3,4$.
\end{lem}

\begin{pf}
By Lemma \ref{lem6.6.1}, $b(G_i-e)\geqslant b(G_i)-1$ for any edge
$e\in E(G_i)$. Note that $G_i-e$ is a planar graph with
$g(G_i-e)\geqslant g(G_i)$ and $\Delta(G_i-e)\leqslant\Delta(G_i)$.
Thus, if $b(G_i-e)\geqslant b(G_i)$, then $G_i-e$ is also a
counterexample, a contradiction to the minimum of $G_i$. Hence
$b(G_i-e)=b(G_i)-1$ for any edge $e\in E(G)$.
\end{pf}

\vskip6pt

From Lemma~\ref{lem6.6.2} we obtain the following conclusion
immediately.

{\color{blue}\begin{thm} It is unable to construct minimum
counterexamples to Conjecture~\ref{con6.1} and Conjecture
\ref{con6.11} by the operation of an edge deletion.
\end{thm}}

Next we consider the effect of the edge contraction on the bondage
number. Given a graph $G$, the {\it contraction} of $G$ by the edge
$e=xy$, denoted by $G/xy$, is the graph obtained from $G$ by
contracting two vertices $x$ and $y$ to a new vertex $n_{xy}$ and
then deleting all multi-edges. It is easy to observe that
$\upsilon(G/xy)+\varepsilon(G/xy)<\upsilon(G)+\varepsilon(G)$ and
$G/xy$ is also planar if $G$ is planar.

First, we investigate the influence of the edge contraction on the
domination and the bondage numbers for general graphs.

\begin{lem}\label{lem3.3} Let $G$ be any graph. Then
$\gamma(G)-1\leqslant\gamma(G/xy)\leqslant\gamma(G)$ for any edge
$xy$ of $G$.
\end{lem}

\begin{pf}
Let $S$ be a $\gamma$-set of $G$. If neither $x$ nor $y$ belongs to
$S$, then $S$ is a dominating set in $G/xy$. If
$S\cap\{x,y\}\ne\emptyset$, then $(S\setminus\{x,y\})\cup\{n_{xy}\}$
is a dominating set in $G/xy$, since $n_{xy}$ dominates all
neighbors of $x$ and $y$. It follows that
$\gamma(G/xy)\leqslant|S|=\gamma(G)$.

On the other hand, let $S'$ be a $\gamma$-set in $G/xy$. If
$n_{xy}\in S'$, then $S=S'\setminus\{n_{xy}\}\cup\{x,y\}$ is a
dominating set of $G$. If $n_{xy}\notin S'$, then $S'$ contains a
vertex $z$ such that $zn_{xy}\in E(G/xy)$. By the definition of edge
contraction, $zx\in E(G)$ or $yz\in E(G)$, which implies that
$S'\cup\{y\}$ or $S'\cup\{x\}$ is a dominating set of $G$. Thus
$\gamma(G)\leqslant|S'|+1=\gamma(G/xy)+1$.
\end{pf}

\begin{thm}\label{thm3.4}
Let $G$ be any graph and $xy$ be any edge in $G$. If $N_G(x)\cap
N_G(y)=\emptyset$ and $\gamma(G/xy)=\gamma(G)$, then
$b(G/xy)\geqslant b(G)$.
\end{thm}

\begin{pf}
Let $E'\subseteq E(G/xy)$ with $|E'|=b(G/xy)$ such that
 \begin{equation}\label{e6.6.4}
 \gamma(G/xy-E')>\gamma(G/xy).
 \end{equation}
Since $N_G(x)\cap N_G(y)=\emptyset$, the set of edges incident with
$x$ and $y$ except $xy$ in $G$ is the set of edges incident with
$n_{xy}$ in $G/xy$. It is easy to verify  that
 \begin{equation}\label{e6.6.5}
 (G-E')/xy=G/xy-E'.
 \end{equation}
Then by Lemma \ref{lem3.3}, (\ref{e6.6.5}) and (\ref{e6.6.4}), we
have
 $$
 \gamma(G-E')\geqslant\gamma((G-E')/xy)=\gamma(G/xy-E')>\gamma(G/xy)=\gamma(G),
 $$
and so $b(G)\leqslant |E'|=b(G/xy)$.
\end{pf}

\vskip6pt

Theorem \ref{thm3.4} is best possible. To illustrate this fact, we
recall two simple examples.

\begin{exa}\textnormal{(see Theorem~\ref{thm2.1.1})}\label{ex3.5}\
   $b(C_n)=\left\{\begin{array}{cl}%
   3&\ {\rm if}\ n\equiv 1(mod\ 3),\\
   2&\ {\rm otherwise},
   \end{array}\right.$ for $n\geqslant3$,
\end{exa}

\begin{exa}\textnormal{(see Theorem~\ref{thm2.1.1})}\label{ex3.6}\
$b(K_n)=\lceil \frac n2\rceil$ for $n\geqslant2$.
\end{exa}

\vskip6pt

The above examples show that the conditions of Theorem \ref{thm3.4}
are necessary. Clearly $\gamma(C_n)=\lceil n/3\rceil$ and
$\gamma(K_n)=1$; for any edge $xy$, $C_n/xy=C_{n-1}$ and
$K_n/xy=K_{n-1}$. By Example \ref{ex3.5}, if $n\equiv1$ (mod $3$),
then $\gamma(C_n/xy)<\gamma(G)$ and $b(C_n/xy)=2<3=b(C_n)$. Thus the
result in Theorem \ref{thm3.4} is generally invalid without the
hypothesis $\gamma(G/xy)=\gamma(G)$. Furthermore, the condition
$N_G(x)\cap N_G(y)=\emptyset$ can not be omitted even if
$\gamma(G/xy)=\gamma(G)$, since for odd $n$,
$b(K_n/xy)=\lceil\frac{n-1}2\rceil<\lceil\frac n2\rceil=b(K_n)$, by
Example \ref{ex3.6}.

On the other hand, the above examples also show that the equality in
$b(G/xy)\geqslant b(G)$ may hold ($b(C_n/xy)=b(C_n)=2$ if $n\equiv0$
(mod $3$), $b(K_n/xy)=b(K_n)$ if $n$ is even). Thus the bound in
Theorem \ref{thm3.4} is tight. However, provided all the conditions,
$b(G/xy)$ can be arbitrarily larger than $b(G)$. Given a graph $H$,
let $G$ be the graph formed from $H\circ K_1$ by adding a new vertex
$x$ and joining it to an vertex $y$ of degree one in $H\circ K_1$.
Then $G/xy=H\circ K_1$, $\gamma(G)=\gamma(G/xy)$ and $N_G(x)\cap
N_G(y)=\emptyset$. But $b(G)=1$ since $\gamma(G-xy)=\gamma(G/xy)+1$,
and $b(G/xy)=\delta(H)+1$ by Theorem~\ref{thm5.1}. The gap between
$b(G)$ and $b(G/xy)$ is $\delta(H)$.

\vskip6pt

Now we apply Theorem \ref{thm3.4} to $G_2$.

\begin{cor}\label{cor3.7}
If $\gamma(G_2/xy)=\gamma(G_2)$ for some edge $xy$, then
$N_{G_2}(x)\cap N_{G_2}(y)\ne\emptyset$.
\end{cor}

\begin{pf}
Note that $G_2/xy$ is a simple planar graph. If
$\gamma(G_2/xy)=\gamma(G_2)$ for some edge $xy$ and $N_{G_2}(x)\cap
N_{G_2}(y)=\emptyset$, then it follows from Theorem \ref{thm3.4}
that $b(G_2/xy)\geqslant b(G_2)$. Hence $G_2/xy$ is a
counterexample. But
$\varepsilon(G_2/xy)+\upsilon(G_2/xy)<\varepsilon(G_2)+\upsilon(G_2)$,
a contradiction.
\end{pf}

\vskip6pt

By Theorem~\ref{thm6.9}, $G_2$ must contain triangles, and so there
is an edge $xy$ with $N_{G_2}(x)\cap N_{G_2}(y)\ne\emptyset$. From
Theorem~\ref{thm3.4} and Corollary~\ref{cor3.7} we obtain the
following conclusion immediately.

{\color{blue}\begin{thm}\label{thm6.6.9} It is unable to construct a
minimum counterexample to the first conjecture in Conjecture
\ref{con6.11}, that is $G_2$, by the operation of an edge
contraction.
\end{thm}}

\vskip6pt

Finally we consider $G_1$.

\begin{lem}\label{lem3.8} Let $G$ be any graph and $xy$ be any edge in
$G$. Then $\gamma(G/xy)\leqslant\gamma(G-x)$.
\end{lem}

\begin{pf}
Let $S$ be a $\gamma$-set of $G-x$. If $y\notin S$, then there
exists a vertex $z\in S$ such that $yz\in E(G-x)$. Thus $zn_{xy}\in
E(G/xy)$, i.e., $n_{xy}$ is dominated by $z\in S$. Therefore $S$ is
also a dominating set of $G/xy$ and
$\gamma(G/xy)\leqslant|S|=\gamma(G)$.

Now assume $y\in S$ and let $S'=(S\setminus\{y\})\cup\{n_{xy}\}$. If
$yz\in E(G-x)$, then $zn_{xy}\in E(G/xy)$, which means that the
vertices dominated by $y$ in $G-x$ are all dominated by $n_{xy}$ in
$G/xy$. Thus $S'$ is a dominating set of $G/xy$ and
$\gamma(G/xy)\leqslant|S'|=\gamma(G)$.
\end{pf}

\begin{lem}\label{lem3.10}
$\gamma(G_1/xy)=\gamma(G_1)-1$ for every edge $xy$ in $G_1$.
\end{lem}

\begin{pf}
Let $xy$ be an edge of $G_1$. Since $b(G_1)>\Delta(G_1)+1$, then
$G_1$ is not $\gamma$-critical by Theorem~\ref{thm3.1.2}, that is,
$\gamma(G_1-u)<\gamma(G_1)$ for any vertex $u\in V(G_1)$. By Lemma
\ref{lem3.8}, we have
$\gamma(G_1/xy)\leqslant\gamma(G_1-x)<\gamma(G_1)$.
\end{pf}

\section{Results on Crossing Number Restraints}

It is quite natural to generalize the known results on the bondage
number for planar graphs to for more general graphs in terms of
graph-theoretical parameters. In this section, we consider graphs
with crossing number restraints.

The {\it crossing number} $cr(G)$ of $G$ is the smallest number of
pairwise intersections of its edges when $G$ is drawn in the plane.
If $cr(G)=0$, then $G$ is a planar graph. A spanning subgraph $H$ of
$G$ is called a {\it maximum planar subgraph} of $G$ if $H$ is
planar and contains as many edges of $G$ as possible. we can easily
observe the following property of maximum planar subgraphs.

\begin{lem} \label{lem7.1} Let $G$ be a graph with $cr(G)>0$ and $H$ a
maximum planar subgraph of $G$. Then

(1)\ $0<|E(G)|-|E(H)|\leqslant cr(G)$;

(2)\ $H$ contains a cycle;

(3)\ both $G$ and $H$ have the same number of components.
\end{lem}                                                     % Lemma 7.1

\subsection{General Methods}

Use $\upsilon_i(G)$ to denote the number of vertices of degree $i$
in $G$ for $i=1,2,\ldots,\Delta(G)$. By using Lemma~\ref{lem7.1}
with some more effort of computations, Huang and Xu obtained the
following results.

\begin{thm}\label{thm7.2} \textnormal{(Huang and
Xu~\cite{hx07}, 2007)}\ Let G be a connected graph, then
$$
 b(G)\leqslant \left\{\begin{array}{ll}
 6\ &{\rm if}\ g(G)\geqslant 4\ {\rm and}\
2cr(G)<\upsilon_1+2\upsilon_2+2\upsilon_3+\sum\limits_{i=8}^{\Delta}(i-7)\upsilon_i+8;\\
5\ &{\rm if}\ g(G)\geqslant 5\ {\rm and}\
6cr(G)<3\upsilon_1+6\upsilon_2+5\upsilon_3+ \sum\limits_{i=7}^{\Delta}(3i-18)\upsilon_i+20;\\
4\ &{\rm if}\ g(G)\geqslant 6\ {\rm and}\
4cr(G)<\upsilon_1+2\upsilon_2+\sum\limits_{i=6}^{\Delta}(2i-10)\upsilon_i+12;\\
3\ &{\rm if}\ g(G)\geqslant 8\ {\rm and}\
6cr(G)<\sum\limits_{i=5}^{\Delta}(3i-12)\upsilon_i+16.
\end{array}\right.
$$
\end{thm}                                                 % Theorem 6.2

\begin{cor}\label{cor7.3} \textnormal{(Huang and
Xu~\cite{hx07}, 2007)}\ For is a connected graph $G$,
$$
b(G)\leqslant \left\{\begin{array}{ll}
 6, \ & {\rm if}\  g(G)\geqslant4\ {\rm and}\ cr(G)\leqslant 3;\\
 5, \ & {\rm if}\  g(G)\geqslant5\ {\rm and}\ cr(G)\leqslant 4;\\
 4, \ & {\rm if}\  g(G)\geqslant6\ {\rm and}\ cr(G)\leqslant 2;\\
 3, \ & {\rm if}\  g(G)\geqslant8\ {\rm and}\ cr(G)\leqslant 2.
 \end{array}\right.
$$
\end{cor}                                       % Corollary 6.3

\begin{cor}\label{cor7.4} \textnormal{(Huang and
Xu~\cite{hx07}, 2007)}\ Let G be a connected graph. Then

{\rm (a)}\ $b(G)\leqslant 6$ if $G$ is not $4$-regular, $cr(G)=4$
and $g(G)\geqslant 4$;

{\rm (b)}\ $b(G)\leqslant\Delta(G)+1$ if $G$ is not $3$-regular,
$cr(G)\leqslant 4$ and $g(G)\geqslant 5$;

{\rm (c)}\ $b(G)\leqslant 4$ if $G$ is not $3$-regular, $cr(G)=3$
and $g(G)\geqslant 6$;

{\rm (d)}\ $b(G)\leqslant 3$ if $cr(G)=3$, $g(G)\geqslant 8$ and
$\Delta(G)\geqslant 5$.
\end{cor}                                 % Corollary 6.4

These corollaries generalize some known results for planar graphs.
For example, Corollary~\ref{cor7.3} contains Theorem~\ref{thm6.9};
the conclusion (b) in Corollary~\ref{cor7.4} contains
Corollary~\ref{con6.10}.

In order to generalized other known results for planar graphs, we
make the following construction. For a connected graph $G$, let
$G_0$ be a subgraph of $G$ without isolated vertices, $H_0$ be a
maximum planar subgraph of $G_0$ and let $E'=E(G_0)\setminus
E(H_0)$.
% and $X_i=\{x\in V(G)|\ d_G(x)=i\}$.

Suppose $X=\{x_1,x_2,\ldots,x_k\}$ is a given independent set of
$G_0$ with $d_{G_0}(x_i)\geqslant2$ for each $i=1,2,\ldots,k$. By
Lemma~\ref{lem6.2}, there exists $F_i\subseteq \{yz:\ y,z\in
N_{H_0}(x_i), y\neq z, yz\notin E(H_{i-1})\}$ such that
$H_i=H_{i-1}+F_i$ is planar and $H_i[N_{H_0}(x_i)]$ is connected
when $d_{H_0}(x_i)=2$ and $2$-connected when $d_{H_0}(x_i)\geqslant
3$ for $i=1,2,\ldots,k$. Let $G_k=H_k+E'\setminus
E(H_k)=H_k+E'\setminus \cup_{i=1}^k F_i$.

After these constructions we obtained a planar graph $H_k$. Then
$\varepsilon(H_k)\leqslant 3\upsilon(H_k)-6=3\upsilon(G_k)-6$.
Applying this inequality, Huang and Xu~\cite{hx07} obtained the
following results for graphs with crossing number restraints, which
generalize some known results for planar graphs.

\begin{thm} \textnormal{(Huang and
Xu~\cite{hx07}, 2007)}\ If G is a connected graph with
$cr(G)\leqslant\upsilon_3(G)+\upsilon_4(G)+3$, then $b(G)\leqslant
8$.
 \end{thm}                                                      %Theorem 6.5

\begin{cor}\label{cor7.6} \textnormal{(Huang and
Xu~\cite{hx07}, 2007)}\ If G is a connected graph with
$cr(G)\leqslant 3$, then $b(G)\leqslant 8$.
\end{cor}                                                     % Corollary 6.6

{\bf Remark:}\ Perhaps being unaware of this result, in 2010, Ma
{\it et al.}~\cite{mzw10} proved that $b(G)\leqslant 12$ for any
graph $G$ with $cr(G)=1$.

\begin{thm} \textnormal{(Huang and
Xu~\cite{hx07}, 2007)}\ Let G be a connected graph, and $I=\{x\in
V(G):\ d_G(x)=5,\ d_G(x,y)\geqslant 3$ if $d_G(y)\leqslant 3$, and
$d_G(y)\neq 4$ for every $y\in N_G(x)\}$. If $I$ is independent, has
no vertex adjacent to vertices of degree 6 and
 $$
  \begin{array}{rl}
cr(G)<\max\left\{\frac{5\upsilon_3(G)+|\,I|-2\upsilon_4(G)+28}{11},\frac{7\upsilon_3(G)+40}{16}\right\},
 \end{array}
 $$
then $b(G)\leqslant 7$.
\end{thm}                                                       % Theorem 6.7

\begin{cor} \textnormal{(Huang and
Xu~\cite{hx07}, 2007)}\ Let G be a connected graph with
$cr(G)\leqslant2$. If $I=\{x\in V(G):\ d_G(x)=5, d_G(y,x)\geqslant3$
if $d_G(y)\leqslant 3$ and $d_G(y)\neq 4$ for every $y\in N_G(x)\}$
is independent, and has no vertices adjacent to vertices of degree
6, then $b(G)\leqslant 7$.
\end{cor}                                                     % Corollary 6.8

\begin{thm}\label{lem7.9} \textnormal{(Huang and
Xu~\cite{hx07}, 2007)}\ Let $G$ be a connected graph. If $G$ satisfies\\
$1)\ 5cr(G)+\upsilon_5<2\upsilon_2+3\upsilon_3+2\upsilon_4+12$; or\\
$2)\ 7cr(G)+2\upsilon_5<3\upsilon_2+4\upsilon_4+24$,\\ then
$b(G)\leqslant7$.
\end{thm}                                                       % Theorem 6.9

\begin{prop} \textnormal{(Huang and
Xu~\cite{hx07}, 2007)}\ Let $G$ be a connected graph with no
vertices of degree four and five. If $cr(G)\leqslant 2$, then
$b(G)\leqslant 6$.
\end{prop}                                                   % Proposition 6.10

\begin{thm}\label{lem7.11} \textnormal{(Huang and
Xu~\cite{hx07}, 2007)}\ If $G$ is a connected graph with
$cr(G)\leqslant4$ and not $4$-regular when $cr(G)=4$, then
$b(G)\leqslant\Delta(G)+2$.
\end{thm}                                                       % Theorem 6.11

The above results generalize some known results for planar graphs.
For example, Corollary~\ref{cor7.6} and Theorem~\ref{lem7.11}
contain Theorem~\ref{thm6.3}; the first condition in
Theorem~\ref{lem7.9} contains the second condition in
Theorem~\ref{thm6.6}.

\vskip6pt From Corollary~\ref{cor7.4} and Theorem~\ref{lem7.11} we
suggest following questions.

\begin{que} Is there a

\noindent(a) $4$-regular graph $G$ with $cr(G)=4$ such that
$b(G)\geqslant7$?

\noindent(b) $3$-regular graph $G$ with $cr\leqslant4$ and
$g(G)\geqslant5$ such that $b(G)=5$?

\noindent(c) $3$-regular graph $G$ with $cr=3$ and $g(G)=6$ or $7$
such that $b(G)=5$?
\end{que}

\subsection{Carlson and Develin's Methods}

In this subsection, we introduce an elegant method presented by
Carlson and Develin~\cite{cd06} to obtain some upper bounds for the
bondage number of a graph.

Suppose that $G$ is a connected graph. We say that $G$ has {\it
genus} $\rho$ if $G$ can be embedded in a surface $S$ with $\rho$
handles such that edges are pairwise disjoint except possibly for
end-vertices. Let $D$ be an embedding of $G$ in surface $S$, and let
$\phi(G)$ denote the number of regions in $D$. The boundary of every
region contains at least three edges and every edge is on the
boundary of at most two regions (the two regions are identical when
$e$ is a cut-edge). For any edge $e$ of $G$, let $r_G^1(e)$ and
$r_G^2(e)$ be the numbers of edges comprising the regions in $D$
which the edge $e$ borders. It is clear that every $e=xy\in E(G)$,
  \begin{equation}\label{e7.1}
 \left\{\begin{array}{ll}
 r_G^2(e)\geqslant r_G^1(e)\geqslant 4\ &\ {\rm if}\ |N_G(x)\cap N_G(y)|=0,\\
 r_G^2(e) \geqslant 4,\ r_G^1(e)\geqslant 3\  &\ {\rm if}\ |N_G(x)\cap N_G(y)|=1,\\
 r_G^2(e)\geqslant r_G^1(e)\geqslant 3\ &\ {\rm if}\ |N_G(x)\cap N_G(y)|\geqslant2.
 \end{array}\right.
  \end{equation}
Following Carlson and Develin~\cite{cd06}, for any edge $e=xy$ of
$G$, we define
 \begin{equation}\label{e7.2}
 D_G(e)=\frac{1}{d_G(x)}+\frac{1}{d_G(y)}+\frac{1}{r_G^1(e)}+\frac{1}{r_G^2(e)}-1.
 \end{equation}
By the well-known Euler's Formula
 $$%\begin{equation}\label{e2.10}
 \upsilon(G)-\varepsilon(G)+\phi(G)=2-2\rho(G),
 $$%\end{equation}
it is easy to see that
 \begin{equation}\label{e7.3}
 \sum_{e\in E(G)}D_G(e)=\upsilon(G)-\varepsilon(G)+\phi(G)=2-2\rho(G).
 \end{equation}
If $G$ is a connected planar graph, that is, $\rho(G)=0$, then
 \begin{equation}\label{e7.4}
 \sum_{e\in E(G)}D_G(e)=\upsilon(G)-\varepsilon(G)+\phi(G)=2.
 \end{equation}

Combining these formulas with Theorem~\ref{thm3.1.7}, Carlson and
Develin \cite{cd06} gave a simple and intuitive proof of
Theorem~\ref{thm6.3} and obtained the following result.

\begin{thm}\label{thm7.13} \textnormal{(Carlson and Develin
\cite{cd06}, 2006)} Let $G$ be a connected graph which can be
embedded on a torus. Then $b(G)\leqslant\Delta(G)+3$.
\end{thm}                                                       % Theorem 5.13

\begin{pf}Suppose that $G$ is a graph which can be embedded on a
torus, for which $\rho(G)=2$. By Corollary~\ref{con3.1.4}, if
$\delta(G)\leqslant 4$ then theorem holds, so we can assume
$\delta(G)\geqslant 5$. For the sake of contradiction, assume
$b(G)\geqslant \Delta (G) + 4$. Let $e=xy$ be an arbitrary edge of
$G$ and let, without loss of generality, $d_G(x)\leqslant d_G(y)$.
By Theorem~\ref{thm3.1.7} we should have,
 \begin{equation}\label{e7.5}
 \Delta(G)+4\leqslant b(G)\leqslant
 d_G(x) + d_G(y)-1-|N_G(x)\cap N_G(y)|.
 \end{equation}
If $d_G(x)=5$, then (\ref{e7.5}) implies $d_G(y)=\Delta(G)$ and
$|N_G(x)\cap N_G(y)|=0$, $r^2_ e\geqslant r^1_e\geqslant 4$ in
(\ref{e7.1}), and so $D_e<0$ in (\ref{e7.2}).

If $d_G(x)=6$, then $|N_G(x)\cap N_G(y)|=1$ by (\ref{e7.5}),
$r^2_e\geqslant 4$ and $r^1_e\geqslant 3$ in (\ref{e7.1}), and so
hence $D_e< 0$ in (\ref{e7.2}).

If $d_G(x)\geqslant 7$, then clearly $D_e<0$.

Therefore $\sum\limits_{e\in E(G)}D_e<0$, which is a contradiction
to (\ref{e7.3}) when $\rho(G)=2$.
\end{pf}

\vskip6pt

By the way, Cao, Xu and Xu~\cite{cxx09} generalized the result in
(\ref{e6.1}) to a connected graph $G$ that can be embedded on a
torus, that is,
 $$
 b(G)\leq\left\{
\begin{array}{rl}
&6, \ \ \mbox{if}\ \   g(G)\geq 4\ {\rm and}\ G\ {\rm is\ not}\ 4{\rm-regular};\\
&5, \ \ \mbox{if}\ \   g(G)\geq 5;\\
&4, \ \ \mbox{if}\ \   g(G)\geq 6\ {\rm and}\ G\ {\rm is\ not}\ 3{\rm-regular};\\
&3, \ \ \mbox{if}\ \   g(G)\geq 8.
\end{array} \right.
 $$

Several authors used this method to obtain some results on the
bondage number. For example, Fischermann {\it et al.}~\cite{frv03}
used this method to prove the second conclusion in
Theorem~\ref{thm6.7}. Recently, Cao, Huang and Xu~\cite{chx08} have
used this method to deal with more general graphs with small
crossing numbers. First, they found the following property by
Lemma~\ref{lem7.1}.

\begin{lem}\label{lem7.14}
$\delta(G)\leqslant5$ for any graph $G$ with $cr(G)\leqslant5$.
\end{lem}

This result implies that $b(G)\leqslant\Delta(G)+4$ by
Corollary~\ref{con3.1.4} for any graph $G$ with $cr(G)\leqslant5$.
Cao et al established the following relation in terms of maximum
planar subgraphs.

\begin{lem}\label{lem7.15}
\textnormal{(Cao {\it et al.}~\cite{chx08}, 2008)} Let $G$ be a
connected graph with crossing number $cr(G)$ and let $H$ be a
maximum planar subgraph of $G$. Then
 $$
\sum\limits_{e\in E(H)}D_G(e)\geqslant2-\frac{2cr(G)}{\delta(G)}.
$$
\end{lem}

Using Carlson and Develin's method and combining Lemma~\ref{lem7.14}
with Lemma~\ref{lem7.15}, Cao {\it et al.} proved the following
results.

\begin{thm}\label{thm7.16} \textnormal{(Cao {\it et al.}~\cite{chx08}, 2008)}
Let $G$ be a connected graph. Then
$b(G)\leqslant\Delta(G)+2$ if $G$ satisfies one of the following
conditions:

{\rm (a)}\ $cr(G)\leqslant3$,

{\rm (b)}\ $cr(G)=4$ and $G$ is not $4$-regular,

{\rm (c)}\ $cr(G)=5$ and $G$ contains no vertices of degree $4$.
\end{thm}                                                       % Theorem 6.16

\begin{thm}\label{thm7.17}\textnormal{(Cao {\it et al.}~\cite{chx08}, 2008)}
Let $G$ be a connected graph with $\Delta(G)=5$ and
$cr(G)\leqslant4$. If no triangles contain two vertices of degree
$5$, then $b(G)\leqslant6=\Delta(G)+1$.
\end{thm}                                                       % Theorem 6.17

\begin{thm}\label{thm7.18} \textnormal{(Cao {\it et al.}~\cite{chx08}, 2008)}
Let $G$ be a connected graph with $\Delta(G)\geqslant6$ and
$cr(G)\leqslant3$. If $\Delta(G)\geqslant7$ or if $\Delta(G)=6$,
$\delta(G)\ne 3$ and every edge $e=xy$ with $d_G(x)=5$ and
$d_G(y)=6$ is contained in at most one triangle, then
$b(G)\leqslant\min\{8,\Delta(G)+1\}$.
\end{thm}                                                   % Theorem 6.18

Using Carlson and Develin's method, it can be proved that
$b(G)\leqslant\Delta(G)+2$ if $cr(G)\leqslant3$ (see the first
conclusion in Theorem~\ref{thm7.16}), and not yet proved that
$b(G)\leqslant 8$ if $cr(G)\leqslant3$, although it has been proved
by using other method (see Corollary~\ref{cor7.6}).
%Recall our
%proofs of Theorem~\ref{thm7.16} and Theorem~\ref{thm7.17}, we used
%the conclusion that $D_G(e)\leqslant0$ for any $e=xy\in E(G)$.
%However, there are many graphs

\vskip6pt

By using Carlson and Develin's method, Samodivkin~\cite{s09}
obtained some results on the bondage number for graphs with some
given properties.

Kim~\cite{k04} showed that $b(G)\leqslant\Delta(G)+2$ for a
connected graph $G$ with genus $1$, $\Delta(G)\leqslant 5$ and
having a torodial embedding of which at least one region is not
$4$-sided. Recently, Gagarin and Zverovich~\cite{gz10} further have
extended Carlson and Develin's ideas to establish a nice upper bound
for arbitrary graphs that cab be embedded on orientable or
nonorientable topological surface.

\begin{thm}\label{thm7.2.7} \textnormal{(Gagarin and Zverovich~\cite{gz10}, 2010)}
Let $G$ be a graph
embeddable on an orientable surface of genus $h$ and a
non-orientable surface of genus $k$. Then
$b(G)\leqslant\min\{\Delta(G)+h+2,\Delta(G)+k+1\}$.
\end{thm}

%{gz10} A. Gagarin and V. Zverovich, Upper bounds for the bondage
%number of graphs on topological surfaces, ArXiv: 1012.4117, 2010.

This result generalizes the corresponding upper bounds of
Theorems~\ref{thm6.3} and Theorems~\ref{thm7.13} for any orientable
or nonorientable topological surface.

By investigating the proof of Theorem~\ref{thm7.2.7},
Huang~\cite{h11} found that the issue of the orientability can be
avoided by using the Euler characteristic $\chi
(=\upsilon(G)-\varepsilon(G)+\phi(G))$ instead of the genera $h$ and
$k$, the relations are $\chi=2-2h$ and $\chi=2-k$. To have the best
result from Theorem~\ref{thm7.2.7}, one wants $h$ and $k$ as small
as possible, this is equivalent to having $\chi$ as large as
possible.

According to Theorem~\ref{thm7.2.7}, if $G$ is planar $(h=0,
\chi=2)$ or can be embedded on the real projective plane $(k=1,
\chi=1)$, then $b(G)\leqslant\Delta(G)+2$. In all other cases,
Huang~\cite{h11} had the following improvement for
Theorem~\ref{thm7.2.7}, the proof is based on the technique
developed by Carlson-Develin and Gagarin-Zverovich, and includes
some elementary calculus as a new ingredient, mainly the
intermediate value theorem and the mean value theorem.

\begin{thm}\label{thm7.2.8} \textnormal{(Huang~\cite{h11}, 2011)}
Let $G$ be a graph embeddable on a surface whose Euler
characteristic $\chi$ is as large as possible. If $\chi\leqslant 0$
then $b(G)\leqslant\Delta(G)+\lfloor r\rfloor$, where $r$ is the
largest real root of the following cubic equation in $z$:
 $$
 z^3+2z^2+(6\chi-7)z+18\chi-24=0.
 $$
In addition, if $\chi$ decreases then $r$ increases.
\end{thm}

The following result is asymptotically equivalent to
Theorem~\ref{thm7.2.8}.

\begin{thm}\label{thm7.2.9} \textnormal{(Huang~\cite{h11}, 2011)}
Let $G$ be a graph embeddable on a surface whose Euler
characteristic $\chi$ is as large as possible. If $\chi\leqslant 0$
then $b(G)<\Delta(G)+\sqrt{12-6\chi}+1/2$, or equivalently,
$b(G)\leqslant\Delta(G)+\lceil\sqrt{12-6\chi}-1/2\rceil$.
\end{thm}

\section{Other Types of Bondage Numbers}

Since the concept of the bondage is based upon the domination, all
sorts of dominations, which are generalizations of the normal
domination or adding restricted conditions to the normal dominating
set, can develop a new ``bondage number'' as long as a variation of
domination number is given. In this section, we will survey some
results on the bondage number under some restricted conditions.

%variant various variation

\subsection{Restrained Bondage Numbers}

A dominating set $S$ of $G$ is called to {\it restrained} if the
induced subgraph $G[\bar S]$ contains no isolated vertices, where
$\bar S=V(G)\setminus S$. The {\it restrained domination number} of
$G$, denoted by $\gamma^r(G)$, is the smallest cardinality of a
restrained dominating set of $G$. It is clear that $\gamma^r(G)$
exists and $\gamma(G)\leqslant\gamma^r(G)$ for any non-empty graph
$G$.

The concept of restrained domination was introduced by Telle and
Proskurowski~\cite{tp97} in 1997, albeit indirectly, as a vertex
partitioning problem. Here conditions are imposed on a set $S$, the
complementary set $\bar S$ and on edges between the sets $S$ and
$\bar S$. For example, if we require that every vertex in $\bar S$
should be adjacent to some other vertex of $\bar S$ (the condition
on the set $\bar S$) and to some vertex in $S$ (the condition on
edges between the sets $S$ and $\bar S$), then S is a restrained
dominating set. Concerning the study of the restrained domination
numbers, the reader is referred to \cite{dhhm00a, dhhm00b, dhhlm99,
h99, tp97}.

In 2008, Hattingh and Plummer~\cite{hp08} defined the {\it
restrained bondage number} $b^r(G)$ of a nonempty graph $G$ to be
the minimum cardinality among all subsets $B\subseteq E(G)$ for
which $\gamma^r(G-B)>\gamma^r(G)$.

For some simple graphs, their restrained domination numbers can be
easily determined, and so restrained bondage numbers have been also
determined. For example, it is clear that $\gamma^r(K_n)=1$. Domke
{\it et al}.~\cite{dhhlm99} showed that $\gamma^r(C_n)=n-2\lfloor
n/3\rfloor$ for $n\geqslant 3$, and $\gamma^r(P_n)=n-2\lfloor
(n-1)/3\rfloor$ for $n\geqslant 1$. Using this results, Hattingh and
Plummer~\cite{hp08} obtained the restricted domination numbers for
$K_n, C_n, P_n$ and $G=K_{n_1,n_2,\cdots,n_t}$.

\vskip6pt\begin{thm}\label{thm8.1.1} \textnormal{(Hattingh and
Plummer~\cite{hp08}, 2008)}\ For a complete graph $K_n\ (n\geqslant
3)$,
  $$
  b^r(K_n)=\left\{\begin{array}{ll}
  \ 1 \ & {\rm if}\ n=3;\\
  \left\lceil\frac{n}2\right\rceil\ &{\rm otherwise}.
  \end{array}\right.
  $$

For a cycle $C_n$ of order $n\ (\geqslant 3)$,
 $$
 b^r(C_n)=\left\{\begin{array}{ll}
 1\ & {\rm if}\  n\equiv 0\,({\rm mod}\, 3;\\
 2\ & {\rm otherwise}.\end{array}\right.
 $$

For a path $P_n$ of order $n\,(\geqslant 4)$, $b^r(P_n)=1$.

For a complete $t$-partite graph $G=K_{n_1,n_2,\cdots,n_t}$ with
$n_1\leqslant n_2\leqslant \cdots\leqslant n_t\, (t\geqslant 2)$,
%where $n_i\geqslant 2$ for some $1\leqslant i\leqslant t$,
 $$
 b^r(G)=\left\{\begin{array}{ll}
 \left\lceil\frac m2\right\rceil\ & {\rm if}\ n_m=1\ {\rm and}\ n_{m+1}\geqslant 2\, (1\leqslant m<t);\\
 2t-2\ & {\rm if}\ n_1=n_2=\cdots =n_t=2\,(t\geqslant 2);\\
 2 & {\rm if}\ n_1=2\ {\rm and}\ n_2\geqslant 3\,(t=2);\\
 \sum\limits_{i=1}^{t-1}\ n_i-1 & {\rm otherwise}.
 \end{array}\right.
 $$
\end{thm}

\begin{pf} We first show the conclusion for a complete graph $K_n$.

If $n=3$ then $\gamma^r(K_3)=1$ clearly. Now, removing any edge from
$K_3$ yields $P_3$. Since $\gamma^r(P_3)=3$, it follows that
$b^r(K_3)=1$. Let $n\geqslant 4$ and let $H$ be a spanning subgraph
of $K_n$ that is obtained by removing fewer than $\lceil n/2\rceil$
edges from $K_n$. Then H contains a vertex of degree $n-1$.
Moreover, for every $x\in V(H)$, $d_H(x)\geqslant 2$. Hence,
$\gamma^r(H)=1$. It follows that $b^r(K_n)\geqslant \lceil
n/2\rceil$.

Let $H$ be the graph obtained by removing a perfect matching $M$
from $K_n$. If $n$ is even, then $|M|=n/2=\lceil n/2\rceil$. Thus,
for every $x\in V(H)$, $d_H(x)= n-2$, whence $\gamma^r(H)=2$. If $n$
is odd then, then $|M|=(n-1)/2=\lceil n/2\rceil-1$. there is exactly
one vertex $x\in V(H)$ such that $d_{H}(x)=n-1$. Let $H'$ be the
graph obtained by removing from $H$ one edge incident with $x$. It
follows that $\gamma^r(H')=2$.

In either case, we can obtain a graph, by removing $\lceil
n/2\rceil$ edges from $K_n$, whose restrained bondage number is
lager than $\gamma^r(K_n)$. Thus $b^r(K_n)\leqslant \lceil
n/2\rceil$, whence $b^r(K_n)=\lceil n/2\rceil$.

We now show the conclusion for a cycle $C_n$.

Assume $n\equiv 0\,({\rm mod}\, 3)$. Since $\gamma^r(C_n)<\gamma^r
(P_n)$, $b^r(C_n)=1$. Thus, assume $n\equiv i\,({\rm mod}\, 3)$
$(i=1,2)$. Since $\gamma^r(C_n)=\gamma^r(P_n)$, it follows that
$b^r(C_n)\geqslant 2$. Let $H$ be the graph obtained by the removal
of two edges from $C_n$ such that $P_3$ and $P_{n-3}$ are formed.
Then
 $$
 \begin{array}{rl}
 \gamma^r(H)&=\gamma^r(P_{n-3})+\gamma^r(P_3)\\
 &=(\lceil (n-3)/3\rceil +i-1)+3\\
 &=(\lceil n/3-1\rceil+i-1)+3\\
 &=(\lceil n/3\rceil+i-1)+2\\
 &=\gamma^r(C_n)+2\\
 &>\gamma^r(C_n).
 \end{array}
 $$
Thus, $b^r(C_n)\leqslant 2$. Hence, $b^r(C_n)=2$.

We next show the conclusion for a path $P_n$, $n\geqslant 4$.

Assume $n\equiv i\,({\rm mod}\, 3)$ $(i=1,2)$. Since
$\gamma^r(P_n)=\gamma^r(C_n)$, by reasoning similar to that in the
previous proof, we have $b^r(P_n)\leqslant 1$, whence $b^r(P_n)=1$.
Assume $n\equiv 0\,({\rm mod}\, 3)$. Let $H$ be the graph obtained
by the removal of one edge from $P_n$ such that $P_3$ and $P_{n-3}$
are formed. Then
 $$
 \begin{array}{rl}
 \gamma^r(H)&=\gamma^r(P_{n-3})+\gamma^r(P_3)\\
 &=(\lceil (n-3)/3\rceil+2)+3\\
 &=(\lceil n/3-1\rceil+2)+3\\
 &=\lceil n/3\rceil -1+2+3\\
 &=(\lceil n/3\rceil +2)+2\\
 &=\gamma^r(P_n)+2\\
 &>\gamma^r(P_n).
 \end{array}
 $$
Thus, $b^r(P_n)\leqslant 1$. Hence, $b^r(P_n)=1$.

The proof of the last conclusion is little complex and omitted here.
\end{pf}

\vskip6pt\begin{thm}\label{thm8.1.2} \textnormal{(Hattingh and
Plummer~\cite{hp08}, 2008)}\ For a tree $T$ of order $n\,(\geqslant
4)$, $T\ncong S_n$ if and only if $b^r(T)=1$, where $S_n$ is a star
of order $n$.
\end{thm}

Theorem~\ref{thm8.1.2} shows that the restrained bondage number of a
tree can be computed in constant time. However, the decision problem
for $b^r(G)$ is NP-complete even for bipartite graphs.

\vskip6pt\begin{prob} Consider the decision problem:

Restrained Bondage Problem

Instance: A graph $G$ and a positive integer $b$.

Question: Is $b^r(G)\leqslant b$?
\end{prob}

\begin{thm}\label{thm8.1.4} \textnormal{(Hattingh and
Plummer~\cite{hp08}, 2008)}\ The restrained bondage problem is
NP-complete, even for bipartite graphs.
\end{thm}

Consequently, it is significative to establish some sharp bounds of
the restrained bondage number of a graph in terms of some other
graphic parameters.

\vskip6pt\begin{thm}\label{thm8.1.5} \textnormal{(Hattingh and
Plummer~\cite{hp08}, 2008)}\ For any graph $G$ with
$\delta(G)\geqslant 2$,
 $$
 b^r(G)\leqslant \min\limits_{xy\in E(G)}\{d_G(x)+d_G(y)-2\}.
 $$
\end{thm}

\begin{pf}
Let $b^r=\min\{d_G(x) + d_G(y)-2:\ xy\in E(G)\}$, and let $xy\in
E(G)$ such that $d_G(x)+d_G(y)-2=b^r$. Suppose to the contrary that
$b^r(G)>b^r$. Let $E_{xy}$ denote the set of edges that are incident
with at least one of $x$ and $y$, but not both. Then $|E_{xy}|=b^r$
and $\gamma^r(G-E_{xy})=\gamma^r(G)$ since $b^r(G)>b^r$. Since $x$
and $y$ are vertices of degree one in $G-E_{xy}$, it follows that
$\gamma^r(G-x-y)=\gamma^r(G)-2$. Let $R$ be a $\gamma^r$-set of
$G-x-y$ and let $N_{xy}=N_G(x)\cup N_G(y)-\{x,y\}$. Since
$\delta(G)\geqslant 2$, it follows that $N_{xy}\ne\emptyset$. If
$N_{xy}\subseteq R$, then $R$ is a restrained dominating set of $G$
of cardinality $\gamma^r(G-x-y)=\gamma^r(G)-2$, a contradiction.
Hence, $N_{xy}\nsubseteq R$ and there is a vertex $z\in N_{xy}$ such
that $z\notin R$. Without loss of generality, assume $z$ is adjacent
to $x$. Then $R\cup\{y\}$ is a restrained dominating set of $G$ of
cardinality $\gamma^r(G-x-y)+1=\gamma^r(G)-1$, a contradiction.
\end{pf}

%\vskip6pt
\begin{cor}\label{cor8.1.6}
$b^r(G)\leqslant\Delta(G)+\delta(G)-2$ for any graph $G$ with
$\delta(G)\geqslant 2$.
\end{cor}

Notice that the bounds stated in Theorem~\ref{thm8.1.5} and
Corollary~\ref{cor8.1.6} are sharp. Indeed the class of cycles whose
orders are congruent to 1, 2 (mod 3) have a restrained bondage
number achieving these bounds (see Theorem~\ref{thm8.1.1}).

\vskip6pt

{\color{blue}\noindent{\bf Remarks}\ Theorem~\ref{thm8.1.5} is an
analogue of Theorem~\ref{thm3.1.3}. A quite natural problem is
whether or not, for restricted bondage number, there are analogues
of Theorem~\ref{thm3.1.5}, Theorem~\ref{thm3.1.7},
Theorem~\ref{thm3.2.1} and so on. Indeed, we should consider all
results for the bondage number whether or not there are analogues
for restricted bondage number.}

\vskip6pt\begin{thm}\label{thm8.1.7} \textnormal{(Hattingh and
Plummer~\cite{hp08}, 2008)}\ For any graph $G$ with
$\gamma^r(G)\geqslant 2$,
 $$
 b^r(G)\leqslant(\gamma^r(G)-1)\Delta(G) +1.
 $$
\end{thm}

\begin{pf}
We proceed by induction on $\gamma^r(G)$. Let $\gamma^r(G)=2$, and
suppose $b^r(G)\geqslant \Delta(G)+2$. Let $x\in V(G)$ be of maximum
degree. It follows that $\gamma^r(G-x)=\gamma^r(G)-1=1$ and
$b^r(G-x)\geqslant 2$. Since $\gamma^r(G)=2$ and $\gamma^r(G-x)=1$,
there is a vertex $y\in V(G-x)$ that is adjacent to every vertex in
$V(G)-\{x\}$. Furthermore, $x$ is adjacent to every vertex in
$V(G)-\{y\}$. Let $e$ be any edge incident with $y$, and let
$H=(V(G-x), E(G-x-e))$. Since $b^r(G-x)\geqslant 2$, it follows that
$\gamma^r(H)=1$. Hence, there is a vertex $z\in V(G-x)$ such that
$z\ne y$ and $z$ is adjacent to every vertex in $V(G-x)$. Since $y$
is the only vertex not in $N_G(x)$, we have $z\in N_G(x)$. Hence,
$d_G(z)=|V(G)|-1$, a contradiction. Thus,
$b^r(G)\leqslant\Delta(G)+1$, for $\gamma^r(G)=2$.

Now, assume that, for any graph $G'$ such that
$\gamma^r(G')=k\geqslant 2$, $b^r(G') \leqslant (k-1)\Delta(G') +
1$. Let $G$ be a graph such that $\gamma^r(G) = k + 1$. Suppose to
the contrary that $b^r(G) > k\Delta(G) + 1$. Let $x\in V(G)$ and
notice that $\gamma^r(G-x)=\gamma^r(G)-1=k$. Furthermore, $b^r(G)
\leqslant b^r(G-x) + d_G(x)$. By the inductive hypothesis we have
 $$
 \begin{array}{rl}
 b^r(G) &\leqslant [(k-1)\Delta(G-x)+1]+d_G(x)\\
 & \leqslant [(k-1)\Delta(G-x)+1]+\Delta(G)\\
 & = k\Delta(G)+1.
 \end{array}
 $$
Thus, $b^r(G)\leqslant k\Delta(G)+1$, contradicting our assumption
that $b^r(G) > k\Delta(G) + 1$. By induction the proof is complete.
\end{pf}

\vskip6pt

We now consider the relation between $b^r(G)$ and $b(G)$.

\vskip6pt\begin{thm}\label{thm8.1.8} \textnormal{(Hattingh and
Plummer~\cite{hp08}, 2008)}\ If $\gamma^r(G)=\gamma(G)$ for some
graph $G$, then $b^r(G) \leqslant b(G)$.
\end{thm}
%
% \begin{equation}\label{e8.1.1}
% b^r(G) \leqslant b(G)\ \ {\rm if}\ \gamma^r(G)=\gamma(G).
% \end{equation}

\begin{pf}\
Indeed, assume $\gamma^r(G)=\gamma(G)$. Let $B$ be a set of edges
such that $\gamma(G-B)>\gamma(G)$, where $|B|=b(G)$. Then
 $$
 \gamma^r(G)=\gamma (G)< \gamma (G-B) \leqslant \gamma^r (G-B),
 $$
whence $b^r (G) \leqslant |B|=b(G)$.
\end{pf}

\vskip6pt

However, we do not have $b^r(G)=b(G)$ for any graph $G$ even if
$\gamma^r(G)=\gamma (G)$. Observe that $\gamma^r(K_3)=\gamma(K_3)$,
yet $b^r(K_3)=1$ and $b(K_3)=2$. We still may not claim that
$b^r(G)=b(G)$ even in the case that every $\gamma(G)$-set is a
$\gamma^r(G)$-set. The example $K_3$ again demonstrates this.

\vskip6pt\begin{thm}\label{thm8.1.9} \textnormal{(Hattingh and
Plummer~\cite{hp08}, 2008)}\ There is an infinite class of graphs in
which each graph $G$ satisfies $b(G)<b^r(G)$.
\end{thm}

\begin{pf}
%Furthermore, there is an infinite class of graphs in which each
%graph $G$ satisfies $b(G)<b^r(G)$.
Such a graph $G$ can be obtained
from some connected graph $H$ by attaching $\ell\geqslant 2$ pendant
vertices to each vertex in $H$, denoted by $BC(H)$. Let $\mathscr
B=\{G: G=BC(H)$ for some graph $H$ such that $\delta(H)\geqslant
2\}$.

For any $G \in \mathscr B$, there is a graph $H$ such that
$G=BC(H)$, where $\ell$ is number of pendant vertices attached to
each vertex in $H$. Let $L$ denote the set of pendant vertices of
$G$. Clearly, $L$ is the unique $\gamma^r$-set of $G$ and $V(H)$ is
the unique $\gamma$-set of $G$. It follows immediately that $b(G)=1$
by Theorem~\ref{thm2.2}, and $b^r(G)=\min\{\delta(H),
\ell\}\geqslant 2$. This fact shows that $b(G)< b^r (G)$ for any $G
\in \mathscr B$.
% \begin{equation}\label{e8.1.2}
% b(G)< b^r (G) \ {\rm for\ any }\ G \in \mathscr B.
% \end{equation}
 \end{pf}

\vskip6pt

In fact, there exists a graph $G \in \mathscr B$ such that $b^r(G)$
can be much larger than $b(G)$, which is stated as the following
theorem.

\vskip6pt\begin{thm}\label{thm8.1.10} \textnormal{(Hattingh and
Plummer~\cite{hp08}, 2008)}\ For each positive integer $k$ there is
a graph $G$ such that $k=b^r (G)-b(G)$.
\end{thm}

Combining Theorem~\ref{thm8.1.8} with Theorem~\ref{thm8.1.9}, we
have the following corollary.

{\color{blue}\vskip6pt\begin{cor}\label{cor8.1.11}
% \textnormal{(Hattingh and Plummer~\cite{hp08}, 2008)}
\ The bondage number and the restrained bondage number are
unrelated.
\end{cor}}

{\color{blue}\noindent{\bf Comments}\ Corollary~\ref{cor8.1.6} and
Theorem~\ref{thm8.1.8} provides a possibility to attack
Conjecture~\ref{con6.1} by characterizing planar graphs with
$\gamma^r=\gamma$ and $b^r=b$ when $3\leqslant\Delta\leqslant 6$. }

\vskip16pt

\subsection{Total Bondage Numbers}

A dominating set $S$ of a graph $G$ without isolated vertices is
called to be {\it total} if the induced subgraph $G[S]$ contains no
isolated vertices. The minimum cardinality of a total dominating set
is called the {\it total domination number} of $G$ and denoted by
$\gamma^t(G)$. It is clear that $\gamma(G)\leqslant\gamma^t(G)
\leqslant 2\gamma(G)$ for any graph $G$ without isolated vertices.

The total domination in graphs was introduced by Cockayne {\it et
al.}~\cite{cdh80} in 1980. Pfaff {\it et at.}~\cite{plh83, lphh84}
in 1983 showed that the problem determining total domination number
for general graphs is NP-complete, even for bipartite graphs, and
chordal graphs. Even now, total domination in graphs has been
extensively studied in the literature. In 2009, Henning~\cite{h09}
gave a survey of selected recent results on total domination in
graphs.

The {\it total bondage number} of $G$, denoted by $b^t(G)$, is the
smallest cardinality of a subset $B\subseteq E(G)$ with the property
that $G-B$ contains no isolated vertices and
$\gamma^t(G-B)>\gamma^t(G)$.

From definition, $b^t(G)$ may not exist for some graphs, for
example, $G=K_{1,n}$. We put $b^t(G)=\infty$ if $b^t(G)$ does not
exist. In fact, $b^t(G)$ is finite for any connected graph $G$ other
than $K_1, K_2, K_3$ and $K_{1,n}$. Since there is a path of length
$3$ in $G$, we can find $B_1\subseteq E(G)$ such that $G-B_1$ is a
spanning tree $T$ of $G$, containing a path of length $3$. So
$\gamma^t(T)=\gamma^t(G-B_1)\geqslant\gamma^t(G)$. For the tree $T$
we find $B_2\subseteq E(T)$ such that
$\gamma^t(T-B_2)>\gamma^t(T)\geqslant\gamma^t(G)$. Thus, we have
$\gamma^t(G-B_2-B_1)>\gamma^t(G)$, and so
$b^t(G)\leqslant|B_1|+|B_2|$.

In 1991, Kulli and Patwari~\cite{kp91} first studied the total
bondage number of a graph and calculated the exact values of
$b^t(G)$ for some standard graphs.

\begin{thm}\label{thm8.12}\textnormal{(Kulli and Patwari~\cite{kp91}, 1991)}
For a cycle $C_n$ and a path $P_n$, $n\geqslant 4$,
 $$
 \begin{array}{rl}
  & b^t(C_n)=\left\{\begin{array}{ll}
3&{\rm if}\ n\equiv 2\,({\rm mod}\,4);\\
2&{\rm otherwise},
\end{array}\right.
\quad {\rm and}\quad \\
 & b^t(P_n)=\left\{\begin{array}{ll}
2&{\rm if}\ n\equiv 2\,({\rm mod}\,4);\\
1&{\rm otherwise}.
\end{array}\right.
\end{array}
 $$
For a complete bipartite graph $K_{m,n}$ and a complete $K_n$,
 $$
 \begin{array}{rl}
 & b^t(K_{m,n})=m \ {\rm with}\ 2\leqslant m\leqslant n;\\
 &
 b^t(K_n)=\left\{\begin{array}{ll}
 4\ & {\rm for}\ n=4;\\
 2n-5\ & {\rm for}\ n\geqslant 5.
 \end{array}\right.
 \end{array}
 $$
\end{thm}

Recently, Hu, Lu and Xu~\cite{hlx09} have obtained some results on
the total bondage number of the Cartesian product $P_m\times P_n$ of
two paths $P_m$ and $P_n$.

\begin{thm}\label{thm8.13}\textnormal{(Hu, Lu and Xu~\cite{hlx09}, 2009)}
For the Cartesian product $P_m\times P_n$ of two paths $P_m$ and
$P_n$,
 $$
 \begin{array}{rl}
 & b^t(P_2\times P_n)=\left\{ \begin{array}{l}
 1 \ {\rm if}\ n\equiv 0\,({\rm mod}\, 3),  \\
2 \ {\rm if}\ n\equiv 2\,({\rm mod}\, 3),\\
3 \ {\rm if}\ n\equiv 1\,({\rm mod}\, 3);
\end{array}
 \right.\\
 & b^t(P_3\times P_n)=1;\\
 & b^t(P_4\times P_n)\left\{ \begin{array}{ll}
 =1\ & {\rm if}\  n\equiv 1\,({\rm mod}\, 5);\\
 =2\ & {\rm if}\  n\equiv 4\,({\rm mod}\, 5);\\
 \leqslant 3\ & {\rm if}\  n\equiv 2\,({\rm mod}\, 5);\\
 \leqslant 4\ & {\rm if}\  n\equiv 0,3\,({\rm mod}\, 5).\\
 \end{array} \right.
 \end{array}
 $$
 \end{thm}

Generalized Petersen graphs are an important class of commonly used
interconnection networks and have been studied recently. By
constructing a family of minimum total dominating sets, Cao {\it et
al.}~\cite{cym10} determined the total bondage number of the
Generalized Petersen graphs.

From Theorem~\ref{thm2.3}, we know that $b(T)\leqslant 2$ for a
nontrivial tree $T$. But given any position integer $k$, Sridharan
{\it et al.}~\cite{ses07} constructed a tree $T$ for which
$b^t(T)=k$. Let $H_k$ be the tree obtained from the star $K_{1,k+1}$
by subdividing $k$ edges twice. The tree $H_7$ is shown in
Figure~\ref{f9}. It can be easily verified that $b^t(H_k)=k$. We
state this fact as the following theorem.

\begin{figure}[ht]  % h: here; t: top; b: bottom
%\psset{unit=.9}     % scale
\psset{arrowsize=.14}

\begin{center}
\begin{pspicture}(-3.5,0)(3.5,4.5)

\cnode(0,0){3pt}{1} \cnode(0,1){3pt}{2} \ncline{1}{2}
\cnode(-1,1){3pt}{11} \cnode(-2,1){3pt}{12} \cnode(-3,1){3pt}{13}
\ncline{2}{11} \ncline{11}{12} \ncline{12}{13}

\cnode(-.866,1.5){3pt}{21} \cnode(-1.732,2){3pt}{22}
\cnode(-2.589,2.5){3pt}{23} \ncline{2}{21} \ncline{21}{22}
\ncline{22}{23}

\cnode(-.5,1.866){3pt}{31} \cnode(-1,2.732){3pt}{32}
\cnode(-1.5,3.589){3pt}{33} \ncline{2}{31} \ncline{31}{32}
\ncline{32}{33}

\cnode(0,2){3pt}{41} \cnode(0,3){3pt}{42} \cnode(0,4){3pt}{43}
\ncline{2}{41} \ncline{41}{42} \ncline{42}{43}

\cnode(1,1){3pt}{51} \cnode(2,1){3pt}{52} \cnode(3,1){3pt}{53}
\ncline{2}{51} \ncline{51}{52} \ncline{52}{53}

\cnode(.866,1.5){3pt}{61} \cnode(1.732,2){3pt}{62}
\cnode(2.589,2.5){3pt}{63} \ncline{2}{61} \ncline{61}{62}
\ncline{62}{63}

\cnode(.5,1.866){3pt}{71} \cnode(1,2.732){3pt}{72}
\cnode(1.5,3.589){3pt}{73} \ncline{2}{71} \ncline{71}{72}
\ncline{72}{73}

\end{pspicture}

\caption{\label{f9}\footnotesize $H_7$.}
\end{center}

\end{figure}

%\vskip6pt
\begin{thm}\label{thm8.22a} \textnormal{(Sridharan
{\it et al.}~\cite{ses07}, 2007)}\ For any positive integer $k$,
there exists a tree $T$ with $b^t(T)=k$.
\end{thm}

Combining Theorem~\ref{thm2.1.1} with Theorem~\ref{thm8.22a}, we
have what follows.

{\color{blue}\vskip6pt\begin{cor}\label{cor8.23a}
\textnormal{(Sridharan {\it et al.}~\cite{ses07}, 2007)}\ The
bondage number and the total bondage number are unrelated, even for
trees.
\end{cor}}

However, Sridharan {\it et al.}~\cite{ses07} gave an upper of the
total bondage number of a tree in terms of its maximum degree: For
any tree $T$ of order $n$, if $T\ne =K_{1,n-1}$, then
$b^t(T)=\min\{\Delta(T), \frac 13(n-1)\}$. Rad and
Raczek~\cite{rr11} improved this upper and gave a constructive
characterization of a certain class of trees attaching the upper
bound.

\begin{thm}\label{thm8.2.5}{\rm (Rad and Raczek~\cite{rr11}, 2011)}
$b^t(T)\leqslant\Delta (T)-1$ for any tree $T$ with maximum degree
at least three.
\end{thm}

In general, the decision problem for $b^t(G)$ is NP-complete for any
graph $G$. We state the decision problem for the total bondage as
follows.

\vskip6pt\begin{prob} Consider the decision problem:

Total Bondage Problem

Instance: A graph $G$ and a positive integer $b$.

Question: Is $b^t(G)\leqslant b$?
\end{prob}

\begin{thm}\label{thm8.15}
\textnormal{(Hu and Xu~\cite{hx12}, 2012)} The total bondage problem
is NP-complete.
\end{thm}

Consequently, it is significative to establish
% the exact value of the total bondage number for some special graphs, and
some sharp
bounds of the total bondage number of a graph in terms of some other
graphic parameters.

\begin{thm}\label{thm8.16}\textnormal{(Kulli and Patwari~\cite{kp91}, 1991)}
$b^t(G)\leqslant 2\upsilon -5$ for any graph $G$ with order
$\upsilon\geqslant 5$.
\end{thm}

In 2007, Sridharan {\it et al.}~\cite{ses07} improved this result as
follows.

\begin{thm}\label{thm8.17}
%\textnormal{(Sridharan, Elias and Subramanian~\cite{ses07}, 2007)}
Let $G$ be a graph with order $\upsilon\geqslant 5$. Then
 $$
 b^t(G)\leqslant\left\{\begin{array}{ll}
  \frac 13 (\upsilon -1)\ &{\rm if}\ G\ {\rm contains\ no\
  cycles};\\
  \upsilon -1\ & {\rm if}\ g(G)\geqslant 5;\\
  \upsilon -2\ & {\rm if}\ g(G)=4.
  \end{array}\right.
  $$
\end{thm}

Rad and Raczek~\cite{rr11} also established some upper bounds of
$b^t(G)$ for a general graph $G$. In particular, they gave an upper
bound of $b^t(G)$ in terms of the girth of a graph.

\begin{thm}\label{thm8.2.10}{\rm (Rad and Raczek~\cite{rr11}, 2011)}
$b^t(G)\leqslant 4\Delta (G)-5$ for a graph $G$ with $g(G)\geqslant
4$.
\end{thm}

\vskip20pt

\subsection{Paired Bondage Numbers}

A dominating set $S$ of $G$ is called to be {\it paired} if the
induced subgraph $G[S]$ contains a perfect matching. The {\it paired
domination number} of $G$, denoted by $\gamma^p(G)$, is the minimum
cardinality of a paired dominating set of $G$. Clearly,
$\gamma^t(G)\leqslant\gamma^p(G)$ for every connected graph $G$ with
order at least two, where $\gamma^t(G)$ is the total domination
number of $G$, and $\gamma^p(G)\leqslant 2\gamma(G)$ for any graph
$G$ without isolated vertices. Paired domination was introduced by
Haynes and Slater~\cite{hs95, hs98}, and further studied in
\cite{h08, phs01, qkcd03, skh04}.

%\bibitem{r08}%Joanna Raczek
%J. Raczek, Paired bondage in trees. {\it Discrete Mathematics}, {\bf
%308} (2008), 5570-5575.

%\bibitem{wf04}
%[13] Y. Wu, Q. Fan, The bondage number of four domination parameters
%for trees, J. Math., Wuchan Univ. 24 (2004) 267¨C270.

The {\it paired bondage number} of $G$ with $\delta(G)\geqslant 1$,
denoted by $b^p(G)$, is the minimum cardinality among all sets of
edges $B\subseteq E$ such that $\delta(G-B)\geqslant 1$ and
$\gamma^p(G-B)>\gamma^p(G)$.

%We say that $G$ is a $p$-strongly stable graph if, for all $B\subseteq E$,
%either $\gamma^p(G-B)=\gamma^p(G)$ or $\delta(G-B)=0$, and write
%$b^p(G)=0$.

The concept of the paired bondage number was first proposed by
Raczek~\cite{r08} in 2008. The following observations follow
immediately from the definition of the paired bondage number.

\vskip6pt\begin{prop}\ Let $G$ be a graph with $\delta(G)\geqslant
1$.

(a)\ If $H$ is a subgraph of $G$ such that $\delta(H)\geqslant 1$,
then $b^p(H)\leqslant b^p(G)$.

(b)\ If $H$ is a subgraph of $G$ such that $b^p(H)=1$ and $k$ is the
number of edges removed to form $H$, then $1\leqslant
b^p(G)\leqslant k+1$.

(c)\ If $xy\in E(G)$ such that $d_G(x), d_G(y)>1$, and $xy$ belongs
to each perfect matching of each minimum paired dominating set of
$G$, then $b^p(G)=1$.
\end{prop}

Based on these observations, $b^p(P_n)\leqslant 2$. In fact, the
paired bondage number of a path has been determined.

\vskip6pt\begin{thm}\label{thm8.19} \textnormal{(Raczek~\cite{r08},
2008)}\ Let $P_n$ be a path of order $n\geqslant 2$, and let $k$ be
a positive integer. Then
 $$
 b^p(P_n)=\left\{\begin{array}{ll}
 0  & {\rm if}\ n=2,3\ {\rm or}\ 5;\\
 1  & {\rm if}\ n=4k, 4k + 3\ {\rm or}\ 4k + 6;\\
 2 & {\rm otherwise}.
 \end{array}
 \right.
 $$
\end{thm}

Foe a cycle $C_n$ of order $n\geqslant 3$, since
$\gamma^p(C_n)=\gamma^p(P_n)$, from Theorem~\ref{thm8.19} we obtain
the following result.

\vskip6pt\begin{cor}\textnormal{ (Raczek~\cite{r08}, 2008)}\ Let
$C_n$ be a cycle of order $n\geqslant 3$, and let $k$ be a positive
integer. Then
 $$
 b^p(C_n)=\left\{\begin{array}{ll}
 0  & {\rm if}\ n=3\ {\rm or}\ 5;\\
 2  & {\rm if}\ n=4k, 4k + 3\ {\rm or}\ 4k + 6;\\
 3 & {\rm otherwise}.
 \end{array}
 \right.
 $$
\end{cor}

A wheel $W_n$, where $n\geqslant 4$, is a graph with $n$ vertices,
formed by connecting a single vertex to all vertices of a cycle
$C_{n-1}$. Of course $\gamma^p(W_n)=2$.

\vskip6pt\begin{thm}\label{thm8.21}\textnormal{ (Raczek~\cite{r08},
2008)}\ For a complete bipartite graph $K_{m,n}$, where
$1<m\leqslant n$, $b^p(K_{m,n} )=m$. For a wheel $W_n$,
 $$
 b^p(W_n)=\left\{\begin{array}{ll}
 4  & {\rm if}\ n=4;\\
 3  & {\rm if}\ n=5;\\
 2 & {\rm otherwise}.
 \end{array}
 \right.
 $$
\end{thm}

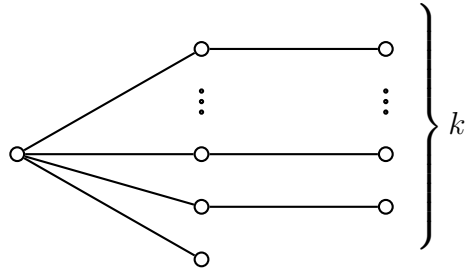
\begin{figure}[ht]  % h: here; t: top; b: bottom
\psset{unit=.7}     % scale
\psset{arrowsize=.14}

\begin{center}
\begin{pspicture}(-4,0)(5,4)
\cnode(-3.5,2){3pt}{1} \cnode(0,0){3pt}{2} \cnode(0,1){3pt}{3}
\cnode(0,2){3pt}{4} \cnode(0,4){3pt}{5} \cnode(3.5,1){3pt}{6}
\cnode(3.5,2){3pt}{7} \cnode(3.5,4){3pt}{8} \cnode(0,2.8){1pt}{9}
\cnode(0,3){1pt}{10} \cnode(0,3.2){1pt}{11} \cnode(3.5,2.8){1pt}{12}
\cnode(3.5,3){1pt}{13} \cnode(3.5,3.2){1pt}{14} \ncline{1}{2}
\ncline{1}{3} \ncline{1}{4} \ncline{1}{5} \ncline{3}{6}
\ncline{4}{7} \ncline{5}{8}

\put(4,2.38){$\left.\rule{0mm}{18mm}\right\}k$}
\end{pspicture}

\caption{\label{f10} \footnotesize A tree $T$ with $b_p(T)=k$.}
\end{center}

\end{figure}

Theorem~\ref{thm3.1.5} shows that, if $T$ is a non-trivial tree,
then $b(T)\leqslant 2$. However, no similar result exists for paired
bondage. For any non-negative integer $k$, let $T_k$ be a tree
obtained by subdividing all but one edge of the star $K_{1,k+1}$ (as
shown in Figure~\ref{f10}). It is easy to see that $b^p(T_k)=k$. We
state this fact as the following theorem.

\vskip15pt

%We state this fact as the following theorem.

%\vskip6pt
\begin{thm}\label{thm8.22} \textnormal{(Raczek~\cite{r08}, 2008)}\ For
any non-negative integer $k$, there exists a tree $T$ with
$b^p(T)=k$.
\end{thm}

Consider the tree defined in Figure~\ref{f6} for $k=3$. Then
$b(T_3)\leqslant 2$ and $b^p(T_3)=3$. On the other hand, $b(P_4)=2$
and $b^p(P_4)=1$. Thus, we have what follows.

{\color{blue}\vskip6pt\begin{cor}\label{cor8.23}
%\textnormal{(Raczek~\cite{r08}, 2008)}
\ The bondage number and the
paired bondage number are unrelated, even for trees.
\end{cor}}

A constructive characterization of trees with $b(T)=2$ is given by
Hartnell and Rall in \cite{hr92}. Raczek~\cite{r08} provided a
constructive characterization of trees with $b^p(T)=0$. In order to
state the characterization, we define a labeling and three simple
operations on a tree $T$. Let $y\in V(T)$ and let $\ell(y)$ be the
label assigned to $y$.

{\it Operation} $\mathscr{T}_1$. If $\ell(y)=B$, add a vertex $x$
and the edge $yx$, and let $\ell(x)=A$.

{\it Operation} $\mathscr{T}_2$. If $\ell(y)=C$, add a path $(x_1,
x_2)$ and the edge $yx_1$, and let $\ell(x_1)=B$, and $\ell(x_2)=A$.

{\it Operation} $\mathscr{T}_3$. If $\ell(y)=B$, add a path ($x_1,
x_2, x_3)$ and the edge $yx_1$, and let $\ell(x_1)=C, \ell(x_2)=B$
and $\ell(x_3)=A$.

Let $P_2=(u, v)$ with $\ell(u)=A$ and $\ell(v)=B$. Let$\mathscr T$
be the class of all trees obtained from the labeled $P_2$ by a
finite sequence of Operations $\mathscr{T}_1$, $\mathscr{T}_2$,
$\mathscr{T}_3$.

%We will show that T 2 T if and only if bp(T )=0.

\begin{figure}[h]  % h: here; t: top; b: bottom
\psset{unit=.8}     % scale
\psset{arrowsize=.14}

\begin{center}
\begin{pspicture}(-4,.0)(5,3.5)
\cnode(-3.5,0){3pt}{1}\rput(-3.8,0.5){\scriptsize$A$}
\cnode(-3.5,1.3){3pt}{2}\rput(-3.8,1.3){\scriptsize$A$}
\cnode(-3.5,2.6){3pt}{3}\rput(-3.8,2.6){\scriptsize$A$}
\cnode(-2.2,1.3){3pt}{4}\rput(-2,1){\scriptsize$B$}
\cnode(-0.9,1.3){3pt}{5}\rput(-0.9,0.9){\scriptsize$C$}
\cnode(0.2,1.3){3pt}{6}\rput(0.2,1.7){\scriptsize$B$}
\cnode(1.5,1.3){3pt}{7}\rput(1.5,1.7){\scriptsize$C$}
\cnode(2.8,1.3){3pt}{8}\rput(2.8,1.7){\scriptsize$B$}
\cnode(4.1,1.3){3pt}{9}\rput(4.4,1.3){\scriptsize$A$}
\cnode(1.5,0){3pt}{10}\rput(1.8,0){\scriptsize$A$}
\cnode(0.2,2.6){3pt}{11}\rput(.2,3.0){\scriptsize$B$}
\cnode(1.5,3.4){3pt}{12}\rput(1.9,3.6){\scriptsize$A$}

\ncline{1}{4} \ncline{2}{4} \ncline{3}{4} \ncline{4}{5}
\ncline{5}{6} \ncline{5}{11} \ncline{6}{7} \ncline{6}{10}
\ncline{7}{8} \ncline{8}{9} \ncline{11}{12}
\end{pspicture}
\caption{\label{f11} \footnotesize A tree $T$ belong to the family
$\mathscr{T}$.}
\end{center}

\end{figure}
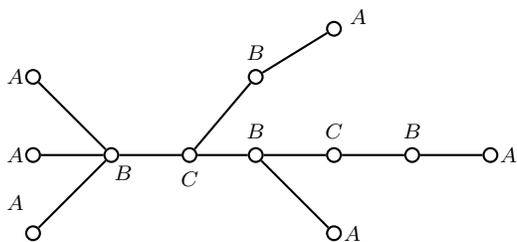

A tree $T$ in Figure~\ref{f11} belongs to the family $\mathscr{T}$.

%\vskip20pt

Raczek~\cite{r08} obtained the following characterization of all
trees $T$ with $b^p(T)=0$.

\vskip6pt\begin{thm}\label{thm8.24} \textnormal{(Raczek~\cite{r08},
2008)}\ Let $T$ be a tree. Then $b^p(T)=0$ if and only if $T$ is in
$\mathscr{T}$.
\end{thm}

We state the decision problem for the paired bondage as follows.

\vskip6pt\begin{prob}\label{prob8.25} Consider the decision problem:

Paired Bondage Problem

Instance: A graph $G$ and a positive integer $b$.

Question: Is $b^p(G)\leqslant b$?
\end{prob}

\begin{con}\label{con8.26}
 The paired bondage problem is NP-complete.
\end{con}

\vskip6pt

\subsection{Independence Bondage Numbers}

A subset $I\subseteq V(G)$ is called an {\it independent set} if no
two vertices in $I$ are adjacent in $G$. The maximum cardinality
among all independent sets is called the {\it independence number}
of $G$, denoted by $\alpha(G)$.

A dominating set $S$ of a graph $G$ is called to be {\it
independent} if $S$ is an independent set of $G$. The minimum
cardinality among all independent dominating set is called the {\it
independence domination number} of $G$ and denoted by $\gamma^i(G)$.

Since an independent dominating set is not only a dominating set but
also an independent set, $\gamma(G)\leqslant\gamma^i(G)\leqslant
\alpha(G)$ for any graph $G$.

It is clear that a maximal independent set is certainly a dominating
set. Thus, an independent set is maximal if and only if it is an
independent dominating set, and so $\gamma^i(G)$ is the minimum
cardinality among all maximal independent sets of $G$. This
graph-theoretical invariant has been well studied in the literature,
see for example Haynes, Hedetniemi and Slater~\cite{hhs98}.
%[T. W. Haynes, S. T. Hedetniemi, P. J. Slater, Fundamentals of Domination in Graphs, Marcel Dekker, New York, 1998].

%Õž²»ª, Áõº£Áú, ËïÁ¼ £¨±±¾©¹¤Òµ´óѧ£©

In 2003, Zhang, Liu and Sun~\cite{zls03} defined the {\it
independence bondage number} $b^i(G)$ of a nonempty graph $G$ to be
the minimum cardinality among all subsets $B\subseteq E(G)$ for
which $\gamma^i(G-B)>\gamma^i(G)$. For some ordinary graphs, their
independence domination numbers can be easily determined, and so
independence bondage numbers have been also determined. Clearly,
$b^i(K_n)=1$ if $n\geqslant 2$.

\vskip6pt\begin{thm}\label{thm8.3b} \textnormal{(Zhang, Liu and
Sun~\cite{zls03}, 2003)}\ For a cycle $C_n$ and a path $P_n$,
$n\geqslant 4$,
 $$
 \begin{array}{rl}
  & b^i(C_n)=\left\{\begin{array}{ll}
1&{\rm if}\ n\equiv 1\,({\rm mod}\,2);\\
2&{\rm if}\ n\equiv 0\,({\rm mod}\,2),
\end{array}\right.
\quad {\rm and}\quad \\
 & b^i(P_n)=\left\{\begin{array}{ll}
1&{\rm if}\ n\equiv 0\,({\rm mod}\,2);\\
2& {\rm if}\ n\equiv 1\,({\rm mod}\,2).
\end{array}\right.
\end{array}
 $$
For a complete bipartite graph $K_{m,n}$,
$b^i(K_{m,n})=\max\{m,n\}$.
\end{thm}

\begin{pf}
We will give the proof of an assertion, say for $P_n$ to show a
basic method. Let $P_n=(x_1,x_2,\ldots, x_n)$ be a path. Clearly,
$\gamma^i(P_n)=\lceil\frac n2\rceil$. We compute $b^i(P_n)$
according as $n$ is even or odd.

If $n$ is even, then $\gamma^i(P_n)=\lceil\frac n2\rceil=\frac n2$,
and
 $$
 \begin{array}{rl}
 \gamma^i(P_n-x_1x_2)=1+\left\lceil\frac {n-1}2\right\rceil=1+\frac
 n2>\gamma^i(P_n).
 \end{array}
 $$
Thus, $b^i(P_n)=1$.

If $n$ is odd, then $\gamma^i(P_n)=\lceil\frac n2\rceil=\frac
{n+1}2$. Let $e$ be any edge in $P_n$. Then the edge $e$ partitions
$P_n$ into two subpaths $P_{k}$ and $P_{\ell}$, where $k+\ell=n$.
Since $n$ is odd, $k$ and $\ell$ are of have different parity.
Without loss of generality, let $k$ be even and $\ell$ odd. Then
 $$
  \begin{array}{rl}
 \gamma^i(P_n-e)=\gamma^i(P_k)+\gamma^i(P_\ell)=\frac
 k2+\frac{\ell+1}2=\frac{n+1}2=\gamma^i(P_n),
 \end{array}
 $$
which implies $b^i(P_n)>1$. On the other hand,
 $$
  \begin{array}{rl}
 \gamma^i(P_n-x_1x_2-x_2x_3)=2+\left\lceil\frac {n-2}2\right\rceil=2+\frac
 {n-2}2=1+\frac{n+1}2>\gamma^i(P_n),
 \end{array}
 $$
which implies $b^i(P_n)\leqslant 2$. It follows that $b^i(P_n)=2$ if
$n$ is odd.
\end{pf}

\vskip6pt

{\color{blue}\noindent{\bf Comments}\ Apart from the above-mentioned
results, as far as we know, there are no other results on the
independence bondage number. We never so much as know any result
on this parameter for a tree.}%(never so much as know-ÎÒÃÇÉõÖÁ²»ÖªµÀ)

\vskip6pt

\vskip6pt

\section{Generalized Types of Bondage Numbers}

There are various generalizations of the classical domination, such
as distance domination, fractional domination and so on. Every such
a generalization can lead to a corresponding bondage. In this
section, we introduce some of them.

\vskip6pt

\subsection{$p$-Bondage Numbers}

In 1985, Fink and Jacobson~\cite{fj85} introduced the concept of
$p$-domination. Let $p$ be a positive integer. A subset $S$ of
$V(G)$ is a {\it $p$-dominating set} of $G$ if $|S\cap
N_G(y)|\geqslant p$ for every $y\in \bar S$. The {\it $p$-domination
number} $\gamma_p(G)$ is the minimum cardinality among all
$p$-dominating sets of $G$. Any $p$-dominating set of $G$ with
cardinality $\gamma_p(G)$ will be called a $\gamma_p$-set of $G$.
Note that the $\gamma_1$-set is the classic minimum dominating set.
%For any $S, T\subseteq V(G)$, $S$ $p$-dominates $T$ in $G$ if, for every $v\in
%T-S$, $|S\cap N_G(v)|\geq p$.
Notice that every graph has a $p$-dominating set since the vertex
set $V(G)$ is such a set. We also note that the 1-dominating set is
a dominating set, and so $\gamma(G)=\gamma_1(G)$. The $p$-domination
number has received much research attention, see a state-of-the-art
survey articles by Chellali {\it et. al.}~\cite{cfhv12}.
% (see, for example, \cite{bcf05, bcv06, f85, fj85, sv93}).

It is clear from definition that every $p$-dominating set of a graph
certainly contains all vertices of degree at most $p-1$. By this
simple observation, to avoid happening the trivial case, we always
assume $\Delta (G)\geqslant p$. For $p\geqslant 2$, Lu {\it et
al}.~\cite{lhxl10} gave a constructive characterization of trees
with unique minimum $p$-dominating sets.

Recently, Lu and Xu~\cite{lx11} have introduced the concept to the
{\it $p$-bondage number} of $G$, denoted by $b_p(G)$, as the minimum
cardinality among all sets of edges $B\subseteq E(G)$ such that
$\gamma_p(G-B)> \gamma_p(G)$. Clearly, $b_1(G)=b(G)$.

%Blidia, Chellali and Volkmann~\cite{bcv06} established some bounds
%on the $p$-domination number for trees.
Lu and Xu~\cite{lx11} established a lower bound and an upper bound
of $\gamma_p(T)$ for any integer $p\geqslant 2$ and any tree $T$
with $\Delta(T)\geqslant p$, and characterized all trees achieving
the lower bound and the upper bound, respectively.

\begin{thm}\label{thm9.3a}\textnormal{(Lu and Xu~\cite{lx11}, 2011)} For
any integer $p\geqslant 2$ and any tree $T$ with $\Delta(T)\geqslant
p$,
 $$1 \leqslant b_p(T)\leqslant \Delta(T)-p+1.$$
\end{thm}

Let $S$ be a given subset of $V(G)$ and, for any $x\in S$, let
 $$
 N_p(x,S,G)=\{y\in \overline S\cap N_G(x):\ |N_G(y)\cap S|=p\}.
 $$
Then, all trees achieving the lower bound $1$ can be characterized
as follows.

%Let $G$ be a graph and $D$ be a subset of $V(G)$. For any $x\in D$,
%a vertex $y$ not in $D$ is said to be a {\it $p$-private neighbor}
%of $x$ with respect to $D$ if $y$ is a neighbor of $x$ and $|D\cap
%N_G(y)|=p$. The {\it $p$-private neighborhood} of $x$ with respect
%to $D$, denoted by $PN_p(x,D,G)$, is the set of all $p$-private
%neighbors of $x$ with respect to $D$ in $G$.

\begin{thm}\label{thm9.3b}\textnormal{(Lu and Xu~\cite{lx11}, 2011)}
Let $T$ be a tree with $\Delta(T)\geqslant p\geqslant 2$. Then
$b_p(T)=1$ if and only if for any $\gamma_p$-set $S$ of $T$ there
exists an edge $xy\in (S,\overline S)$ such that $y\in N_p(x,S,T)$
% or $y\in S$ and $x\in N_p(y,S,T)$.
\end{thm}

The symbol $S(a,b)$ denotes the {\it double star} obtained by adding
an edge between the central vertices of two stars $K_{1,a-1}$ and
$K_{1,b-1}$. And the vertex with degree $a$ (resp., $b$) in $S(a,b)$
is called the {\it $L$-central vertex} (resp., {\it $R$-central
vertex}) of $S(a,b)$.

To characterize all trees attaining the upper bound given in
Theorem~\ref{thm9.3a}, we define three types of operations on a tree
$T$ with $\Delta(T)=\Delta\geqslant p+1$.

\vskip6pt

\noindent
 $$
 \begin{array}{rl}
\mbox{\emph{{Type 1:}}}
    &\mbox{Attach a pendant edge to a vertex $y$ with $d_T(y)\leqslant p-2$ in $T$.}\\
\mbox{\emph{{Type 2:}}}
    &\mbox{Attach a star $K_{1,\Delta-1}$ to a vertex $y$ of $T$ by joining its central vertex to $y$,}\\
    &\mbox{where $y$ in a $\gamma_p$-set of $T$ and $d_T(y)\leqslant \Delta-1$.}\\
\mbox{\emph{{Type 3:}}}
    &\mbox{Attach a double star $S(p,\Delta-1)$
                                 to a pendant vertex $y$ of $T$ by coinciding its}\\
    &\mbox{$R$-central vertex with $y$, where the unique neighbor of $y$ is in a $\gamma_p$-set of $T$.}
\end{array}
 $$

Let $\mathscr{B}=\{T :\ T$ is a tree that can be obtained from
$K_{1,\Delta}$ or $S(p,\Delta)$ by a finite sequence of operations
of Type $1,2,3\}$.

\begin{thm}\label{thm9.4a} \textnormal{(Lu and Xu~\cite{lx11}, 2011)}
A tree with the maximum $\Delta\geqslant p+1$ has $p$-bondage number
$\Delta-p+1$ if and only if it belongs to $\mathscr B$.
\end{thm}

\begin{thm}\label{t5} {\rm (Lu and Xu~\cite{lx12}, 2012)}
Let $G_{m,n}=P_m\times P_n$. Then
 $$
 \begin{array}{rl}
 & b_2(G_{2,n})=1\ \  for\ n\geqslant 2,\\
 & b_2(G_{3,n})=\left\{\begin{array}{ll}
 2& \ { if}\ n\equiv 1\, { (mod\ 3)}; \\
 1& \ { otherwise.}\\
 \end{array}\right. \ \ for\ n\geqslant 2, \\
 & b_2(G_{4,n})=\left\{\begin{array}{ll}
 1& \ { if}\ n\equiv 3\, { (mod\ 4)};\\
 2& \ { otherwise.}
 \end{array} \right. \ \ for\ n\geqslant 7.
 \end{array}
  $$
\end{thm}

\vskip6pt

\subsection{Distance Bondage Numbers}

A subset $S$ of vertices of a graph $G$ is said to be a {\it
distance $k$-dominating set}
%, in short, a {\it distance $k$-dominating set}
for $G$ if every vertex in $G$ not in $S$ is at distance at most $k$
from some vertex of $S$. The minimum cardinality of all distance
$k$-dominating sets is called the {\it distance $k$-domination
number} of $G$ and denoted by $\gamma_k(G)$ (does not confuse with
above-mentioned $\gamma_p(G)$!). When $k=1$, a distance
$1$-dominating set is a normal dominating set, and so
$\gamma_1(G)=\gamma(G)$ for any graph $G$. Thus, the distance
$k$-domination is a generalization of the classical domination.

A subset $I\subseteq V(G)$ is called a {\it distance $k$-independent
set} if $d_G(x,y)>k$ for any two distinct vertices $x$ and $y$ in
$I$. When $k=1$, a distance $1$-independent set is a classical
independent set. The maximum cardinality among all distance
$k$-independent sets is called the {\it distance $k$-independence
number} of $G$, denoted by $\alpha_k(G)$.
% A distance $k$-independent set $I$ is called an {\it $\alpha_k$-set} if
%$|I|=\alpha_k(G)$.
The relation between $\gamma_k$ and $\alpha_{k}$ for a tree obtained
by Meir and Moon~\cite{mm75}, who proved that
$\gamma_k(T)=\alpha_{2k}(T)$ for any tree $T$. The further research
results can be found in Henning {\it et al}~\cite{hos91}, Tian and
Xu~\cite{tx04, tx08, tx09a, tx09b} and Liu {\it et
al}.~\cite{ltx09}.

In 1998, Hartnell {\it et al.}~\cite{hjvw98} defined the {\it
distance $k$-bondage number} of $G$, denoted by $b_k(G)$, to be the
cardinality of a smallest subset $B$ of edges of $G$ with the
property that $\gamma_k(G-B)>\gamma_k(G)$. From Theorem \ref{thm2.3}
it is clear that if $T$ is a nontrivial tree, then $1\leqslant
b_1(T)\leqslant2$. Hartnell {\it et al.}~\cite{hjvw98} generalized
this result to any integer $k\geqslant 1$.

\begin{thm} \label{thm8.1}\textnormal{(Hartnell {\it et al.}~\cite{hjvw98}, 1998)}
For every nontrivial tree $T$ and positive integer $k$, $1\leqslant
b_k(T)\leqslant2$.
\end{thm}                                                       % Theorem 7.1

Hartnell {\it et al.}~\cite{hjvw98} and Topp and
Vestergaard~\cite{tv00} also characterized the trees having distance
$k$-bondage number $2$. In particular, the class of trees for which
$b_1(T)=2$ are just those which have a unique maximum
$2$-independent set (see Theorem~\ref{thm2.6}).
%This characterization is the basis of a linear algorithm for
%determining the $1$-bondage number of a tree presented by Hartnell {\it et al.}

Since, when $k=1$, the distance $1$-bondage number $b_1(G)$ is the
classical bondage number $b(G)$, Theorem~\ref{thm2.8} gives the
NP-completeness of deciding the distance $k$-bondage number of
general graphs.

\begin{thm}\label{thm8.2}
Given a nonempty undirected graph $G$ and positive integers $k$ and
$b$ with $b\leqslant \varepsilon(G)$, determining wether or not
$b_k(G)\leqslant b$ is NP-complete.
\end{thm}

%\vskip6pt

\noindent{\bf Comments}\ For a vertex $x \in V(G)$, the {\it open
$k$-neighborhood} $N_k(x)$ of $x$ is defined as $N_k(x)=\{y\in
V(G):\ 1\leqslant d_G(x,y)\leq k\}$. The {\it closed
$k$-neighborhood} $N_k[x]$ of $x$ in $G$ is defined as
$N_k(x)\cup\{x\}$. Let
 $$
 \Delta_k(G)={\rm{max}}\{|N_k(x)|:\
{\rm{for\ any}}\ x\in V(G)\}.
 $$
Clearly, $\Delta_1(G)=\Delta(G)$. The {\it $k$-th power} of a graph
$G$ is the graph $G^k$ with vertex set $V(G^k)=V(G)$ and edge set
$E(G^k)=\{xy\,:\,1\leqslant d_G(x,y)\leqslant k\}$. The following
lemma holds directly from the definition of $G^k$.

\vskip6pt{\bf Lemma 2.2}\ (Tian and Xu~\cite{tx08})\
$\Delta(G^k)=\Delta_k(G)$ and $\gamma(G^k)=\gamma_k(G)$ for any
graph $G$ and each $k\geq 1$.

\vskip6pt

A graph $G$ is {\it $k$-distance domination-critical}, or {\it
$\gamma_k$-critical} for short, if $\gamma_k(G-x)<\gamma_k(G)$ for
every vertex $x$ in $G$, proposed by Henning, Oellermann and
Swart~\cite{hos03}
%[Henning, M. A., Oellermann, O. R. and
%Swart, H. C., Distance domination critical graphs. {\it J. Combin.
%Math. Combin. Comput}. 44 (2003), 33-45.]

\vskip6pt{\bf Lemma 2.3}\ (Tian and Xu~\cite{tx08})\quad For each
$k\geq 1$, a graph $G$ is $\gamma_k(G)$-critical if and only if
$G^k$ is $\gamma(G^k)$-critical.

\vskip4pt {\it Proof}\quad This is clear for $k=1$, so we assume
$k\geq 2$ below.

Suppose that $G$ is a $\gamma_k$-critical graph. Let $x$ be any
vertex in $G$. From definition, a $k$-dominating set of $G-x$ is a
dominating set of $(G-x)^k$. Since $(G-x)^k$ is a spanning subgraph
of $G^k-x$, it follows that $G^k$ is $\gamma(G^k)$-critical.

For converse, suppose that $G^k$ is $\gamma(G^k)$-critical. Then
$\gamma(G^k-x)<\gamma(G^k)$ for any vertex $x$ in $G$. Thus,
%there must exist such a $\gamma(G^k)$-set $S$ in $G^k$ that contains $x$
%and that $S\setminus\{x\}$ can dominate all vertices in $\overline
%S$ and $x$ is not adjacent to any vertex in $S$ rather than $x$.
%This implies that $S\setminus\{x\}$ is a dominating set of $G^k-x$.
there must exist a dominating set $S$ of $G^k-x$ such that $S$
contains no vertex $y$ such that $d_G(x,y)\leqslant k$. Therefore,
no edge of $G^k$ joining a vertex of $S$ to a vertex of
$V(G^k)-(S\cup\{x\})$ arises in $G^k$ from a path of length at most
$k$ that contains $x$. It follows that $S$ is a dominating set of
$(G-x)^k$, and hence a $k$-dominating set of $G-x$. This completes
the proof. \hfill\rule{1mm}{2mm}

\vskip6pt

{\color{red}By the above facts, can we generalize the results on the
bondage for $G$ to $G^k$? In particular, do the following
propositions hold?

(a)\ $b(G^k)=b_k(G)$ for any graph $G$ and each $k\geq 1$.

(b)\ $b(G^k)\leqslant\Delta_k(G)$ if $G$ is not
$\gamma_k(G)$-critical.}

\vskip6pt

Let $k$ and $p$ be positive integers. A subset $S$ of $V(G)$ is
defined to be a $(k,p)$-dominating set of $G$ if, for any vertex
$x\in V(G)\setminus S$, $|N_k(x)\cap S|\geq p$. The
$(k,p)$-domination number of $G$, denoted by $\gamma_{k,p}(G)$, is
the minimum cardinality among all $(k,p)$-dominating sets of $G$.
Clearly, for a graph $G$, a $(1,1)$-dominating set is a classical
dominating set,  a $(k,1)$-dominating set is a distance
$k$-dominating set, and a $(1,p)$-dominating set is the
above-mentioned $p$-dominating set. This,
$\gamma_{1,1}(G)=\gamma(G)$, $\gamma_{k,1}(G)=\gamma_k(G)$ and
$\gamma_{1,p}(G)=\gamma_p(G)$.

The concept of $(k,p)$-domination in a graph $G$ is a generalized
domination which combines distance $k$-domination and $p$-domination
in $G$. So the investigation of $(k,p)$-domination of $G$ is more
interesting and has received the attention of many researchers, see
for example, \cite{bhs94, fv05, kmv07,lhx10}.

{\color{red}It is quite natural to propose the concept of bondage
number for $(k,p)$-domination. However, as far as we known, none has
proposed this concept until today. This is a worth-while topic for
us.}

\vskip6pt

\subsection{Fractional Bondage Numbers}

% ÒÔÏÂÄÚÈÝÀ´×Ô Domke and Laskar~\cite{dl97}%(´æÓдËÎÄ)

If $\sigma$ is a function mapping the vertex-set $V$ into some set
of real numbers, then for any subset $S\subseteq V$, let
$\sigma(S)=\sum\limits_{x\in S}\sigma(x)$. Also let
$|\sigma|=\sigma(V)$.

A real-value function $\sigma: V\to [0,1]$ is a dominating function
of a graph $G$ if for every $x\in V(G)$, $\sigma(N_G[x])\geqslant
1$. Thus, if $S$ is a dominating set of $G$, $\sigma$ is a function,
where
 $$
 \sigma(x)=\left\{\begin{array}{ll}
 1\ & {\rm if}\ x\in S,\\
 0\ & {\rm otherwise},\end{array}\right.
 $$
then $\sigma$ is a dominating function of $G$.
%A dominating function is
%minimal if for every $x\in V$ with $\sigma(x)>0$, there exists a vertex
%$y\in N_G[x]$ such that $\sigma(N_G[y])=1$.
The fractional domination
number of $G$, denoted by $\gamma^f(G)$, is defined as follows.
 $$
 \gamma^f(G)=\min\{|\sigma|:\ \sigma\ {\rm is\ a \ domination\ fuction\ of}\ G \}.
 $$

The $\gamma^f$-bondage number of $G$, denoted by $b^f(G)$, is
defined as the minimum cardinality of a subset $B\subseteq E$ whose
removal results in $\gamma^f(G-B)>\gamma^f(G)$.

Hedetniemi {\it et al.}\footnote{S. M. Hedetniemi, S. T. Hedetniemi
and T. V. Wimer, Linear time resource allocation algorithms for
trees. Tech. Report URI-014, Dept. Math. Sci. Clemson University,
1987.} were the first to study fractional domination although
Farber\footnote{M. Farber, Domination, independent domination and
duality in strongly chordal graphs. Discrete Applied Mathematics, 7
(1986), 115-130.} introduced the idea indirectly. The concept of the
fractional domination number was proposed by Domke and
Laskar~\cite{dl97}, in 1997. The fractional domination numbers for
some ordinary graphs are determined.

\vskip6pt\begin{prop}\label{prop8.29}\textnormal{(Domke {\it et
al.}~\cite{dhl98}, 1998; Domke and Laskar~\cite{dl97}, 1997)}

(a)\ If $G$ is a $k$-regular graph with order $n$, then
$\gamma^f(G)=n/{k+1}$;

(b)\ $\gamma^f(C_n)=n/3$, $\gamma^f(P_n)=\lceil n/3\rceil$,
$\gamma^f(K_n)=1$;

(c)\  $\gamma^f(G)=1$ if and only if $\Delta(G)=n-1$;

(d)\ $\gamma^f(T)=\gamma(T)$ for any tree $T$;

(e)\ $\gamma^f(K_{n,m})=\frac{n(m+1)+m(n+1)}{nm-1}$, where
$\max\{n,m\}>1$.
\end{prop}

The assertions (a)-(d) are due to Domke {\it et al.}~\cite{dhl98}
and the assertion (e) is due to Domke and Laskar~\cite{dl97}.

According to these results, Domke and Laskar~\cite{dl97} determined
the the fractional domination number for these graphs.

\begin{thm}\label{thm8.29}
\ $b^f(K_n)=\left\lceil\frac{n}{2}\right\rceil;$\ $b^f(K_{n,m})=1$
where $\max\{n,m\}>1$,
  $$
   b^f(C_n)=\left\{\begin{array}{cl}%
   2&\ {\rm if}\ n\equiv 0({\rm mod}\ 3),\\
   1&\ {\rm otherwise};
   \end{array}\right.
  $$
  $$
   b^f(P_n)=\left\{\begin{array}{cl}%
   2&\ {\rm if}\ n\equiv 1({\rm mod}\ 3),\\
   1&\ {\rm otherwise};
   \end{array}\right.
   $$
\end{thm}

It is easy to see that for any tree $T$, $b^f(T)\leqslant 2$. In
fact, since $\gamma^f(T)=\gamma(T)$ for any tree $T$,
$\gamma^f(T')=\gamma(T')$ for any subgraph $T'$ of $T$. It follows
that
 $$
 \begin{array}{rl}
 b^f(T)&=\min\{B\subseteq E(T):\ \gamma^f(T-B)>\gamma^f(T)\}\\
       &=\min\{B\subseteq E(T):\ \gamma(T-B)>\gamma(T)\}\\
       &=b(T)\leqslant 2.
       \end{array}
       $$

%\vskip10pt

{\color{red} There are several bounds on the parameter involving the
degree of vertices. A study of these bounds will appear elsewhere.}

%´ý²é
{\color{blue} Ghoshal {\it et al.}\footnote{J. Ghoshal, R. Laskar,
D. Pillone and C. Wallis, Strong bondage and reinforcement number of
graphs. Preprint.} defined and studied similar parameters involving
the bondage and reinforcement numbers associated with the strong
domination number.}

%\vskip20pt

\newpage

\subsection{Roman Bondage Numbers}

A {\it Roman dominating function} on a graph $G$ is a labeling
$f:V\rightarrow \{0, 1, 2\}$ such that every vertex with label 0 has
at least one  neighbor with label 2. The weight of a Roman
dominating function is the value $f(V(G))=\sum_{u\in V(G)}f(u)$,
denoted by $f(G)$. The minimum weight of a Roman dominating function
on a graph $G$ is called the {\it Roman domination number}, denoted
by $\gamma_{R}(G)$.

A Roman dominating function $f : V\rightarrow \{0, 1, 2\}$ can be
represented by the ordered partition $(V_0, V_1, V_2)$ (or
$(V_{0}^{f},V_{1}^{f},V_{2}^{f})$ to refer to $f$) of $V$, where
$V_i=\{v\in V\mid f(v) = i\}$. In this representation, its weight is
$\omega(f)=|V_1|+2|V_2|$. It is clear that $V_1^f\cup V_2^f$ is a
dominating set of $G$, called {\it the Roman dominating set},
denoted by $D^f_{\rm R}=(V_1,V_2)$. Since $V_1^f\cup V^f_2$ is a
dominating set when $f$ is a Roman dominating function, and since
placing weight 2 at the vertices of a dominating set yields a Roman
dominating function, in \cite{cdhh04}, it was observed that
\begin{equation}\label{e9.4.1}
\gamma(G)\leqslant \gamma_{R}(G)\leqslant 2\gamma(G).
\end{equation}
A graph $G$ is called to be {\it Roman} if $\gamma_{\rm
R}(G)=2\gamma(G)$.

The definition of the Roman dominating function was given implicitly
by Stewart \cite{s99} and ReVelle and Rosing \cite{rr00}. Roman
dominating numbers have been studied. In particular, Bahremandpour
{\it et al.} showed~\cite{bhsx12} that the problem determining the
Roman domination number is NP-complete even for bipartite graphs.

Let  $G$ be a graph with maximum degree at least two. The {\em Roman
bondage number} $b_{R}(G)$ of $G$ is the minimum cardinality of all
sets $E'\subseteq E$ for which $\gamma_{R}(G-E')>\gamma_{R}(G)$.
Since in the study of Roman bondage number the assumption
$\Delta(G)\geqslant 2$ is necessary, we always assume that when we
discuss $b_R(G)$, all graphs involved satisfy $\Delta(G)\geqslant
2$. The Roman bondage number $b_R(G)$ was introduced by Jafari Rad
and Volkmann in \cite{rv11a}.
%, and has been further studied for example in \cite{AQ,
%DKSV, DSV, EP, HX, RV2}.

Recently, Bahremandpour {\it et al.} have showed~\cite{bhsx12} that
the problem determining the Roman bondage number is NP-hard  even
for bipartite graphs.

\begin{center}
\begin{minipage}{130mm}
\setlength{\baselineskip}{24pt}

\vskip6pt\noindent {\bf Roman bondage number problem:}

\noindent {\bf Instance:}\ {\it A nonempty bipartite graph $G$ and a
positive integer $k$.}

\noindent {\bf Question:}\ {\it Is $b_{\rm R}(G)\leqslant k$?}

\end{minipage}
\end{center}

\begin{thm}\label{thm9.4.1}{\rm (Bahremandpour {\it et al.}~\cite{bhsx12}, 2012)}\
The Roman bondage number problem is NP-hard even for bipartite
graphs.
\end{thm}

The exact value of $b_R(G)$ is known only for a few family of graphs
including the complete graphs, cycles and paths.

\begin{thm}\label{thm9.4.2}
{\rm (Jafari Rad and Volkmann~\cite{rv11a}, 2011)} For a path $P_n$
and a cycle $C_n$ of order $n$,
$$
b_{\rm R}(P_n)= \left\{ \begin{array}{ll}
2, & {\rm if}\ n\equiv 2\, ({\rm mod}\, 3);\\
1, & {\rm otherwise}.
 \end{array} \right.
 $$
$$
b_{\rm R}(C_n)= \left\{ \begin{array}{ll}
3, & {\rm if}\ n=2\, ({\rm mod}\, 3);\\
2, & {\rm otherwise}.
 \end{array} \right.
 $$
For a complete graph $K_n$ ($n\geq 3$), $b_{\rm
R}(K_n)=\lceil\frac{n}{2}\rceil$.
\end{thm}

\begin{lem}\label{lem9.4.3} {\rm (Cockayne et al. \cite{cdhh04})}\
If $G$ is a graph of order $n$ and contains vertices of degree
$n-1$, then $\gamma_{\rm R}(G)=2$.
\end{lem}

Using Lemma~\ref{lem9.4.3}, the third conclusion in
Theorem~\ref{thm9.4.2} can be generalized to more general case,
which is similar to Lemma~\ref{lem3.5.2}.

\begin{prop}\label{prop9.4.4}
Let $G$ be a graph with order $n\geq 3$ and
%$\gamma_{\rm R}(G)=2$,
$t$ be the number of vertices of degree $n-1$ in $G$. If $t\geq 1$
then $b_{\rm R}(G)=\lceil\frac{t}{2}\rceil$.
\end{prop}

\begin{pf}
Let $H$ be a spanning subgraph of $G$ obtained by removing fewer
than $\lceil\frac{t}{2}\rceil$ edges from $G$. Then $H$ contains
vertices of degree $n-1$ and, hence, $\gamma_{\rm
R}(H)=2=\gamma_{\rm R}(G)$ by Lemma~\ref{lem9.4.3}, which implies
$b_{\rm R}(G)\geq \lceil\frac{t}{2}\rceil$.

Since $G$ contains $t$ vertices of degree $n-1$, it contains a
complete subgraph $K_t$ induced by these $t$ vertices. We can remove
$\lceil\frac{t}{2}\rceil$ edges such that no vertices have degree
$n-1$ and, hence, $\gamma_{\rm R}(H)\geq 3>2=\gamma_{\rm R}(G)$
since $n\geq 3$. Thus $b_{\rm R}(G)\leqslant
\lceil\frac{t}{2}\rceil$, whence $b_{\rm
R}(G)=\lceil\frac{t}{2}\rceil$.
\end{pf}

\vskip6pt

For a complete bipartite graph $K_{m,n}$, where $1\leqslant
m\leqslant n$. Ebadi and PushpaLatha~\cite{ep10} determined that
 $$
 b_R(K_{m,n})=\left\{\begin{array}{ll}
 1 & \ {\rm if}\ m=1\ {\rm and}\ n\ne 1;\\
 5 & \ {\rm if}\ m=n=3;\\
 m & \ {\rm otherwise}.
 \end{array}\right.
 $$

The above second conclusion shows $b_R(K_{3,3})=5$. However, for a
complete $t$-partite graph, when $t\geqslant 3$, we have the
following result.

\begin{thm}\label{thm9.4.3}
{\rm (Hu and Xu \cite{hx}, 2011)}\ Let $G = K_{m_1,m_2,\ldots,m_t}$
be a complete $t$-partite graph with $m_1=\ldots=m_i<m_{i+1}\leq
\ldots \leq m_t$, $t\geq 2$ and $n=\sum\limits_{j=1}^{t} m_j$. Then
 $$
 b_{\rm R}(G)=\left\{\begin{array}{ll}
 \lceil\frac{i}{2}\rceil\  & {\rm if}\ m_i=1\  {\rm and}\ n\geq 3;\\
 2 \ & {\rm if}\ m_i=2\ {\rm and}\ i=1;\\
 i\  & {\rm if}\ m_i=2\ {\rm and}\ i\geq 2;\\
 n-1 \ & {\rm if}\ m_i=3\ {\rm and}\ i=t\geq 3;\\
 n-m_t  \ & {\rm if}\ m_i\geq 3\  {\rm and}\ m_t\geq 4.
 \end{array}\right.
 $$
\end{thm}

Consider $K_{3,3,\ldots,3}$ of order $n\ge 9$, which is an
$(n-3)$-regular graph. The above result means that $b_{\rm
R}(K_{3,3,\ldots,3})=n-1$. In the same paper, Hu and Xu further
determined that $b_{\rm R}(G)=n-2$ for any $(n-3)$-regular graph $G$
of order $n\ge 5$ and $G\ne K_{3,3,\ldots,3}$.

\begin{thm}\label{thm9.4.4a} {\rm (Hu and Xu \cite{hx}, 2011)}\
Let $G$ be an {\rm ($n-3$)}-regular graph of order $n \ge 5$  but
$G\ne K_{3,3,\ldots,3}$. Then $b_{\rm R}(G)=n-2$.
\end{thm}

For a tree $T$ with order $n\geqslant 3$, Ebadi and
PushpaLatha~\cite{ep10}, and Jafari Rad and Volkmann~\cite{rv11a},
independently, obtained an upper bound of $b_R(T)$.

\begin{thm}\label{thm9.4.3a} \ $b_R(T)\leqslant
3$ for any tree $T$ with order $n\geqslant 3$.
\end{thm}

\begin{thm}{\rm (Bahremandpour {\it et al.}~\cite{bhsx12}, 2012)}\
$b_{\rm R}(P_2\times P_n)=2$ for $n\geqslant 2$.
\end{thm}

\begin{thm}\label{thm9.4.5}
{\rm (Jafari Rad and Volkmann~\cite{rv11a}, 2011)}\ Let $G$ be a
graph of order $n\geqslant 3$.

(a)\ If G is a graph and $(x,y,z)$ a path of length $2$ in $G$, then
$b_R(G)\leqslant d_G(x)+d_G(y)+d_G(z)-3-|N_G(x)\cap N_G(y)|$.

(b)\ If $G$ is connected, then
$b_R(G)\leqslant\lambda(G)+2\Delta(G)-3$.
\end{thm}

Theorem~\ref{thm9.4.5} (a) implies
$b_R(G)\leqslant\delta(G)+2\Delta(G)-3$. Note that for a planar
graph $G$, $\delta(G)\leqslant 5$, moreover, $\delta(G)\leqslant 3$
if the girth at least 4 and $\delta(G)\leqslant 2$ if the girth at
least 6. These two facts show that $b_R(G)\leqslant 2\Delta(G)+2$
for connected planar graphs $G$. Jafari Rad and
Volkmann~\cite{rv11b} improved this bound.

\begin{thm}\label{thm9.4.6}
{\rm (Jafari Rad and Volkmann~\cite{rv11b}, 2011)}\ Let $G$ be a
connected planar graph of order $n\geqslant 3$ with girth $g(G)$.
Then
 $$
 b_R(G)\leqslant\left\{\begin{array}{ll}
 2\Delta(G);\\
 \Delta(G)+6;\\
 \Delta(G)+5 &\ {\rm if}\ G\ {\rm contains\ no\ vertices\ of\
degree\ five};\\
 \Delta(G)+4 & \ {\rm if}\ g(G)\geqslant 4;\\
 \Delta(G)+3 & \ {\rm if}\ g(G)\geqslant 5;\\
 \Delta(G)+2 & \ {\rm if}\ g(G)\geqslant 6;\\
 \Delta(G)+1 & \ {\rm if}\ g(G)\geqslant 8.
 \end{array}\right.
 $$
\end{thm}

According to Theorem~\ref{thm9.4.3}, $b_R(C_n)=3$ for a cycle $C_n$
of length $n\geqslant 8$ with $n\equiv\,2\,{\rm (mod\, 3)}$, and
therefore the last result in Theorem~\ref{thm9.4.6} is best
possible, at least for $\Delta=2$.

Combining the fact that every planar graph $G$ with minimum degree 5
contains an edge $xy$ with $d_G(x)=5$ and $d_G(y)\in\{5, 6\}$ with
Theorem~\ref{thm9.4.5} (a), Akbari, Khatirinejad and
Qajar~\cite{akq12} obtained the following result.

\begin{thm}\label{thm9.4.7}
{\rm (Akbari, Khatirinejad and Qajar~\cite{akq12}, 2012)}\
$b_R(G)\leqslant 15$ for every planar graph $G$.
\end{thm}

It remains open to show whether the bound in Theorem~\ref{thm9.4.7}
is sharp or not. Though finding a planar graph $G$ with $b_R(G)=15$
seems to be difficult, Akbari, Khatirinejad and Qajar~\cite{akq12}
constructed an infinite family of planar graphs with Roman bondage
number equal to 7 by proving the following result.

\begin{figure}[ht]
\psset{unit=.8}
\begin{center}
\begin{pspicture}(-4,-3.5)(4,3.5)

\psellipse(0,0)(2.5,1.6) \cnode*(-2.49,0){2pt}{1}
\cnode*(-2,.93){2pt}{2} \cnode*(0,1.57){2pt}{3}
\cnode*(2,.93){2pt}{4} \cnode*(2.49,0){2pt}{5}
\cnode*(2,-.93){2pt}{6} \cnode*(0,-1.57){2pt}{7}
\cnode*(-2,-.93){2pt}{8}

\cnode*(-3.2,0.7){2pt}{11} \cnode*(-3.9,.7){2pt}{12}
\cnode*(-3.2,-.7){2pt}{13} \cnode*(-3.9,-.7){2pt}{14} \ncline{1}{11}
\ncline{11}{12} \ncline{1}{13}  \ncline{13}{14}

\cnode*(-2.7,1.){2pt}{21} \cnode*(-3.3,1.6){2pt}{22}
\cnode*(-2,1.7){2pt}{23} \cnode*(-2.6,2.3){2pt}{24} \ncline{2}{21}
\ncline{21}{22} \ncline{2}{23}  \ncline{23}{24}

\cnode*(-2.7,-1){2pt}{81} \cnode*(-3.3,-1.6){2pt}{82}
\cnode*(-2,-1.7){2pt}{83} \cnode*(-2.6,-2.3){2pt}{84} \ncline{8}{81}
\ncline{81}{82} \ncline{8}{83}  \ncline{83}{84}

\cnode*(3.2,0.7){2pt}{51} \cnode*(3.9,.7){2pt}{52}
\cnode*(3.2,-.7){2pt}{53} \cnode*(3.9,-.7){2pt}{54} \ncline{5}{51}
\ncline{51}{52} \ncline{5}{53}  \ncline{53}{54}

\cnode*(2.7,1){2pt}{41} \cnode*(3.3,1.6){2pt}{42}
\cnode*(2,1.7){2pt}{43} \cnode*(2.6,2.3){2pt}{44} \ncline{4}{41}
\ncline{41}{42} \ncline{4}{43}  \ncline{43}{44}

\cnode*(2.7,-1){2pt}{61} \cnode*(3.3,-1.6){2pt}{62}
\cnode*(2,-1.7){2pt}{63} \cnode*(2.6,-2.3){2pt}{64} \ncline{6}{61}
\ncline{61}{62} \ncline{6}{63}  \ncline{63}{64}

\cnode*(.6,2.4){2pt}{31} \cnode*(.6,3.1){2pt}{32}
\cnode*(-.6,2.4){2pt}{33} \cnode*(-.6,3.1){2pt}{34} \ncline{3}{31}
\ncline{31}{32} \ncline{3}{33}  \ncline{33}{34}

\cnode*(.6,-2.4){2pt}{71} \cnode*(.6,-3.1){2pt}{72}
\cnode*(-.6,-2.4){2pt}{73} \cnode*(-.6,-3.1){2pt}{74} \ncline{7}{71}
\ncline{71}{72} \ncline{7}{73}  \ncline{73}{74}

\rput(0,0){$G$} \cnode*(1.7,2.3){1pt}{a} \cnode*(1.5,2.45){1pt}{b}
\cnode*(1.3,2.6){1pt}{c} \cnode*(-1.7,2.3){1pt}{a}
\cnode*(-1.5,2.45){1pt}{b} \cnode*(-1.3,2.6){1pt}{c}
\cnode*(1.7,-2.3){1pt}{a} \cnode*(1.5,-2.45){1pt}{b}
\cnode*(1.3,-2.6){1pt}{c} \cnode*(-1.7,-2.3){1pt}{a}
\cnode*(-1.5,-2.45){1pt}{b} \cnode*(-1.3,-2.6){1pt}{c}

%\psellipticarc(0,0)(2.5,1.8){146.7}{180}  %ÍÖÔ²Ïß
\end{pspicture}
\caption{\label{f14} \footnotesize The graph $\widehat{G}$ is
constructed from $G$.}
\end{center}
\end{figure}
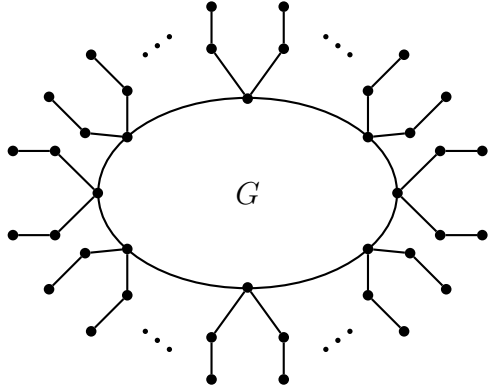

\begin{thm}\label{thm9.4.8}
{\rm (Akbari, Khatirinejad and Qajar~\cite{akq12}, 2012)}\ Let $G$
be a graph of order $n$ and $\widehat{G}$ is the graph of order $5n$
obtained from $G$ by attaching the central vertex of a copy of a
path $P_5$ to each vertex of $G$ (see Figure~\ref{f14}). Then
$\gamma_R(\widehat{G})=4n$ and $b_R(\widehat{G})=\delta(G)+2$.
\end{thm}

By Theorem~\ref{thm9.4.8}, infinitely many planar graphs with Roman
bondage number $7$ by considering any planar graph $G$ with
$\delta(G)=5$ (e.g. the icosahedron graph).

\begin{con}{\rm (Akbari, Khatirinejad and Qajar~\cite{akq12}, 2012)}\
The Roman bondage number of every planar graph is at most 7.
\end{con}

For general bounds, the following observation is directly.

\begin{obs}\label{ob2}
Let $G$ be a graph of order $n$ with maximum degree at least two.
Assume that $H$ is a spanning subgraph of $G$ with $\gamma_{\rm
R}(H)=\gamma_{\rm R}(G)$. If $K=E(G)-E(H)$, then $b_{\rm
R}(H)\leqslant b_{\rm R}(G)\leqslant b_{\rm R}(H)+|K|$.
\end{obs}

\begin{thm}\label{thm9.4.11}{\rm (Bahremandpour {\it et al.}~\cite{bhsx12}, 2012)}\
Let $G$ be a connected graph of order $n\geqslant 3$ with
$\gamma_{\rm R}(G)=\gamma(G)+1$. Then $$b_{\rm R}(G)\leqslant \min
\{b(G),n_{\Delta}\},$$ where $n_{\Delta}$ is the number of vertices
with maximum degree $\Delta$ in $G$.
\end{thm}

\begin{obs}
If $\gamma(G)=\gamma_R(G)$ for a graph $G$, then $b_R(G)\leqslant
b(G)$. \end{obs}

\begin{pf}
Let $F$ be a minimum edge-set of $G$ for which
$\gamma(G-F)>\gamma(G)$, ie, $|F|=b(G)$. By (\ref{e9.4.1}),
$\gamma(G)\leqslant \gamma_R(G)$. By assumption,
$\gamma_R(G)=\gamma(G)<\gamma(G-F)\leqslant \gamma_R(G-F)$. Hence,
$b_R(G)\leqslant b(G)$.
\end{pf}

\begin{thm}\label{thm9.4.12}{\rm (Bahremandpour {\it et al.}~\cite{bhsx12}, 2012)}\
For every Roman graph $G$,
 $$
 b_{\rm R}(G) \geqslant b(G).
 $$
The bound is sharp for cycles on $n$ vertices where $n\equiv
0\;({\rm mod}\;3)$.
\end{thm}

The strict inequality in Theorem~\ref{thm9.4.12} can hold, for
example, $b(C_{3k+2})=2<3=b_{\rm R}(C_{3k+2})$ by
Theorem~\ref{thm9.4.3}.

A graph $G$ is called to be {\it vertex Roman domination-critical}
if $\gamma_{\rm R}(G-x)<\gamma_{\rm R}(G)$ for every vertex $x$ in
$G$. If $G$ has no isolated vertices, then $\gamma_{\rm
R}(G)\leqslant 2\gamma(G) \leqslant 2\beta(G)$. If $\gamma_{\rm
R}(G)=2\beta(G)$, then $\gamma_{\rm R}(G)=2\gamma(G)$ and hence $G$
is a Roman graph. In~\cite{v94}, Volkmann gave a lot of graphs with
$\gamma(G)=\beta(G)$.

The following result is similar to Theorem~\ref{thm4.2}.

\begin{thm}\label{thm9.4.13}{\rm (Bahremandpour {\it et al.}~\cite{bhsx12}, 2012)}\
Let $G$ be a graph. If $\gamma_{\rm R}(G)=2\beta(G)$, then

(a) $b_{\rm R}(G)\geqslant \delta(G)$;

(b) $b_{\rm R}(G)\geqslant \delta(G)+1$ if $G$ is a vertex Roman
domination-critical graph.
\end{thm}

Dehgardi, Sheikholeslami and Volkmann \cite{dsv11} posed the
following problem: If $G$ is a connected graph of order $n\geqslant
4$ with $\gamma_R(G)\geqslant 3$, then
\begin{equation}\label{e9.4.2}
b_{R}(G)\leqslant (\gamma_R(G)-2)\Delta(G).
\end{equation}

Theorem~\ref{thm9.4.5} (a) shows that the inequality (\ref{e9.4.2})
holds if $\gamma_R(G)\geqslant 5$. Thus the bound in (\ref{e9.4.2})
is of interest only when $\gamma_R(G)$ is 3 or 4.

% {\rm (Cockayne et al. \cite{cd04})}\
%If $G$ is a graph of order $n$ and contains vertices of degree
%$n-1$, then $\gamma_{\rm R}(G)=2$.

\begin{lem}\label{lem9.4.14}{\rm (Bahremandpour {\it et al.}~\cite{bhsx12}, 2012)}\
Let $G$ be a nonempty graph of order $n\geqslant 3$, then
$\gamma_{\rm R}(G)=3$ if and only if $\Delta(G)=n-2$.

\end{lem}

The following result shows that (\ref{e9.4.2}) holds for all graphs
$G$ of order $n\geqslant 4$ with $\gamma_R(G)=3,4$, which improves
Theorem~\ref{thm9.4.5} (a).

\begin{thm}\label{tmm9.4.15}{\rm (Bahremandpour {\it et al.}~\cite{bhsx12}, 2012)}\
If $G$ is a connected graph of order $n\geqslant 4$, then
 $$
 b_{R}(G)\le\left\{\begin{array}{ll} \Delta(G)=n-2 & \ {\rm if}\ \gamma_R(G)=3;\\
 \Delta(G)+\delta(G)-1 & \ {\rm if}\ \gamma_R(G)=4
 \end{array}\right.
 $$
with the first equality if and only if $G\cong C_4$.
\end{thm}

Dehgardi et al. \cite{dsv11} proved that  for any connected graph
$G$ of order $n\geqslant 3$, $b_{\rm R}(G)\leqslant n-1$  and posed
the following problems.

{\it Prove or disprove: For any connected graph $G$ of order $n\ge
3$, $b_{\rm R}(G)=n-1$ if and only if $G\cong K_3$.}

{\it Prove or disprove: If $G$ is a connected graph of order $n\ge
3$, then}
 $$
 b_{\rm R}(G)\leqslant n-\gamma_{\rm R}(G)+1.
 $$

Since $\gamma_R(K_{3,3,\ldots,3})=4$, Theorem~\ref{thm9.4.2} shows
that the above two problems are false. Recently Akbari and Qajar
\cite{aq12} proved the following result.

\begin{thm}
If $G$ is a connected graph of order $n\geqslant 3$, then
$$b_{R}(G)\leqslant n-\gamma_R(G)+5.$$
\end{thm}

In \cite{ep10}, Ebadi and PushpaLatha conjectured that
$b_R(G)\leqslant n-1$ for any graph $G$ of order $n\geqslant 3$.
Akbari and Qajar~\cite{aq12} showed that this conjecture is true.

\begin{thm} {\rm (Akbari and Qajar~\cite{aq12}, 2012)}\
$b_R(G)\leqslant n-1$ for any connected  graph $G$ of order
$n\geqslant 3$.
\end{thm}

Theorem~\ref{thm9.4.3} shows this upper is best. We conclude this
section with the following problems.

\begin{prob} Characterize all connected graphs
$G$ of order $n\geqslant 3$ for which $b_R(G)=n-1$.
\end{prob}

\begin{prob} Prove or disprove: If $G$ is a
connected graph of order $n\geqslant 3$, then
$$b_{R}(G)\leqslant n-\gamma_R(G)+3.$$
\end{prob}

%$$ ?????? $$

%\bibitem{akq12)% Saieed Akbari ¡¤ Mahdad Khatirinejad ¡¤ Sahar Qajar,
%S. Akbari, M. Khatirinejad and S. Qajar, A note on the Roman bondage
%number of planar graphs. {\it Graphs and Combinatorics} DOI
%10.1007/s00373-011-1129-8.

%\bibitem{AQ}
%S. Akbari and S. Qajar, A note on Roman bondage number of graphs.
%{\it  Ars Combin.}, (to appear).

%\bibitem {DKSV}
%N. Dehgardi, O. Favaron, B. Kheirfam and S.M. Sheikholeslami,
%Roman fractional bondage number of a graph. {\it J. Combin. Math.
%Combin. Comput.}, (to appear).

%\bibitem {DSV}
%N. Dehgardi, S.M. Sheikholeslami and L. Volkman, On the Roman
%$k$-bondage number of a graph. {\it AKCE Int. J. Graphs Comb.},
%{\bf 8} (2011), 169-180.

%\bibitem{EP}
%K. Ebadi and L. PushpaLatha, Roman bondage and Roman
%reinforcement numbers of a graph. {\it Int. J. Contemp. Math. Sci.},
%{\bf 5} (2010), 1487-1497.

%\bibitem {HX}
%F.-T. Hu and J.-M. Xu, Roman bondage numbers of
%some graphs. (Submitted).

%\bibitem {RV2}
%N.J. Rad and L. Volkmann, On the Roman bondage number of planar  graphs.
%{\it Graphs Combin.}, {\bf 27} (2011), 531-538.

\vskip20pt

\subsection{Remarks and Comments}

There are many variants of domination except mentioned-above ones.

Generally speaking, the concept of $k$-domination has an analog for
all dominations with various restrained conditions.

For example, Hattingh and Henning proposed the concept of connected
$k$-domination [J. H. Hattingh, M. A. Henning, The ratio of the
distance irredundance and domination numbers of a graph, J. Graph
Theory, 18 (1994), 1-9]. A $k$-dominating set $S$ of $G$ is called a
{\it connected $k$-dominating set} of $G$ if the subgraph $G[S]$
induced by $S$ is connected.

Also for example, Henning, Oellermann and Swart proposed the concept
of total $k$-domination [M. A. Henning, O. R. Oellermann, H. C.
Swart, Relations between distance domination parameters, Math.
Pannon., 5 (1) (1994), 69-79]. A dominating set $S$ of $G$ is called
a {\it total $k$-dominating set} of $G$ if every vertex in $G$ is
within distance $k$ from some vertex of $S$ other than itself.

It is quite natural to propose bondage numbers for $k$-dominations
of these types. However, we have not yet seen any research results
on these topics.

\section{Results on Digraphs}

Although domination has been extensively studied in undirected
graphs, it is natural to think of a dominating set as a one-way
relationship between vertices of the graph. Indeed, among the
earliest literature on this subject, J. van Neumman and O.
Morgenstern [Theory of Games and Economic Behavior, Princeton Univ.
Press, Princeton, NJ, 1944] used what is now called domination in
digraphs to find solution (or kernels, which are independent
dominating set) for cooperative $n$-person games. Most likely, the
first formulation of domination by C. Berge [Theor\'ie des Graphes
et ses Applications (Dunod, Paris,1958)] was given in the context of
digraphs and, only some years latter by O. Ore [Theory of Graphs,
Amer. Math. Soc. Colloq. Publ., Vol. 38 (AMC, Providence, RI, 1962)]
for undirected graphs. Despite this history, examination of
domination and its variants in digraphs has been essentially
overlooked (see~\cite{hhs97b} for an overview of the domination
literature). Thus, there are few, if any, such results on domination
for digraphs in the literature.

The bondage number and its related topics for undirected graph have
become one of major areas both in theoretical and applied
researches. However, until recently, Carlson and Develin
\cite{cd06}, Shan and Kang~\cite{sk07}, Huang and Xu~\cite{hx06,
hx07a} studied the bondage number for digraphs, independently. In
this section, we will introduce their results for general digraphs.
Results for some special digraphs such as vertex-transitive digraphs
are introduced in the next section.

\vskip10pt \subsection{Upper Bounds for General Digraphs}

Let $G=(V,E)$ be a digraph without loops and parallel edges. A
digraph $G$ is called to {\it asymmetric} if whenever $(x,y)\in
E(G)$ then $(y,x)\notin E(G)$, and to be {\it symmetric} if
$(x,y)\in E(G)$ implies $(y,x)\in E(G)$. For a vertex $x$ of $V(G)$,
the sets of out-neighbors and in-neighbors of $x$ are, respectively,
defined as $N^+_G(x)=\{y\in V(G):\ (x,y)\in E(G)\}$ and
$N^-_G(x)=\{x\in V(G):\ (x,y)\in E(G)\}$, the out-degree and the
in-degree of $x$ are, respectively, defined as $d^+_G(x)=|N^+_G(x)|$
and $d^-_G(x)=|N^+_G(x)|$. Denote the maximum and the minimum
out-degree\,(\,respectively in-degree) of $G$ by $\Delta^+(G)$ and
$\delta^+(G)$\,(\,respectively $\Delta^-(G)$ and $\delta^-(G)$). The
degree $d_G(x)$ of $x$ is defined as $d^+_G(x)+d^-_G(x)$, the
maximum and the minimum degree of $G$ is denoted by $\Delta(G)$ and
$\delta(G)$, respectively, that is $\Delta(G)=\max\{d_G(x): x\in
V(G)\}$, and $\delta(G)=\min\{d_G(x): x\in V(G)\}$. Note that the
definitions here are different from ones in the next-book on
digraphs.

A subset $S$ of $V(G)$ is called a {\it dominating set} if
$V(G)=S\cup N^+_G(S)$, where $N^+_G(S)$ is the set of out-neighbors
of $S$. Then, just as for undirected graphs, $\gamma(G)$ is the
minimum cardinality of a dominating set, and the {\it bondage
number} $b(G)$ is the smallest cardinality of a set $B$ of edges
such that $\gamma(G-B)>\gamma(G)$ if such a subset $B\subseteq E(G)$
exists. Otherwise, we put $b(G)=\infty$.

Some basic results for undirected graphs stated in Section 3 can be
generalized to digraphs. For example, Theorem~\ref{thm3.1.7} is
generalized by Carlson and Develin \cite{cd06}, Huang and
Xu~\cite{hx06}, Shan and Kang~\cite{sk07}, independently, as
follows.

\begin{thm}\label{thm9.1}
%\textnormal{(Carlson and Develin \cite{cd06}, 2006; Huang and Xu~\cite{hx06},
%2006; Shan and Kang~\cite{sk07}, 2007)}
Let $G$ be a digraph and $(x,y)\in E(G)$.
Then $b(G)\leqslant d_G(y)+d^-_G(x)-|N^-_G(x)\cap N^-_G(y)|$.
\end{thm}                                                     % Lemma 8.1

\begin{cor}\label{cor9.2} %\textnormal{\cite{cd06}}
For a digraph $G$, $b(G)\leqslant\delta^-(G)+\Delta(G)$.
\end{cor}                                                 % Corollary 8.2

Since $\delta^-(G)\leqslant\frac{1}{2}\Delta(G)$, from
Corollary~\ref{cor9.2}, Conjecture~\ref{con3.6.3} is valid for
digraphs.

\begin{cor}\label{cor9.3} %\textnormal{\cite{cd06}}
For a digraph $G$, $b(G)\leqslant\frac{3}{2}\Delta(G)$.
\end{cor}                                                 % Corollary 8.3

In the case of undirected graphs, the bondage number of
$G_n=K_n\times K_n$ achieves this bound. However, it was shown by
Carlson and Develin in \cite{cd06} that if we take the symmetric
digraph $G_n$, we have $\Delta(G_n)=4(n-1)$, $\gamma(G_n)=n$ and
$b(G_n)\leqslant4(n-1)+1=\Delta(G_n)+1$. So this family of digraphs
can not show that the bound in Corollary~\ref{cor9.3} is tight.

Corresponding to Conjecture~\ref{con3.6.1}, which is discredited for
undirected graphs and is valid for digraphs, the same conjecture for
digraphs can be proposed as follows.

\begin{con}\label{con9.4} \textnormal{(Carlson and Develin~\cite{cd06}, 2006)}
 $b(G)\leqslant\Delta(G)+1$ for any digraph $G$.
\end{con}                                                % Conjecture 8.5

If Conjecture~\ref{con9.4} is true, then the following results shows
that this upper bound is tight since $b(K_n\times
K_n)=\Delta(K_n\times K_n)+1$ for a complete digraph $K_n$.

In 2007, Shan and Kang~\cite{sk07} gave some tight upper bounds on
the bondage numbers for some asymmetric digraphs. For example,
$b(T)\leqslant \Delta(T)$ for any asymmetric ditree $T$;
$b(G)\leqslant \Delta(G)$ for any asymmetric digraph $G$ with order
at least $4$ and $\gamma(G)\leqslant 2$. For planar digraphs, they
obtained the following results.

\begin{thm} \textnormal{(Shan and Kang~\cite{sk07}, 2007)}
Let $G$ be a asymmetric planar digraph. Then $b(G)\leqslant
\Delta(G)+2$; and $b(G)\leqslant \Delta(G)+1$ if $\Delta(G)\geqslant
5$ and $d^-_G(x)\geqslant 3$ for every vertex $x$ with
$d_G(x)\geqslant 4$.
\end{thm}                                                   % Theorem 8.6

\vskip10pt \subsection{Results for Some Special Digraphs}

The exact values and bounds of $b(G)$ for some standard digraphs
were determined.

\begin{thm}\label{thm9.6}\textnormal{(Huang and Xu~\cite{hx06}, 2006)}
For a directed cycle $C_n$ and a directed path $P_n$,
 $$
   b(C_n)=\left\{\begin{array}{ll}
3&{\rm if}\ n\ {\rm is\ odd};\\
2&{\rm if}\ n\ {\rm is\ even}.
\end{array}\right.
$$
and %\quad {\rm and}\quad
 $$
 b(P_n)=\left\{\begin{array}{ll}
2&{\rm if}\ n\ {\rm is\ odd};\\
1&{\rm if}\ n\ {\rm is\ even}.
\end{array}\right.
 $$
For the de Bruijn digraph $B(d,n)$ and the Kautz digraph $K(d,n)$,
 $$
 \left\{\begin{array}{ll}
 b(B(d,n))=d
 &{\rm if}\ n\ {\rm is\ odd};\\
 d\leqslant b(B(d,n))\leqslant 2d
 &{\rm if}\ n\ {\rm is\ even};
 \end{array}\right.
 $$
 and %\quad {\rm and}\quad
 $$
 b(K(d,n))=d+1.
 $$
\end{thm}                                             % Corollary 3.4

Like undirected graphs, we can define the total domination number
and the total bondage number. On the total bondage numbers for some
special digraphs, the known results are as follows.

\begin{thm}\textnormal{(Huang and Xu~\cite{hx07a}, 2007)}
For a directed cycle $C_n$ and a directed path $P_n$, $b^t(P_n)$ and
$b^t(C_n)$ all do not exist. For a complete digraph $K_n$,
  $$
 \left\{\begin{array}{ll}
 b^t(K_n)=\infty &\ {\rm if}\ n=1,2;\\
 b^t(K_n)=3&\ {\rm if}\ n=3;\\
 n\leqslant b^t(K_n)\leqslant2n-3&\ {\rm if}\ n\geqslant4.
 \end{array}\right.
 $$
 \end{thm}

The extended de Bruijn digraph $EB(d,n; q_1, \ldots, q_p)$ and the
extended Kautz digraph $EK(d,n; q_1, \ldots, q_p)$ were introduced
by Shibata and Gonda\footnote{Y. Shibata, Y. Gonda, Extension of de
Bruijn graph and Kautz graph. Comput. Math. Appl. 30 (1995) 51-61.}.
If $p=1$, then they are the de Bruijn digraph $B(d,n)$ and the Kautz
digraph $K(d,n)$, respectively. Huang and Xu~\cite{hx07a} determined
their total domination numbers. In particular, their total bondage
numbers for general cases are determined as follows:

\begin{thm}\textnormal{(Huang and Xu~\cite{hx07a}, 2007)}
If $d\geqslant 2$ and $q_i\geqslant 2$ for each $i=1,2,\ldots,p$,
then
 $$
 b^t(EB(d,n; q_1, \ldots, q_p))=d^p-1\ \ {\rm and}\ b^t(EK(d,n; q_1, \ldots, q_p))=d^p.
 $$
In particular, for the de Bruijn digraph $B(d,n)$ and the Kautz
digraph $K(d,n)$,
$$
b^t(B(d,n))=d-1 \quad{\rm and}\quad  b^t(K(d,n))=d.
$$
\end{thm}

Zhang {\it {\it et al.}}~\cite{zlm09} determined the bondage number
in complete $t$-partite digraphs.

\begin{thm}\textnormal{(Zhang {\it {\it et al.}}~\cite{zlm09}, 2009)}
For a complete $t$-partite digraph $K_{n_1,n_2,\ldots,n_t}$, where
$n_1\leqslant n_2\leqslant\cdots\leqslant n_t$,
 $$
b(K_{n_1,n_2,\ldots,n_t})=\left\{
\begin{array}{l}
m\  {\rm if}\ n_m=1\ {\rm and}\ n_{m+1}\geqslant 2\ {\rm for\ some}\
m\
(1\leqslant m<t);\\
4t-3\  {\rm if}\ n_1=n_2=\cdots =n_t=2;\\
\sum\limits_{i=1}^{t-1} n_i+2(t-1)\  {\rm otherwise}.
\end{array}
 \right.
 $$
\end{thm}

\noindent{\bf Comments}\ Since an undirected graph can be thought of
a symmetric digraph, any result for digraphs has an analogy for
undirected graphs in general. In view of this point, studying the
bondage number for digraphs is more significant than for undirected
graphs. Thus, we should further study the bondage number of digraphs
and try to generalized known results on the bondage number and
related variants for undirected graphs to digraphs, prove or
disprove Conjecture~\ref{con9.4}. In particular, determine the exact
values of $b(B(d,n))$ for an even $n$, and $b^t(K_n)$ for
$n\geqslant 4$.

\section{Efficient Dominating Sets}

%We will examine the natural extension of efficient domination to
%oriented graphs, an appropriate setting for efficient domination in
%that it yields more satisfying results than efficient domination in
%graphs.

%In order to obtain upper bound of the bondage number, Huang and
%Xu~\cite{hx08} used efficient dominations.

A dominating set $S$ of a graph $G$ is called to be {\it efficient}
if for every vertex $x$ in $G$, $|N_G[x]\cap S|=1$ if $G$ is a
undirected graph or $|N^-_G[x]\cap S|=1$ if $G$ is a directed graph.

From definition, if $S$ is an efficient dominating set of a graph
$G$, then $S$ is certainly an independent set and every vertex not
in $S$ is adjacent to exactly one vertex in $S$.

It is also clear from definition that a dominating set $S$ is
efficient if and only if $\mathscr N_G[S]=\{N_G[x]:\ x\in S\}$ for
the undirected graph $G$ or $\mathscr N^+_G[S]=\{N^+_G[x]:\ x\in
S\}$ for the digraph $G$ is a partition of $V(G)$, where the induced
subgraph by $N_G[x]$ or $N^+_G[x]$ is an star or an out-star with
the root $x$.

The efficient domination has important applications in many areas,
such as error-correcting codes, and receives much attention in the
late years.

The concept of efficient dominating sets is a measure of the
efficiency of domination in graphs and proposed by Bange {\it et
al.}~\cite{bbs88} in 1988. Unfortunately, as shown in \cite{bbs88},
not every graph has an efficient dominating set and, moreover, it is
an NP-complete problem to determine whether a given graph has an
efficient dominating set. In addition, it has been shown by
Clark~\cite{c93} in 1993 that for a wide range of $p$, almost every
random undirected graph $G\in \mathscr G(\upsilon, p)$ has no
efficient dominating sets. This means that undirected graphs
possessing an efficient dominating set are rare. However, it is easy
to show that every undirected graph has an orientation with an
efficient dominating set (see Bange {\it et al.}~\cite{bbh98}).

In 1993, Barkauskas and Host~\cite{bh93} showed that determining
whether an arbitrary oriented graph has an efficient dominating set
is NP-complete. Even so, the existence of efficient dominating sets
for some graphs has been examined, see, for example, Dejter and
Serra~\cite{ds03} and Lee~\cite{l01} for Cayley graph, Gu, Jia and
Shen~\cite{gjs02} for meshes and tori,  Huang and Xu~\cite{hx08} for
circulant graphs, Harary graphs and tori; Van Wieren, Livingston and
Stout~\cite{vls93} for cube-connected cycles.

In this section, we introduce some results of the bondage number for
some graphs with an efficient dominating set.

\subsection{Results for General Graphs}

In this subsection, we introduce some results on bondage numbers
obtained by applying efficient dominating sets, due to Huang and
Xu~\cite{hx08}. We first state the two following lemmas.

\begin{lem}\label{lem10.1}
Let $G$ be a k-regular graph or digraph of order $\upsilon$. Then
$\gamma(G)\geqslant \upsilon/(k+1)$, with equality if and only if
$G$ has an efficient dominating set. In addition, if $G$ has an
efficient dominating set, then every efficient dominating set is
certainly a $\gamma$-set, and vice versa.
\end{lem}

\begin{pf}
Since $G$ is $k$-regular, then $|N^+[x]|=k+1$ for each $x\in V(G)$.
Hence $\gamma(G)\geqslant\lceil \upsilon(G)/(k+1)\rceil$. It is easy
to observe that the equality holds if and only if there exists a
dominating set $D$ such that $\mathscr N^+[D]$ is a partition of
$V(G)$, equivalently, $D$ is an efficient dominating set.

Now suppose that $G$ has an efficient dominating set, i.e.,
$\gamma(G)=\upsilon(G)/(k+1)$. Then a dominating set $D$ is a
$\gamma$-set if and only if $|D|=\upsilon(G)/(k+1)$. On the other
hand, $D$ is efficient if and only if $|D|=\upsilon(G)/(k+1)$. The
lemma follows.
\end{pf}

\vskip6pt

Let $e$ be an edge and $S$ a dominating set in $G$. We say $e$ {\it
supports} $S$ if $e\in (S,\bar S)$, where $(S,\bar S)=\{(x,y)\in
E(G):\ x\in S, y\in \bar S\}$. Denote by $s(G)$ the minimum number
of edges which support all $\gamma$-sets in $G$.

\begin{lem}\label{lem10.2}
\ For any graph or digraph $G$, $b(G)\geqslant s(G)$, with equality
if $G$ is regular and has an efficient dominating set.
\end{lem}

\begin{pf}
Assume $E'\subseteq E(G)$ with $|E'|<s(G)$. Then $E'$ can not
support all $\gamma$-sets in $G$. Let $D$ be a $\gamma$-set not
supported by $E'$. We prove by contradiction that $D$ is still a
dominating set in $G-E'$.

Suppose to the contrary that there exists some $y\in V(G)\setminus
D\,$ such that $D$ can not dominate it in $G-E'$. Since $D$ is a
dominating set in $G$, there exists a vertex $x\in D$ which
dominates $y$ in $G$. Hence $(x,y)\in E(G)$ supports $D$, which
implies that $(x,y)\notin E'$. It follows that $x$ dominates $y$ in
$G-E'$, a contradiction. Thus, $\gamma(G-E')=\gamma(G)$ for any set
$E'\subseteq E(G)$ with $|E'|<s(G)$, and so $b(G)\geqslant s(G)$.

Now let $G$ be a regular graph with an efficient dominating set, and
$E'$ a set of $s(G)$ edges which supports all $\gamma$-sets. We show
that any $\gamma$-set $D$ is not a dominating set in $H=G-E'$. Since
$E'$ supports $D$, there exists an edge $(x,y)\in E'$ such that
$x\in D$ and $y\notin D$. Hence $y$ is not dominated by $x$ in $H$.
By lemma \ref{lem10.1}, $D$ is efficient, which implies that $D$
dominate $y$ only by $x$. Thus, $D$ can not dominate $y$ in $H$. It
follows that $\gamma(H)>\gamma(G)$, and $b(G)\leqslant|E'|=s(G)$.
The lemma follows.
\end{pf}

\vskip6pt

A graph $G$ is called to be {\it vertex-transitive} if its
automorphism group $Aut(G)$ acts transitively on its vertex-set
$V(G)$. A vertex-transitive graph is regular. Applying
Lemma~\ref{lem10.1} and Lemma~\ref{lem10.2}, Huang and Xu obtained
some results on bondage numbers for vertex-transitive graphs or
digraphs.

\begin{thm}\label{thm10.3}\textnormal{(Huang and Xu~\cite{hx08}, 2008)}
Let $G$ be a vertex-transitive graph or digraph. Then
$$b(G)\geqslant
\left\{\begin{array}{ll}
\lceil \upsilon(G)/2\gamma(G)\rceil &\ {\rm if}\ G\ {\rm is\ undirected};\\
\lceil \upsilon(G)/\gamma(G)\rceil &\ {\rm if}\ G\ {\rm is\
directed}.
\end{array}\right.$$
\end{thm}                                                       % Theorem 9.3

\begin{pf}
Assume $V(G)=\{x_1,\ldots,x_\upsilon\}$. Let $\mathscr{D}_i$ be the
family of all $\gamma$-sets that contain $x_i$ in $G$. We first show
that $|\mathscr{D}_i|=|\mathscr{D}_j|$ for any $i$ and $j$. Since
$G$ is vertex-transitive, there exists an automorphism $\phi$ of $G$
such that $\phi(x_i)=x_j$. Clearly $\phi(D_i)\ne\phi(D'_i)$ for any
distinct $D_i,D'_i\in\mathscr{D}_i$. On the other hand, for any
$D_j\in\mathscr{D}_j$, it holds that
$\phi^{-1}(D_j)\in\mathscr{D}_i$ and $\phi(\phi^{-1}(D_j))=D_j$.
Thus, $\phi$ is a bijection from $\mathscr{D}_i$ to $\mathscr{D}_j$,
and so $|\mathscr{D}_i|=|\mathscr{D}_j|=s$ for any
$i,j\in\{1,2,\ldots,\upsilon\}$.

Note that $\cup_{i=1}^\upsilon\mathscr{D}_i$ contains all
$\gamma$-sets of $G$ and every $\gamma$-set appears $\gamma(G)$
times in it. Hence there are exactly $\upsilon(G)s/\gamma(G)$
$\gamma$-sets in $G$.

If $G$ is undirected, then an edge $x_ix_j$ may only support those
$\gamma$-sets in $\mathscr{D}_i$ and $\mathscr{D}_j$ whose number is
at most $2s$. Hence it needs at least $(\upsilon(G)s/\gamma(G))/2s$
edges to support $\cup_{i=1}^\upsilon\mathscr{D}_i$. It follows from
Lemma \ref{lem10.2} that $b(G)\geqslant s(G)\geqslant\lceil
\upsilon(G)/2\gamma(G)\rceil$.

If $G$ is directed, then an edge $(x_i,x_j)$ only supports those
$\gamma$-sets in $\mathscr{D}_i$. Hence $b(G)\geqslant
s(G)\geqslant\lceil \upsilon(G)s/(\gamma(G)s)\rceil=\lceil
\upsilon(G)/\gamma(G)\rceil.$ The theorem follows.
\end{pf}

\begin{thm}\label{thm10.4a}
If $G$ is a $k$-regular graph, then

$b(G)\leqslant k$ if $G$ is undirected and
$\upsilon(G)\geqslant\gamma(G)(k+1)-k+1$, and

$b(G)\leqslant k+1+\ell$ if $G$ is directed and
$\upsilon(G)\geqslant\gamma(G)(k+1)-\ell$ with $0\leqslant
\ell\leqslant k-1.$
\end{thm}

\begin{pf}
First assume $G$ is undirected. For any $y\in V(G)$ let
$E'=\{(x,y)\in E(G):$ $x\in N_G(y)\}$. Then any minimum dominating
set $D$ in $H=G-E'$ must contain $y$. But $y$ dominates only itself
in $H$. If $|D|=\gamma(G)$, then $D$ dominates at most
 $$
 (|D|-1)(k+1)+1=\gamma(G)(k+1)-k<\upsilon(G)=\upsilon(H)
 $$
vertices in $H$, a contradiction. Hence $\gamma(H)=|D|>\gamma(G)$
and $b(G)\leqslant|E'|=k$.

Now assume $G$ is a digraph. For $y\in V(G)$ let
$N^+(y)=\{w_1,\ldots,w_k\}$ and $E'=\{(x,y):\ x\in
N^-(y)\}\cup\{(y,w_i):\ 1\leqslant i\leqslant \ell+1\}$, where
$0\leqslant \ell\leqslant k-1$. Then any minimum dominating set $D$
in $H=G-E'$ must contain $y$. But $y$ dominates only $k-\ell$
vertices in $H$. If $|D|=\gamma(G)$, then in $H$, $D$ dominates at
most
$$(|D|-1)(k+1)+k-\ell=\gamma(G)(k+1)-\ell-1<\upsilon(G)=\upsilon(H)$$
vertices, a contradiction. Hence $\gamma(H)=|D|>\gamma(G)$ and
$b(G)\leqslant|E'|=k+1+\ell$.
\end{pf}

\vskip6pt

Next we will establish a better upper bound of $b(G)$. To this aim,
we introduce the following terminology, which generalizes the
concept of the edge-covering of a graph $G$. For $V'\subseteq V(G)$
and $E'\subseteq E(G)$, we say $E'$ {\it covers} $\,V'\,$ and call
$E'$ an {\it edge-covering\,} for $V'$ if there exists an edge
$(x,y)\in E'$ for any vertex $x\in V'$. For $y\in V(G)$, let
$\beta'[y]$ be the minimum cardinality over all edge-coverings for
$N^-_G[y]$.
%$N[v]$ or $N^-[v]$ by $\beta'[v]$.

\begin{thm}\label{thm10.4}
\textnormal{(Huang and Xu~\cite{hx08}, 2008)}
If $G$ is a $k$-regular graph with order $\gamma(G)(k+1)$, then
$b(G)\leqslant \beta'[y]$ for any $y\in V(G)$.
\end{thm}

\begin{pf}
For any $y\in V(G)$, let $E'$ be the smallest set of edges that
covers $N^-[y]$. To dominate $y$, any $\gamma$-set $D$ in $G$ must
contain some vertex $x$ in $N^-[y]$. Since $E'$ covers $N^-[y]$,
then $x$ dominates at most $k$ vertices in $H=G-E'$. Hence $D$
dominates at most
 $$
 (|D|-1)(k+1)+k<\gamma(G)(k+1)=\upsilon(G)=\upsilon(H)
 $$
vertices, which implies that $D$ is not a dominating set in $H$.
Thus $\gamma(H)>\gamma(G)$ and $b(G)\leqslant|E'|=\beta'[y]$.
\end{pf}

\vskip6pt

The upper bound of $b(G)$ given in Theorem~\ref{thm10.4} is tight in
view of $b(C_n)$ for a cycle or a directed cycle $C_n$ (see
Theorem~\ref{thm2.1.1} and Theorem~\ref{thm9.6}, respectively).

It is easy to see that for a $k$-regular graph $G$,
$\lceil(k+1)/2\rceil\leqslant \beta'[y]\leqslant k$ when $G$ is
undirected and $\beta'[y]=k+1$ when $G$ is directed. By this fact
and Lemma \ref{lem10.1}, the following theorem is merely a simple
combination of Theorem \ref{thm10.3} and Theorem \ref{thm10.4}.

\begin{thm}\label{thm10.5}
Let $G$ be a vertex-transitive graph of degree $k$. If $G$ has an
efficient dominating set, then
$$
\left\{\begin{array}{ll}%
\lceil\frac{k+1}{2}\rceil\leqslant b(G)
%\leqslant\beta'[v]
\leqslant k
&{\rm if}\ G\ {\rm is\ undirected};\\
b(G)=k+1 &{\rm if}\ G\ {\rm is\ directed}.
\end{array}\right.
$$
\end{thm}

\begin{thm}\label{thm10.6a}
If $G$ is an undirected vertex-transitive cubic graph with order
$4\gamma(G)$ and girth $g(G)\leqslant5$, then $b(G)=2$.
\end{thm}                                                     % Corollary 9.6

\begin{pf}
Since $G$ is a cubic graph of order $\upsilon(G)=4\gamma(G)$, then
by Lemma \ref{lem10.1}, any $\gamma$-set in $G$ is efficient. By
Theorem~\ref{thm10.5}, $2\leqslant b(G)\leqslant3$. Thus, we only
need to show $b(G)\leqslant2$. Let $D$ be an efficient dominating
set in $G$. By the proof of Theorem \ref{thm10.3}, there are
$n(G)s/\gamma(G)=4s$ distinct efficient dominating sets in $G$
provided that a vertex of $G$ belongs to $s$ distinct efficient
dominating sets.

If $g=3$, then there exists a cycle $(u_1,u_2,u_3)$ of length $3$.
Suppose that $v_1$ is the neighbor of $u_1$ such that $v_1$ is not
in the cycle. Then $E'=\{(u_1,v_1),(u_2,u_3)\}$ covers $N[u_1]$. By
Theorem \ref{thm10.4}, $b(G)\leqslant\beta'[u_1]\leqslant|E'|=2$.

If $g=4$ or $5$, there exists a cycle $(u_1,u_2,u_3,u_4)$ or
$(u_1,u_2,u_3,u_4,u_5)$. For any $1\leqslant i<j\leqslant4$, it is
easy to observe that $d(u_i,u_j)\leqslant2$. Note that two distinct
vertices $u,v$ in $D$ satisfy $d(u,v)\geqslant3$, since $N[u]\cap
N[v]=\emptyset$. Hence there exists no efficient dominating set
containing both $u_i$ and $u_j$. Suppose that $\mathscr{D}_i$ is the
family of efficient dominating sets containing $u_i$ for
$i=1,2,3,4$. Then $\mathscr{D}_i\cap\mathscr{D}_j=\emptyset$. It
follows that $E'=\{(u_1,u_2),(u_3,u_4)\}$ supports exactly $4s$
efficient dominating sets, i.e., all sets in
$\cup_{i=1}^4\mathscr{D}_i$. Since there are only $4s$ distinct
efficient dominating sets in $G$, then by Lemma \ref{lem10.2},
$b(G)=s(G)\leqslant|E'|=2$.
\end{pf}

\vskip6pt

\noindent{\bf Remarks} The above proof leads to a byproduct. In the
case of $g=5$ we have $\mathscr{D}_i\cap\mathscr{D}_j=\emptyset$ for
$i=1,2,\ldots,5$. Then $G$ has at least $5s$ efficient dominating
sets. But there are only $\upsilon(G)s/\gamma(G)=4s$ distinct
efficient dominating sets in $G$. This contradiction implies that an
undirected vertex-transitive cubic graph with girth five has no
efficient dominating sets. But a similar argument for $g(G)=3,4$ or
$g(G)\geqslant6$ could not give any contradiction. This is
consistent with the result that $CCC(n)$, a vertex-transitive cubic
graph with girth $n$ if $3\leqslant n\leqslant8$, or girth $8$ if
$n\geqslant9$, has efficient dominating sets for all $n\geqslant3$
except $n=5$ (see Theorem~\ref{thm10.10}).

%Theorem~\ref{thm10.4} holds subject to the condition
%$\upsilon(G)=\gamma(G)(k+1)$. Even if this condition is not
%satisfied, we can obtain another upper bound for $b(G)$, provided
%that $\upsilon(G)$ is close to $\gamma(G)(k+1)$, i.e., $\gamma(G)$
%is close to its lower bound $\upsilon(G)/(k+1)$.

\vskip20pt

\subsection{Results for Cayley Graphs}

In this subsection, we will use Theorem~\ref{thm10.5} to determine
the exact values or approximative values of bondage numbers for some
special vertex-transitive graphs by characterizing the existence of
efficient dominating sets in these graphs.

%Theorem \ref{thm10.4} holds subject to the condition
%$\upsilon(G)=\gamma(G)(k+1)$. Even if this condition is not
%satisfied, we can obtain another upper bound for $b(G)$, provided
%that $\upsilon(G)$ is close to $\gamma(G)(k+1)$, i.e., $\gamma(G)$
%is close to its lower bound $\upsilon(G)/(k+1)$.

%\begin{thm}\label{thm9.7}\textnormal{(Huang and Xu~\cite{hx08}, 2008)}
%If $G$ is a $k$-regular graph, then
%$$b(G)\leqslant\left\{\begin{array}{ll}
%k & \ {\rm if\ G\ is\ undirected\ and}\
%\upsilon(G)\geqslant\gamma(G)(k+1)-k+1;\\
%k+1+\ell & \ {\rm if\ G\ is\ directed\ and}\
%\upsilon(G)\geqslant\gamma(G)(k+1)-\ell,\,0\leqslant \ell\leqslant
%k-1.
%\end{array}\right.$$
%\end{thm}

Let $\Gamma$ be a non-trivial finite group, $S$ be a non-empty
subset of $\Gamma$ without the identity element of $\Gamma$. A
digraph $G$ defined as follows
$$
V(G)=\Gamma;\quad (x,y)\in E(G)\Leftrightarrow x^{-1}y\in S\ {\rm
for\ any}\ x,y\in \Gamma.
$$
is called a {\it Cayley digraph} of the group $\Gamma$ with respect
to $S$, denoted by $C_{\Gamma}(S)$. If $S^{-1}=\{s^{-1}:s\in S\}=S$,
then $C_{\Gamma}(S)$ is symmetric, and is called a Cayley undirected
graph, a Cayley graph for short. Cayley graphs or digraphs are
certainly vertex-transitive.

A circulant graph $G(n;S)$ of order $n$ is a Cayley graph
$C_{Z_n}(S)$, where $Z_n=\{0,1, \ldots, n-1\}$ is the addition group
of order $n$ and $S$ is a nonempty subset of $Z_n$ without the
identity element and, hence, is a vertex-transitive digraph of
degree $|S|$. If $S^{-1}=S$, then $G(n;S)$ is an undirected graph.
If $S=\{1,k\}$, where $2\leqslant k\leqslant n-2$, we write $G(n;
1,k)$ for $G(n;\{1,k\})$ or $G(n;\{\pm 1,\pm k\})$, and call it a
double loop circulant graph.

For directed $G=G(n; 1,k)$, we showed that  $\lceil
n/3\rceil\leqslant \gamma(G)\leqslant\lceil n/2\rceil$ and $G$ has
an efficient dominating set if and only if $3|n$ and $k\equiv 2\,
{\rm (mod\, 3)}$. For directed $G=G(n; 1,k)$, $k\ne n/2$, we showed
that $\lceil n/5\rceil\leqslant \gamma(G)\leqslant\lceil n/3\rceil$
and $G$ has an efficient dominating set if and only if $5|n$ and
$k\equiv \pm 2\, {\rm (mod\, 5)}$. By Theorem~\ref{thm10.5}, we can
obtain the bondage number of a double loop circulant graph if it has
an efficient dominating set.

\begin{thm}\label{thm10.7}
Let $G$ be a double loop circulant graph $G(n;1,k)$. If $G$ is
directed with $3\mid n$ and $k\equiv2\ ({\rm mod}\ 3)$, or $G$ is
undirected with $\,5\mid n$ and $k\equiv\pm2\ ({\rm mod}\ 5)$, then
$b(G)=3$.
\end{thm}

The $m\times n$ torus is the cartesian product $C_m\times C_n$ of
two cycles, and is a Cayley graph $C_{Z_m \times Z_n}(S)$, where
$S=\{(0,1), (1, 0)\}$ for directed cycles and $S=\{(0,¡À1), (¡À1,
0)\}$ for undirected cycles and, hence, is vertex-transitive. Gu
{\it et al}~\cite{gjs02} showed that the undirected torus $C_m\times
C_n$ has an efficient dominating set if and only if both $m$ and $n$
are multiples of 5. We showed that the directed torus $C_m\times
C_n$ has an efficient dominating set if and only if both $m$ and $n$
are multiples of 3. Moreover, we found a necessary condition for a
dominating set containing the vertex $(0, 0)$ in $C_m\times C_n$ to
be efficient, and obtained the following result.

\begin{thm}\label{thm10.8}
Let $G=C_m\times C_n$. If $G$ is undirected and both $m$ and $n$ are
multiples of $5$, or if $G$ is directed and both $m$ and $n$ are
multiples of $3$, then $b(G)=3$.
\end{thm}

The hypercube $Q_n$ is the Cayley graph $C_{\Gamma}(S)$, where
$\Gamma=Z_2\times \ldots\times Z_2 =(Z_2)^n$ and $S =\{100\cdots 0,
010\cdots 0, \ldots, 00 \cdots 01\}$. Lee~\cite{l01} showed that
$Q_n$ has an efficient dominating set if and only if $n=2^m-1$ for a
positive integer $m$. Then we obtain the following result by
Theorem~\ref{thm10.5}.

\begin{thm}\label{thm10.9}
If $n=2^m-1$ for a positive integer $m$, then $2^{m-1}\leqslant
b(Q_n)\leqslant 2^m-1$.
\end{thm}

The $n$-dimensional cube-connected cycle, denoted by $CCC(n)$, is
constructed from the $n$-dimensional hypercube $Q_n$ by replacing
each vertex $x$ in $Q_n$ with an undirected cycle $C_n$ of length
$n$ and linking the $i$th vertex of the $C_n$ to the $i$th neighbor
of $x$. It has been proved that $CCC(n)$ is a Cayley graph and,
hence, is a vertex-transitive graph with degree $3$. Van Wieren {\it
et al.}~\cite{vls93} proved that $CCC(n)$ has an efficient
dominating set if and only if $n\ne 5$. Then we derive the following
result from Theorem~\ref{thm10.5} and Theorem~\ref{thm10.6a}.

\begin{thm}\label{thm10.10}
Let $G=CCC(n)$ be the $n$-dimensional cube-connected cycles with
$n\geqslant 3$ and $n\ne 5$. Then $\gamma(G)=n\,2^{n-2}$ and
$2\leqslant b(G)\leqslant 3$. In addition, $b(CCC(3))=b(CCC(4))=2$.
\end{thm}

{\color{red}{\bf Remarks}\ Whether we can determine the exact value
of $b(CCC(n))$ for $n\geqslant 5$.}

%\section{Others}

%It was shown in~\cite{bhns83} (and later in~\cite{hr92}) that the
%%star is the unique graph with the property that the bondage number
%is 1 and the deletion of any edge results in the domination number
%increasing. We conclude by determining when this very special
%property holds for higher bondage number. Let us call a graph
%uniformly bonded if it has bondage number $b$ and the deletion of
%any $b$ edges results in a graph with increased domination number.

%Theorem 5 (Hartnell and Rall~\cite{hr99}). The only uniformly bonded
%graphs with bondage number 2 are $C_3$ and $P_4$. The unique graph
%with bondage number 3 that is uniformly bonded is the graph $C_4$.
%There are no such graphs for bondage number greater than 3.

\subsection{Comments}

\vskip6pt \noindent{\bf Remarks}\

$FQ_n$ is a Cayley graph and hence vertex-transitive.

The $n$-dimensional augmented cube $AQ_n$ is vertex-symmetric.

The $n$-dimensional star graph $S_{n}$ is vertex- and
edge-transitive

The $n$-dimensional pancake graph $P_n$ is a Cayley graph and,
hence, is vertex transitive.

A bubble-sort graph $B_n$ is a Cayley graph on the symmetric group
on $\{1, 2, \cdots, n\}$ with the set of transpositions
$\{(1,2),(2,3),\cdots, (n-1,n)\}$ as the generating set.

the $(n,k)$-star graph $S_{n,k}$ is a generalization of $S_n$. It
has been shown that $S_{n,k}$ is an $(n-1)$-regular
$(n-1)$-connected vertex-transitive graph.

The $(n,k)$-arrangement graph $A_{n,k}$ is a regular graph of degree
$k(n-k)$ with $\frac {n\,!}{(n-k)\,!}$ vertices and diameter
$\lfloor\frac 32 k\rfloor$. $A_{n,1}$ is isomorphic to a complete
graph $K_n$ and $A_{n,n-1}$ is isomorphic to a star graph $S_n$.
Moreover, $A_{n,k}$ is vertex-transitive and edge-transitive.

An $n$-dimensional alternating group graph $AG_n$ is a Cayley graph
and, hence, is vertex-transitive.

\end{document}